\documentclass[12pt]{article}
\vspace{-1mm}
\usepackage{amssymb}
\usepackage{pstricks}
\newcommand{\R}{\mathbb{R}}
\newcommand{\K}{\mathbb{K}}
\newcommand{\Z}{\mathbb{Z}}

\newcommand{\Q}{\mathbb{Q}}

\newcommand{\N}{\mathbb{N}}

\newcommand{\T}{\mathbb{T}}
\newcommand{\sphere}{\mathbb{S}}
\def\build#1_#2^#3{\mathrel{
\mathop{\kern 0pt#1}\limits_{#2}^{#3}}}
\def\llbracket{[\hspace{-.10em} [ }
\def\rrbracket{ ] \hspace{-.10em}]}
\def\tri{\,\triangle\,}
\def\cq{$\hfill \square$}

\def\ind{{\bf 1}_}
\def\al{\alpha}
\def\ve{\varepsilon}
\def\ga{\gamma}
\def\lm{\lambda}

\def\l{{\cal L}}

\def\m{{\cal M}}
\def\tree{{\cal T}_{\ov{\bf e}}}

\def\g{{\cal G}}

\def\u{{\cal U}}
\def\j{{\cal J}}

\def\t{{\cal T}}

\def\r{{\cal R}}

\def\W{{\cal W}}
\def\w{{\rm w}} 
\def\vf{{\varphi}}
\def\ga{{\gamma}}
\def\mm{{\bf m}}
\def\pp{{\bf p}}
\def\eg{{\bf e}}
\def\PP{\hbox{\bf P}}

\def\beq{\begin{equation}}
\def\eeq{\end{equation}}
\def\ba{\begin{eqnarray*}}
\def\ea{\end{eqnarray*}}
\def\ov{\overline}
\def\wh{\widehat}
\def\wt{\widetilde}

\def\la{\longrightarrow}

\def\noi{\noindent}
\def\proof{\vskip 3mm \noindent{\bf Proof:}\hskip10pt}
\def\rem{\noindent{\bf Remark.} }
\def\rems{\noindent{\bf Remarks.} }

\newtheorem{theorem}{Theorem}[section]
\newtheorem{lemma}[theorem]{Lemma}
\newtheorem{proposition}[theorem]{Proposition}
\newtheorem{corollary}[theorem]{Corollary}

\begin{document}

\title{ \bf GEODESICS IN LARGE PLANAR MAPS\\
AND IN THE BROWNIAN MAP}
\author{
Jean-Fran\c cois {\sc Le Gall}\\
{\small Universit\'e Paris-Sud and Institut universitaire de France}}
\date{}
\maketitle

\begin{abstract}
We study geodesics in the random metric space called the Brownian
map, which appears as the scaling limit of large planar maps.
In particular, we completely describe geodesics starting
from the distinguished point called the root, and we characterize the
set $S$ of all points that are connected to the root by more than one
geodesic. The set $S$ is dense in the Brownian map and homeomorphic
to a non-compact real tree. Furthermore, for every 
$x$ in $S$, the number of distinct geodesics from $x$ to the root
is equal to the number of connected components of $S\backslash\{x\}$. 
In particular, points of the Brownian map can be connected
to the root by at most three distinct geodesics. Our results have 
applications to the behavior of geodesics in large planar maps.
\end{abstract}

\section{Introduction}

The main goal of the present work is to study geodesics in the
random metric space called the Brownian map, which appears as
the scaling limit of several classes of discrete planar maps. 
We prove in particular that a typical point of the Brownian map
is connected by a unique geodesic to the distinguished 
point called the root. We are also able to give an explicit description of the set
of those points that are connected to the root by at least
two distinct geodesics. In particular, we obtain that this set is dense in the Brownian map
and is homeomorphic
to a real tree. Moreover we show that countably many points
are connected to the root by three distinct geodesics, but no
point is connected to the root by four distinct geodesics. Because of
the invariance of the distribution of the Brownian map under
uniform re-rooting, similar results hold if we replace the root by
a point randomly chosen according to the volume measure on the Brownian map.
Although the Brownian map is a singular metric space, there are
striking analogies between our analysis and the classical results
known in the setting of Riemannian manifolds.

Our main results have direct applications to 
uniqueness properties of geodesics in discrete planar maps. One indeed conjectures that
the Brownian map is the universal scaling limit of  
discrete random planar maps, in a way similar to the appearance of 
Brownian motion as the scaling limit of all discrete random paths
satisfying mild integrability conditions. If this conjecture is correct,
the present work will provide information about the behavior of 
geodesics in large discrete random maps in a very general
setting. The preceding analogy with Brownian motion also suggests
that the Brownian map should provide the ``right model'' of 
a Brownian random surface. 

To motivate our main results, let us start by describing some typical
applications to discrete models.
Recall that a planar map is a proper embedding of a 
connected graph in the two-dimensional sphere $\sphere^2$. Loops and multiple edges
are a priori allowed. The faces of the map are the connected components of the complement
of the union of edges. A planar map is rooted 
if it has a distinguished oriented edge called the root edge,
whose  origin is called the root vertex. Two rooted planar maps
are said to be equivalent if the second one is the image of the first one
under an orientation-preserving
homeomorphism of the sphere, which also preserves the root edges.
Two equivalent planar maps are identified.
Given an integer $p\geq 2$, a
$2p$-angulation is a planar map where each face has degree $2p$, that is $2p$ adjacent
edges (one should count edge sides, so that if an edge lies entirely inside a face
it is counted twice). A $2p$-angulation is bipartite, so that it cannot have loops, but
it may have multiple edges. We denote by $\m^p_n$ the set of all
rooted $2p$-angulations with $n$ faces. Thanks to the preceding
identification, the set $\m^p_n$ is finite.

Let $M$ be a planar map and let $V(M)$ denote the vertex set of $M$.
A path in $M$ with length $k$ is  a finite sequence $a_0,a_1,\ldots,a_k$
in $V(M)$ such that $a_i$ and $a_{i-1}$ are connected by an edge of the map,
for every $i\in\{1,\ldots,k\}$. The graph distance $d_{gr}(a,a')$
beween two vertices $a$ and $a'$ is the minimal $k$ such that there exists a
path $\ga=(a_0,a_1,\ldots,a_k)$ with $a_0=a$ and $a_k=a'$. A path
$\ga=(a_0,a_1,\ldots,a_k)$ is called a discrete geodesic (from $a_0$ to $a_k$)
if $k=d_{gr}(a_0,a_k)$. If $\ga=(a_0,\ldots,a_k)$ and $\ga'=(a'_0,\ldots,a'_{k'})$
are two paths (possibly with different lengths), the distance between 
$\ga$ and $\ga'$ is defined by
$$d(\ga,\ga')=\max\{d_{gr}(a_{i\wedge k},a'_{i\wedge k'}):i\geq 0\}.$$

Throughout the present work, we fix an integer $p\geq 2$, and we
consider a random $2p$-angulation $M_n$, which is uniformly distributed over the set $\m^p_n$.
We denote the vertex set of $M_n$ by $\mm_n=V(M_n)$  and 
the root vertex of $M_n$ by $\partial_n\in\mm_n$. Note that $\#(\mm_n)=(p-1)n+2$ by Euler's formula.
For every $a,a'\in \mm_n$, we denote by ${\rm Geo}_n(a,a')$ the set of all 
discrete geodesics from $a$ to $a'$ in the map $M_n$.

\begin{proposition}
\label{uniquenessdiscrete}
For every $\delta>0$,
$$\frac{1}{n}\,\#\{a\in\mm_n: \exists \ga,\ga'\in{\rm Geo}_n(\partial_n,a),\,d(\ga,\ga')\geq \delta n^{1/4}\}
\build{\la}_{n\to\infty}^{} 0$$
in probability.
\end{proposition}

Recall that if $R(M_n)=\max\{d_{gr}(\partial_n,a):a\in\mm_n\}$
is the radius of the map $M_n$, it is known that
$n^{-1/4}R(M_n)$ converges in distribution to a positive random variable
(\cite{CS} if $p=2$ and \cite{MaMi},\cite{We} in the general case).
Thus typical distances between vertices of $M_n$ are of order $n^{1/4}$.
Proposition \ref{uniquenessdiscrete} means that when $n$ is large,
for a typical vertex $a$ of $M_n$, there is essentially a unique discrete geodesic
from the root vertex to $a$, up to deviations that are small in comparison
with the diameter of the map.

One can get a stronger version of Proposition \ref{uniquenessdiscrete} by considering approximate
geodesics. Fix a nonnegative function $\theta(n), n\in\N$ such that $\theta(n)=o(1)$
as $n\to\infty$. An approximate geodesic from $a$ to $a'$ is a path from 
$a$ to $a'$ whose length is less than $(1+\theta(n))\,d_{gr}(a,a')$. Denote the
set of all approximate geodesics from $a$ to $a'$ in $M_n$ by $\ov{\rm Geo}_n(a,a')$.
Then  the statement of Proposition \ref{uniquenessdiscrete} still holds if
${\rm Geo}_n(\partial_n,a)$ is replaced by $\ov{\rm Geo}_n(\partial_n,a)$. 

Proposition \ref{uniquenessdiscrete} is concerned with discrete geodesics 
from the root vertex to a typical point of $\mm_n$.
What happens now if we consider exceptional points\,? To answer this question,
fix $\delta >0$, and for every $a,a'\in\mm_n$, let ${\rm Mult}_{n,\delta}(a,a')$ be the 
maximal integer $k$ such that there exist $k$ paths $\ga_1,\ldots,\ga_k\in {\rm Geo}_n(a,a')$
such that $d(\ga_i,\ga_j)\geq \delta n^{1/4}$ if $i\not= j$. Define analogously
$\ov{\rm Mult}_{n,\delta}(a,a')$ by replacing discrete geodesics with
approximate geodesics.

\begin{proposition}
\label{multiplediscrete}
We have, for every $\delta>0$,
$$\lim_{n\to\infty} P(\exists a\in \mm_n: \ov{\rm Mult}_{n,\delta}(\partial_n,a)\geq 4) =0.$$
\end{proposition}

In other words, when $n$ is large, there cannot be more than $3$
``macroscopically different'' discrete geodesics connecting a point of $\mm_n$
to the root vertex. Can the maximal number $3$ be attained~? The next
proposition provides an answer to this question. 

\begin{proposition}
\label{maximaldiscrete}
We have
$$\lim_{\delta\to0}\Big(\liminf_{n\to\infty} 
P(\exists a\in \mm_n: {\rm Mult}_{n,\delta}(\partial_n,a)=3) \Big) = 1.$$
\end{proposition}

It should be emphasized that the root vertex $\partial_n$ plays no special role in the preceding
propositions. These statements remain valid if $\partial_n$
is replaced by a vertex chosen uniformly at random in $\mm_n$.

The proof of Propositions \ref{uniquenessdiscrete}, \ref{multiplediscrete} and \ref{maximaldiscrete}
depends on analogous results concerning the Brownian map.
We will now explain how the Brownian map can be obtained as the continuous limit of rescaled
uniform $2p$-angulations \cite{IM}. We first observe that $(\mm_n,d_{gr})$ can be viewed as
a random compact metric space, or more precisely as a random variable
taking values in the space $\K$ of isometry classes of compact metric spaces.
We equip $\K$ with the Gromov-Hausdorff distance $d_{GH}$ (see \cite{BBI}, \cite{Gro}
and Section 8 below), and the metric space $(\K,d_{GH})$ is then Polish.
It is proved in \cite{IM} that we can find a sequence $(n_k)_{k\geq 1}$ of values of $n$ converging
to $+\infty$ and then construct the random maps $M_{n_k}$ in such a way that along the sequence 
$(n_k)$ we have the almost sure convergence 
\begin{equation}
\label{convGH}
(\mm_n,\kappa_pn^{-1/4}d_{gr})\build{\la}_{n\to\infty}^{\rm(a.s.)} (\mm_\infty, D)
\end{equation}
in the Gromov-Hausdorff sense (here $\kappa_p$ is a positive constant depending on $p$).
The limiting random metric space $(\mm_\infty, D)$ is called the Brownian map. At this point
we need to comment on the necessity of taking a subsequence in order to get 
the convergence (\ref{convGH}). As will be explained below, the space $\mm_\infty$
is obtained as a quotient space of the well-known random real tree called the CRT,
for an equivalence relation which is explicitly defined in terms of Brownian
labels assigned to the vertices of the CRT. However the metric $D$, which
induces the quotient topology on $\mm_\infty$, has not been completely
characterized, and it is conceivable that different subsequences
might give rise to different limiting metrics $D$ in (\ref{convGH}).  Still one conjectures that the space $(\mm_\infty,D)$
does not depend on the chosen subsequence, nor on the integer $p$, and that 
a convergence analogous to (\ref{convGH}) holds for more general classes
of random planar maps. In the present work, we abusively talk about the Brownian 
map, but it should be understood that we deal with one of the possible 
limits in (\ref{convGH}) (our terminology is thus different from \cite{MaMo}
where the Brownian map refers to the same space $\mm_\infty$, but with a specific choice of
the distance, which may or may not coincide with $D$). Despite the lack of uniqueness, the topological or even metric
properties of $(\mm_\infty, D)$ can be investigated in detail and yield interesting
consequences for large planar maps. For instance
the non-existence of small ``bottlenecks'' in large $2p$-angulations was derived 
in \cite{LGP} as a consequence of the convergence (\ref{convGH}) and the fact
that $(\mm_\infty, D)$ is homeomorphic to $\sphere^2$. Our study of geodesics 
in the Brownian map is motivated in part by the same strategy.

Let us now give a more precise description of the space $\mm_\infty$. We use the notation
$(\t_\eg,d_\eg)$ for the random rooted real tree known as the Continuum Random Tree or CRT,
which was introduced and studied by Aldous \cite{Al1},\cite{Al3}. The notation
$\t_\eg$ is justified by the fact that $\t_\eg$ can be defined as the tree coded by a 
normalized Brownian excursion $\eg=(\eg_t)_{0\leq t\leq 1}$. Precisely, $\t_\eg
=[0,1]\,/\!\sim_{\eg}$ is  the quotient space of the interval $[0,1]$ for the equivalence relation $\sim_\eg$
such that $s\sim_\eg t$ if and only if $\eg_s=\eg_t=\min_{[s\wedge t,s\vee t]}
\eg_r$. The distance $d_\eg(a,b)$ is defined for  $a,b\in\t_\eg$ by $d_\eg(a,b)=\eg_s+\eg_t-2\min_{[s\wedge t,s\vee t]}\eg_r$,
where $s$, respectively $t$, is an arbitrary representative of $a$, resp. $b$, in $[0,1]$
(see subsections 2.3 and
2.4 below). The root $\rho_{\eg}$ of $\t_\eg$ is the equivalence class
of $0$ in the quotient  $[0,1]\,/\!\sim_{\eg}$, and the uniform probability measure 
on $\t_\eg$ is the image of Lebesgue measure on $[0,1]$ under the
canonical projection. 

We then introduce 
Brownian labels assigned to the vertices of the CRT: Conditionally
given the tree $\t_\eg$, $(Z_a)_{a\in\t_\eg}$ is the centered Gaussian process whose distribution
is characterized by the properties $Z_{\rho_{\eg}}=0$ and $E[(Z_a-Z_b)^2]=d_\eg(a,b)$
for every $a,b\in\t_{\eg}$. 

We are in fact interested in the
pair $(\t_{\ov\eg},(\ov Z_a)_{a\in\t_{\ov\eg}})$, which may be defined as the pair $(\t_\eg,(Z_a)_{a\in\t_\eg})$
conditioned on the event $\{Z_a\geq 0\;,\;\forall a\in\t_{\eg}\}$. Some care is needed here
since the latter event has zero probability. The paper \cite{LGW} provides a simple 
explicit construction of the pair $(\t_{\ov\eg},(\ov Z_a)_{a\in\t_{\ov\eg}})$: If $a_*$
denotes the (almost surely unique) vertex 
of $\t_\eg$ such that $Z_{a_*}=\min\{Z_a:a\in\t_{\eg}\}$, $\t_{\ov\eg}$
coincides with the tree $\t_\eg$ re-rooted at $a_*$,
and the new labels $\ov Z_a$ are obtained by setting $\ov Z_a=Z_a - Z_{a_*}$,
so that the label of the new root is still $0$. We can again write $\t_{\ov\eg}=[0,1]\,/\!\sim_{\ov\eg}\,$,
where the equivalence relation $\sim_{\ov\eg}$ is defined as above from
a random continuous function $\ov\eg$ which has a simple expression in terms of $\eg$
(see subsection 2.4 below) and the root $\rho=\rho_{\ov\eg}$ of $\t_{\ov\eg}$ is the equivalence
class of $0$ in $[0,1]\,/\!\sim_{\ov\eg}$.
If $a,b\in\t_{\ov\eg}$, we denote by $[a,b]$ the subset of $\t_{\ov\eg}$ which is
the image under the projection $[0,1]\la[0,1]\,/\!\sim_{\ov\eg}$ of the smallest interval
$[s,t]$ such that $s$, respectively $t$, is a representative of $a$, resp. $b$, in $[0,1]$,
and $s\leq t$. If there exists no such interval $[s,t]$, we take $[a,b]=\t_{\ov\eg}$
by convention. Informally, $[a,b]$ corresponds to those points of $\t_{\ov\eg}$ that are visited
when going from $a$ to $b$ ``around the tree'' in ``clockwise order'' and avoiding the root. We define an equivalence
relation on $\t_{\ov\eg}$ by setting, for every $a,b\in\t_{\ov\eg}$,
$$a\approx b\quad\hbox{if and only if}\quad \ov Z_a=\ov Z_b=\min_{c\in[a,b]} \ov Z_c
\quad\hbox{or}\quad \ov Z_a=\ov Z_b=\min_{c\in[b,a]} \ov Z_c.$$

The space $\mm_\infty$ is then defined as the quotient space $\t_{\ov\eg}\,/\!\approx$. 
We denote the canonical projection from $\t_{\ov\eg}$ onto $\mm_\infty$ 
by $\Pi$.
The metric $D$ induces the quotient topology on $\mm_\infty$, and satisfies the
bound
\begin{equation}
\label{boundD}
D(\Pi(a),\Pi(b))\leq \ov Z_a + \ov Z_b -2\min_{c\in[a,b]} \ov Z_c
\end{equation}
for every $a,b\in\t_{\ov\eg}$.
We use the same notation $\rho$ for the root of $\t_{\ov\eg}$,
and for its equivalence class in $\t_{\ov\eg}\,/\!\approx$,
which is a distinguished point of $\mm_\infty$.

If $(E,\delta)$ is a compact metric space and $x,y\in E$,
a geodesic or shortest path from $x$ to $y$ is a continuous path $\ga=(\ga(t))_{0\leq t\leq \delta(x,y)}$ such that
$\ga(0)=x$, $\ga(\delta(x,y))=y$ and $\delta(\ga(s),\ga(t))=|t-s|$ for every $s,t\in[0,\delta(x,y)]$.
The compact metric space $(E,\delta)$ is then called a geodesic space if any two
points in $E$ are connected by (at least) one geodesic. 
From (\ref{convGH}) and the fact that Gromov-Hausdorff limits of geodesic spaces are geodesic
spaces (see \cite{BBI}, Theorem 7.5.1), one gets that $(\mm_\infty,D)$ is almost surely a geodesic
space. Our main goal is to determine the geodesics between the root $\rho$
and an arbitrary point of $\mm_\infty$. 

Before stating our main result, we still need to introduce one more notation.
We define the skeleton ${\rm Sk}(\t_{\ov\eg})$ by saying that a point $a$ of $\t_{\ov\eg}$
belongs to ${\rm Sk}(\t_{\ov\eg})$ if and only if $\t_{\ov\eg}\backslash \{a\}$
is not connected (informally, ${\rm Sk}(\t_{\ov\eg})$ is obtained by removing the leaves of the 
tree $\t_{\ov\eg}$). 
It is proved in Proposition \ref{skeleton} below that, with probability one, 
the restriction of the projection $\Pi$ to ${\rm Sk}(\t_{\ov\eg})$ is a homeomorphism,
and the Hausdorff dimension of $\Pi({\rm Sk}(\t_{\ov\eg}))$ is less than or equal to $2$.
We write ${\rm Skel}_\infty=\Pi({\rm Sk}(\t_{\ov\eg}))$ to simplify notation.
Since the Hausdorff dimension of $\mm_\infty$
is equal to $4$ almost surely \cite{IM}, the set ${\rm Skel}_\infty$
is a relatively
small subset of $\mm_\infty$. The set
${\rm Skel}_\infty$ is dense in $\mm_\infty$ and
from the previous observations it is homeomorphic to
a non-compact real tree. Moreover,  for every $x\in {\rm Skel}_\infty$,
the set ${\rm Skel}_\infty\backslash\{x\}$ is not connected.

The following theorem, which summarizes our main 
contributions, provides a nice geometric interpretation
of the set ${\rm Skel}_\infty$ (see Theorems \ref{main1} 
and \ref{main2} below for more precise statements).

\begin{theorem}
\label{numbergeodesics}
The following properties hold almost surely. For every
$x\in\mm_\infty\backslash {\rm Skel}_\infty$, there is a unique geodesic from
$\rho$ to $x$. On the other hand, for every $x\in {\rm Skel}_\infty$, the number 
of distinct geodesics from $\rho$ to $x$ is equal to the number of connected
components of ${\rm Skel}_\infty\backslash\{x\}$. In particular, the maximal number
of distinct geodesics from $\rho$ to a point of $\mm_\infty$ is equal to $3$. 
\end{theorem}

\rems (i) The invariance of the distribution of the Brownian map under uniform re-rooting
(see Section 8 below) shows that results analogous to Theorem \ref{numbergeodesics}
hold if one replaces the root by a point $z$ distributed uniformly over $\mm_\infty$.
Here the word ``uniformly'' refers to the volume measure $\lambda$ on $\mm_\infty$,
which is the image of the uniform probability measure on $\t_{\ov\eg}$ under the projection $\Pi$. 

\noindent (ii) The last assertion of the theorem easily follows from the previous ones. Indeed,
we already noticed that the (unrooted) tree $\t_{\ov\eg}$ is isometric to the
CRT $\t_\eg$, and standard properties of linear Brownian motion imply that the
maximal number of connected components of ${\rm Sk}(\t_\eg)\backslash\{a\}$, when $a$
varies over $\t_\eg$, is equal to $3$. Furthermore there are countably many
points $a$ for which this number is $3$.

\smallskip
The construction of the Brownian map $(\mm_\infty,D)$ as a quotient space
of the random tree $\t_{\ov\eg}$ may appear artificial, even though it is
a continuous counterpart of the bijections relating labelled trees to 
discrete planar maps (see in particular \cite{BDG}). Theorem \ref{numbergeodesics}
shows that the skeleton of $\t_{\ov\eg}$, or 
rather its homeomorphic image under the canonical projection $\Pi$, 
has an intrinsic geometric meaning: It exactly corresponds to the cut locus 
of $\mm_\infty$ relative to the root $\rho$, provided we define this cut locus as the set of all points 
that are connected to $\rho$ by at least two distinct geodesics. Note that this definition
of the cut locus is slightly different from the one that appears in Riemannian geometry (see e.g. \cite{GHL}), since the latter does not
make sense in our singular setting.

Remarkably enough, the assertions of Theorem \ref{numbergeodesics}
parallel the known results in the setting of differential geometry, which go back
to Poincar\'e \cite{Poin}. Myers \cite{My}
proved that for a complete analytic two-dimensional manifold which is homeomorphic
to the sphere, the cut locus associated
with a given point $A$ is a topological tree, and the number of distinct geodesics joining $A$ to a 
point $M$  of the cut locus is equal to the number of connected 
components of the complement of $\{M\}$ in the cut locus
(see Hebda \cite{Heb} and Itoh \cite{It} and the references therein for more recent related results
under $C^\infty$ or $C^2$-regularity assumptions). On the other hand, 
Shiohama and Tanaka \cite{ST} give examples showing that 
in the (non-differentiable) setting of Alexandrov spaces with curvature
bounded from below, the cut locus may have a fractal structure.

We are able to give explicit formulas for all geodesics from the root 
to an arbitrary point of $\mm_\infty$. Indeed all these geodesics are obtained
as simple geodesics, which had already been considered in \cite{MaMo}. The 
main difficulty in the proof of Theorem \ref{numbergeodesics} is to verify that 
there are no other geodesics from the root. To this end, we define for every
point $x\in\mm_\infty$ a minimal and a maximal geodesic from the root
to $x$. Loosely speaking, these are defined in such a way that any other
geodesic from $\rho$ to $x$ will lie ``between'' the minimal and the maximal
one: See Section 4 for more exact statements. Then one needs to check that
for a given point $x\in\mm_\infty\backslash{\rm Skel}_\infty$,
the minimal and the maximal geodesic from $\rho$ to $x$ coincide,
and therefore also coincide with the simple geodesic from $\rho$ to $x$.
For this purpose, the key step is to prove that a minimal (or maximal)
geodesic cannot visit a point of ${\rm Skel}_\infty$, except 
possibly at its endpoint. The technical estimates of Sections 5 and 6 below are
devoted to the proof of this property. Some of these estimates are of
independent interest. In particular Corollary \ref{unibound} gives the following uniform estimate
on the volume of balls in $(\mm_\infty,D)$: For every $\beta\in\,]0,1[$, there exists
a (random) constant $K_\beta$ such that the volume of any ball of radius 
$r$ in $\mm_\infty$ is bounded above by $K_\beta\,r^{4-\beta}$, for every $r>0$.
In the multifractal formalism, this means that the multifractal spectrum
of the volume measure $\lambda$ on $\mm_\infty$ is degenerate.

A rather surprising consequence of our results is the fact that any two geodesics starting
from the root, or from a typical point of the Brownian map, must coincide over a small interval. More precisely, for every $\ve >0$, there exists a (random) constant $\alpha>0$
such that, if $\ga$ and $\ga'$ are two geodesics starting from the root and with
length greater than $\ve$, we have $\ga(t)=\ga'(t)$ for every $t\in[0,\alpha]$
 (Corollary \ref{coincid}).

Let us comment on recent results related to the present work. The idea of
studying continuous limits of discrete planar maps appeared in the
pioneering paper of Chassaing and Schaeffer \cite{CS}. The problem of 
establishing a convergence of the type (\ref{convGH})
was raised by Schramm \cite{Sch} in the setting of triangulations.
Marckert and Mokkadem \cite{MaMo} considered the case of quadrangulations
and proved a weak form of the convergence (\ref{convGH}). 
See also \cite{MaMi} for a generalization to Boltzmann distributions
on bipartite planar maps. As an important
ingredient of our proofs, we use a bijection between bipartite planar maps
and certain labelled trees called mobiles, which is due to Bouttier, Di Francesco and Guitter
\cite{BDG}. The recent work of Miermont \cite{Mi1}
and Miermont and Weill \cite{MW} strongly suggests that 
a convergence analogous to (\ref{convGH}) should hold for planar maps
that are not bipartite, such as triangulations. In a recent paper \cite{Mi2},
Miermont has obtained, independently of the present work, certain uniqueness results for geodesics in
continuous limits of discrete quadrangulations, in a setting which is however
different from ours. See also the recent work of Bouttier and Guitter \cite{BG} for a detailed
discussion of the number of geodesics connecting two given points 
in a large planar map, and of exceptional pairs of points that can be linked by
``truly'' distinct geodesics. In a different but related direction, the papers of
Angel and Schramm \cite{AS}, Angel \cite{An} and Chassaing and Durhuus \cite{CD}
study various properties of random infinite planar maps that are uniformly distributed
in some sense.

To complete this presentation, let us mention that planar maps are
important objects in several areas of mathematics and physics. They have 
been studied extensively in combinatorics since the pioneering work
of Tutte \cite{Tu}. Planar maps, or maps on more general surfaces, have significant
geometric and algebraic applications: See the book of Lando and Zvonkin \cite{LS}. The interest of
planar maps in theoretical physics first arose from their connections with 
matrix integrals \cite{tH},\cite{BIPZ}. More recently, planar maps have served
as models of random (discrete) surfaces in the theory of two-dimensional
quantum gravity: See in particular the book of Ambj\o rn, Durhuus, 
and Jonsson \cite{ADJ}. Bouttier's recent thesis \cite{Bo} provides a nice account of the
relations between the statistical physics of random surfaces and the combinatorics of
planar maps.

The paper is organized as follows. Section 2 contains a detailed presentation
of the basic objects which are of interest in this work. In particular we discuss the
Bouttier-Di Francesco-Guitter bijection between bipartite planar maps and
the labelled trees called mobiles \cite{BDG}, which plays an important role
in our arguments.
Such bijections between maps and trees were discovered by
Cori and Vauquelin \cite{CV} and then studied in particular by
Schaeffer \cite{Sc}. Theorem \ref{mainIM}
restates the main result of \cite{IM} in a form convenient for our applications.
Section 3 is devoted to some preliminary results. In Section 4, we recall 
the definition of simple geodesics, and we introduce minimal and
maximal geodesics. The main result of this section is Proposition \ref{minigeo},
which shows that the so-called minimal geodesic is indeed a geodesic. 
Section 5 proves the key technical estimate (Lemma \ref{keylemma}).
Loosely speaking, this lemma bounds the probability 
that the range of a minimal geodesic intersects (the image under the canonical projection of) an
interval containing the right end of an excursion interval of $\ov\eg$ with length
greater than some fixed $\delta >0$. The estimate of Lemma \ref{keylemma} is then
used in Section 6 to prove Proposition \ref{leftinc}, which shows that the range of a 
typical minimal geodesic does not
intersect ${\rm Skel}_\infty$. As another ingredient of the proof of
Proposition \ref{leftinc}, we use our uniform estimates for the volume of small
balls in $\mm_\infty$. Section 7 contains the proof
of our main results Theorems \ref{main1} and \ref{main2}, from which Theorem \ref{numbergeodesics}
readily follows. 
Section 8 discusses the re-rooting invariance property of the Brownian map.
Finally Section 9 provides applications to large planar maps, and gives the proof of Propositions \ref{uniquenessdiscrete},
\ref{multiplediscrete} and \ref{maximaldiscrete}. The proof of two technical discrete lemmas is 
presented in the Appendix.

Let us conclude with a comment about  sets of zero probability. As usual in a
random setting, many of the results that are presented in this work hold almost surely,
that is outside a set of zero probability. In a few instances, such as Lemma \ref{leaf}, 
this set of zero probability depends on the choice of a parameter $U\in[0,1]$,
which corresponds to fixing a point of $\mm_\infty$. In such cases, we will always
make this dependence clear: Compare Propositions \ref{keyprop} and \ref{keyprop2}
below for instance.

\medskip
\noi{\it Acknowledgments}. I am indebted to Fr\'ed\'eric Paulin for useful references and
for his comments on a preliminary version of this work, and to Gr\'egory Miermont
for several stimulating conversations. I also thank J\'er\'emie Bouttier and Emmanuel Guitter
for keeping me informed about their recent work on geodesics 
in random planar maps.

\section{The Brownian map}

\subsection{Real trees}

As was already mentioned in the introduction, the Brownian map
is defined as a quotient space of a random real tree.
We start by discussing real trees in a deterministic setting.
A metric space $(\t,d)$ is a real tree if the following two
properties hold for every $a,b\in \t$.
\begin{description}
\item{\rm(a)} There is a unique
isometric map
$f_{a,b}$ from $[0,d(a,b)]$ into $\t$ such
that $f_{a,b}(0)=a$ and $f_{a,b}(
d(a,b))=b$.
\item{\rm(b)} If $q$ is a continuous injective map from $[0,1]$ into
$\t$, such that $q(0)=a$ and $q(1)=b$, we have
$$q([0,1])=f_{a,b}([0,d(a,b)]).$$
\end{description}
\noindent A rooted real tree is a real tree $(\t,d)$
with a distinguished vertex $\rho=\rho(\t)$ called the root.

Let us consider a rooted real tree $(\t,d)$ with root $\rho$. To avoid trivialities,
we assume that $\t$ has more than one point.
For $a\in\t$, the number $d(\rho,a)$ is called the generation
of $a$ in the tree $\t$.
For $a,b\in\t$, the range of the mapping $f_{a,b}$ in (a) is denoted by
$\llbracket a,b\rrbracket=\llbracket b,a\rrbracket$ : This is the line segment between $a$
and $b$ in the tree. We will also use the obvious notation 
$\rrbracket a,b\llbracket$, $\llbracket a, b \llbracket$, $\rrbracket a,b\rrbracket$. 
For every $a\in \t$, $\llbracket \rho,a\rrbracket$ is  interpreted as the ancestral
line of vertex $a$. 

More precisely we can define a partial order on the
tree, called the genealogical order, by setting $a\prec b$
if and only if $a\in\llbracket \rho,b
\rrbracket$. If $a\prec b$, $a$ is called an ancestor of $b$, and
$b$ is a descendant of $a$. If $a,b\in\t$, there is a unique $c\in\t$ such that
$\llbracket \rho,a
\rrbracket\cap \llbracket \rho,b
\rrbracket=\llbracket \rho,c
\rrbracket$. We write $c=a\tri b$ and call $c$ the most recent
common ancestor to $a$ and $b$. Note that 
$\llbracket a,b \rrbracket=\llbracket a\tri b,a\rrbracket \cup \llbracket a\tri b,b\rrbracket$.

The {\it multiplicity} of a vertex $a\in\t$ is the number of
connected components of $\t\backslash\{a\}$. In particular, $a$ is
called a leaf if it has multiplicity one. The {\it skeleton} ${\rm Sk}(\t)$ is the set 
of all vertices $a$ of $\t$ which are not leaves. Note that ${\rm Sk}(\t)$ equipped with
the induced metric is itself a (non-compact) real tree. 
A point $a\in{\rm Sk}(\t)$
 is called simple if $\t\backslash \{a\}$ has exactly two connected
 components.  If $\t\backslash \{a\}$ has (at least) three connected
 components, we say that $a$ is a branching point of 
 $\t$. 

If $a\in \t$, the subtree of descendants of $a$ is denoted by
$\t(a)$ and defined by
$$\t(a)=\{b\in\t:a\prec b\}.$$
If $a\not =\rho$, $a$ belongs to ${\rm Sk}(\t)$ if and only if
$\t(a)\not =\{a\}$.


\subsection{Coding compact real trees}

Compact real trees can be coded by
``contour functions''. Let $\sigma>0$ and let 
$g$ be a continuous function from $[0,\sigma]$
into $[0,\infty[$ such that $g(0)=g(\sigma)=0$. To
avoid trivialities, we will also assume that $g$ is not identically zero.
For every $s,t\in[0,\sigma]$, we set
$$m_g(s,t)=\inf_{r\in[s\wedge t,s\vee t]}g(r),$$
and
$$d_g(s,t)=g(s)+g(t)-2m_g(s,t).$$
It is easy to verify that $d_g$ is a pseudo-metric 
on $[0,\sigma]$. As usual, we introduce the equivalence
relation
$s\sim_g t$ if and only if $d_g(s,t)=0$ (or equivalently if and only if $g(s)=g(t)=m_g(s,t)$).
The function $d_g$ induces a distance on the quotient space $\t_g
:=[0,\sigma]\,/\!\sim_g$, and we keep the
notation $d_g$ for this distance. We denote by
$p_g:[0,\sigma]\longrightarrow
\t_g$ the canonical projection. Clearly $p_g$ is continuous (when
$[0,\sigma]$ is equipped with the Euclidean metric and $\t_g$ with the
metric $d_g$), and therefore $\t_g=p_g([0,\sigma])$ is a compact metric space.
Moreover, it is easy to verify that the topology induced by $d_g$ coincides with the
quotient topology on $\t_g$. 

By Theorem 2.1 of \cite{DuLG}, 
the metric space $(\t_g,d_g)$ is a (compact) real tree. 
Furthermore the mapping $g\la \t_g$
is continuous with respect to the Gromov-Hausdorff distance,
if the set of continuous functions $g$ is equipped with the supremum distance.
We will always view $(\t_g,d_g)$ as a rooted real tree with root $\rho_g=p_g(0)=p_g(\sigma)$.
Note that $d_g(\rho_g,a)=g(s)$ if $a=p_g(s)$.

If $s,t\in[0,\sigma]$, the
property $p_g(s)\prec p_g(t)$ holds if and only if $g(s)=m_g(s,t)$.
Suppose that $0\leq s<t\leq \sigma$. Then $\llbracket p_g(s),p_g(t)\rrbracket
\subset p_g([s,t])$, and in particular
$p_g(s)\tri p_g(t)=p_g(r)$ for any $r\in[s,t]$ such that $p_g(r)=m_g(s,t)$. Such simple
remarks will be used without further comment in the forthcoming proofs.

Let $a\in\t_g$, and let $s(a)$, respectively $t(a)$, denote the smallest, resp. largest,
element in $p_g^{-1}(a)$. Then $\t_g(a)=p_g([s(a),t(a)])$. Hence,
if $g$ is not constant on any non-trivial interval, 
a vertex $a\not =\rho_g$ belongs to ${\rm Sk}(\t_g)$
if and only if $s(a)<t(a)$, that is if $p_g^{-1}(a)$ is not a singleton. Moreover,
if $g$ does not vanish on $]0,\sigma[$, then $\rho_g\notin {\rm Sk}(\t_g)$.
The last two properties will hold a.s. for the (random) functions $g$ that are
considered below.

\subsection{The Brownian snake}

Let $g$ be as in the previous subsection, and also assume that $g$
is H\"older continuous with exponent $\delta$
for some $\delta>0$. We first introduce the Brownian snake
driven by the function $g$.

Let
$\W$ be the space of all finite paths in $\R$. Here a finite path is simply 
a continuous mapping $\w:[0,\zeta]\la \R$, where
$\zeta=\zeta_{(\w)}$ is a nonnegative real number called the 
lifetime of $\w$. The set $\W$ is a Polish space when equipped with the
distance
$$d_\W(\w,\w')=|\zeta_{(\w)}-\zeta_{(\w')}|+\sup_{t\geq 0}|\w(t\wedge
\zeta_{(\w)})-\w'(t\wedge\zeta_{(\w')})|.$$
The endpoint (or tip) of the path $\w$ is denoted by $\wh \w=\w(\zeta_{(\w)})$.
For every $x\in\R$, we set $\W_x=\{\w\in\W:\w(0)=x\}$.

The Brownian snake driven by $g$ is the continuous random process 
$(W^g_s)_{0\leq s\leq \sigma}$ taking values in $\W_0$, whose distribution
is characterized by the following properties:
\begin{description}
\item{(a)} For every $s\in[0,\sigma]$, $\zeta_{(W^g_s)}=g(s)$.
\item{(b)} The process $(W^g_s)_{0\leq s\leq \sigma}$ is time-inhomogeneous
Markov, and its transition kernels are specified as
follows. If $0\leq s\leq
s'$,
\begin{description}
\item{$\bullet$} $W^g_{s'}(t)=W^g_{s}(t)$ for every $t\in [0,m_g(s,s')]$, a.s.;
\item{$\bullet$} the random path $(W^g_{s'}(m_g(s,s')+t)-W^g_{s'}(m_g(s,s')))_{0\leq t\leq
g(s')- m_g(s,s')}$ is independent of $W^g_s$ and distributed as a
one-dimensional Brownian motion started at $0$ and stopped at time 
$g(s')- m(s,s')$.
\end{description}
\end{description}
Informally, the value $W^g_s$ of the Brownian snake at time $s$
is a random path with lifetime $g(s)$. When $g(s)$ decreases,
the path is erased from its tip, and when $g(s)$ increases, the path 
is extended by adding ``little pieces'' of Brownian paths at its tip.
The path continuity of the process $(W^g_s)_{0\leq s\leq \sigma}$
(or rather the existence of a modification with continuous sample paths)
easily follows from the fact that $g$ is H\"older continuous.

Property (b) implies that if $s\sim_g s'$ then 
$W^g_s=W^g_{s'}$ a.s., and this property holds simultaneously for
all pairs $(s,s')$ outside a single set of zero probability, by
a continuity argument. Hence we may view 
$W^g$ as indexed by the tree $\t_g=[0,\sigma]\,/\!\sim_g$. 
We write $Z^g_s=\wh W^g_s$ for the endpoint of $W_s$.
According to the preceding remark, we can view $Z^g$
as indexed by the tree $\t_g$: If $a\in\t_g$,
we interpret $Z^g_a$ as the spatial position of the
vertex $a$. Then it is not difficult to verify that,
for every $r\in[0,d_g(\rho_g,a)]$, $W^g_a(r)$
is the spatial position of the ancestor of $a$
at generation $r$.

The process $(Z^g_a)_{a\in\t_g}$ can be viewed as Brownian 
motion indexed by $\t_g$: Indeed it is a centered Gaussian process
such that $Z^g_{\rho_g}=0$ and $E[(Z^g_a-Z^g_b)^2]=d_g(a,b)$
for every $a,b\in\t_g$.

We now randomize the coding function $g$.
Let $\eg=(\eg_t)_{t\in[0,1]}$ be the normalized Brownian excursion, and take 
$g=\eg$ and $\sigma=1$
in the previous discussion.
The random real tree $(\t_{\eg},d_\eg)$ coded by $\eg$ is the so-called CRT, or
Continuum Random Tree. Using the fact that local minima of Brownian motion
are distinct, one easily checks that points of $\t_{\eg}$ can have multiplicity at most $3$
(equivalence classes for $\sim_\eg$ can contain at most three points).

We then consider
the process $(W^\eg_s)_{s\in[0,1]}$ such that conditionally given $\eg$, 
$(W^\eg_s)_{s\in[0,1]}$ is the Brownian snake driven by $\eg$. 
 Notice that for every $s\in[0,1]$,
$W^\eg_s=(W^\eg_s(t),0\leq t\leq \eg_s)$ is now a random path with a
random lifetime $\eg_s$. As previously, we write $Z^\eg_s=\wh W^\eg_s$ for the
endpoint of $W^\eg_s$. We refer to \cite{Zu} for a detailed discussion 
of the Brownian snake driven by a Brownian excursion, and
its connections with nonlinear partial differential equations.

 \subsection{Conditioning the Brownian snake}
 
In view of our applications, it is important 
to consider the process $(W^\eg_s)_{s\in[0,1]}$
conditioned on the event 
$$W^\eg_{s}(t)\geq 0\ \hbox{ for every }s\in[0,1]\hbox{ and }t\in[0,\eg_s].$$
Here some justification is needed for the conditioning,
since the latter event has probability zero. The paper \cite{LGW}
describes several limit procedures that allow one to make sense 
of the previous conditioning. These procedures all lead to
the same limiting process $\ov W$ which can be described as 
follows from the original process $W^\eg$. Set
$$Z_*=\inf_{t\in[0,1]} Z^\eg_t$$
and let $s_*$ be the (almost surely) unique time in $[0,1]$ such 
that $Z^\eg_{s_*}=Z_*$. The fact that the minimum $ Z_*$ is attained at a unique
time (\cite{LGW} Proposition 2.5) entails that the vertex $p_\eg(s_*)$ is a 
leaf of the tree $\t_{\eg}$. For
every $s,t\in[0,1]$, set $s\oplus t=s+t$ if $s+t\leq 1$ and $s\oplus t=s+t-1$ if $s+t>1$. Then, for every
$t\in[0,1]$, we set
\begin{description}
\item{$\bullet$} 
$\displaystyle{\ov\eg_t=\eg_{s_*}+\eg_{s_*\oplus t}-2\,m_{\eg}(s_*,
{s_*\oplus t})}$;
\item{$\bullet$} $\ov Z_t=Z^\eg_{s_*\oplus t} -Z^\eg_{s_*}$.
\end{description}
Note that $\ov Z_0=\ov Z_1=0$ and $\ov Z_t > 0$ for every $t\in\,]0,1[$.
The function 
$\ov\eg$ is continuous on $[0,1]$, positive on $]0,1[$, and such that $\ov\eg_0=\ov\eg_1=0$.
Hence the tree $\t_{\ov\eg}$ is well defined, and this tree is
isometrically identified with the tree $\t_{\eg}$ re-rooted at the
(minimizing) vertex $p_\eg(s_*)$: See Lemma 2.2 in \cite{DuLG}. 
To simplify notation, we will write $\rho=\rho_{\ov\eg}$ for the root of $\tree$. 

One easily verifies that $s\sim_{\ov \eg} t$ if and only if $s_*\oplus s\sim_{\eg}s_*\oplus t$,
and so $\ov Z_t$ only depends on the equivalence class of $t$
in the tree $\t_{\ov\eg}$. Therefore we may and will often
view $\ov Z$ as indexed by vertices of the tree $\t_{\ov\eg}$.

The conditioned Brownian snake $(\ov W_s)_{s\in[0,1]}$ is now defined 
by the following properties. For every $s\in[0,1]$, $\ov W_s$ is the random element of
$\W_0$ such that:
\begin{description}
\item{\rm(a)} The lifetime of $\ov W_s$ is $\ov\eg_s$.
\item{\rm(b)} We have $\wh{\ov W}_s=\ov Z_s$, and more generally,
for every $r\in[0,\ov\eg_s]$, $\ov W_s(r)=\ov Z_{a_s(r)}$, where $a_s(r)$ is
the ancestor of $p_{\ov\eg}(s)$ at generation $r$ in the tree $\t_{\ov\eg}$.
\end{description}
To interpret this definition, note that $\ov Z$ and $\ov W$
can both be viewed as indexed by the tree $\t_{\ov\eg}$, which is identified with the
tree $\t_{\eg}$ re-rooted at the minimal spatial position. The new spatial positions
$\ov Z_a$ on the re-rooted tree are obtained by shifting the original spatial 
positions $Z^\eg_a$ in such a way that the position of the new root is still zero, and the path $\ov W_a$
just gives the (new) spatial positions along the ancestral line of $a$ in the
re-rooted tree. See \cite{LGW} for more details.

\subsection{The Brownian map}

To simplify notation we write $\sim$ instead of $\sim_{\ov\eg}$ in the remaining part
of this work. Then $\sim$ is a (random) equivalence relation on $[0,1]$ whose graph
is closed.

We now use the process $(\ov Z_t)_{t\in[0,1]}$
of the previous subsection
to define one more (random) equivalence relation on $[0,1]$. 
For every $s,t\in[0,1]$, we write  
$$s\approx t\quad\hbox{if and only if}\quad \ov Z_s=\ov Z_t=\min_{s\wedge t\leq r\leq s\vee t}
\ov Z_r.$$
Notice the obvious similarity with the definition of $\sim$ (in fact the equivalence
relation $\approx$ is nothing but $\sim_{\ov Z}$). From the fact that local minima
of $Z$ are distinct (Lemma 3.1 in \cite{LGP}), one easily obtains that equivalence 
classes of $\approx$ contain one, two or three points at most: See the discussion
in Section 2 of \cite{LGP}.

We say that $t\in\,]0,1]$
is a left-increase time of ${\ov\eg}$ (respectively of $\ov Z$) if there exists $\ve\in\,]0,t]$ such that ${\ov\eg}_r\geq {\ov\eg}_t$ (resp. $\ov Z_r\geq \ov Z_t$) for every
$r\in[t-\ve,t]$. We similarly define the notion of a right-increase time.
For $t\in\,]0,1[$, $p_{\ov\eg}(t)\in {\rm Sk}(\t_{\ov\eg})$ if and only if $t$ is a 
left-increase or a right-increase time of $\ov\eg$.

\begin{lemma}
\label{equirel}
With probability one, any point $t\in\,]0,1[$ which is a right-increase or a left-increase
time of $\ov\eg$ is neither a right-increase nor a left-increase time of $\ov Z$.
Consequently, it is almost surely true that, for every $s,t,r\in\,]0,1[$, 
the properties $s\sim t$ and $s\approx r$
imply that $s=t$ or $s=r$.
\end{lemma}

In other words, if the equivalence class of $s\in\,]0,1[$ for $\sim$ is not a singleton, then
its equivalence class for $\approx$ must be a singleton, and conversely
(we need to exclude the values $s=0$ and $s=1$, since clearly the pair
$\{0,1\}$ is an equivalence class for both $\sim$ and $\approx$). Lemma
\ref{equirel} is proved in \cite{LGP} (Lemma 3.2). This lemma plays a very important
role in what follows. We will systematically discard the negligible set on which
the conclusion of Lemma \ref{equirel} does not hold.

Thanks to Lemma \ref{equirel}, we can define a new equivalence relation $\simeq$
on $[0,1]$ whose graph is the union of the graphs of $\sim$ and $\approx$
respectively: $s\simeq t$ if and only if  $s\sim t$ or $s\approx t$. 

By definition, the Brownian map $\mm_\infty$ is the quotient space $[0,1]/\simeq$, equipped 
with a random distance $D$ which is obtained as a weak limit of the graph
distance on approximating discrete maps. We will explain this in greater detail below,
but we already record certain properties that can be found in
\cite{IM}. It is more convenient to view $D$
as a pseudo-distance on $[0,1]$. Precisely, the random process $(D(s,t))_{s,t\in[0,1]}$
is continuous, takes values in $[0,\infty[$, and satisfies the following:
\begin{description}
\item{(i)} $D(s,t)=D(t,s)$ and $D(r,t)\leq D(r,s)+D(s,t)$ for every $r,s,t\in[0,1]$;
\item{(ii)} $D(s,t)=0$ if and only $s\simeq t$, for every $s,t\in[0,1]$;
\item{(iii)} $D(0,t)=\ov Z_t$ for every $t\in[0,1]$;
\item{(iv)} $\displaystyle{D(s,t)\leq \ov Z_s+\ov Z_t -2\min_{s\wedge t\leq r\leq s\vee t} \ov Z_r}$,
for every $s,t\in[0,1]$.
\end{description}

Thanks to (i) and (ii), $D$ induces a distance on the quotient space $\mm_\infty
=[0,1]\,/\!\simeq$,
which is still denoted by $D$. The (random) metric space
$([0,1]\,/\!\simeq,D)$ appears as a limit of rescaled planar maps,
in the sense of the Gromov-Hausdorff convergence: See subsection 2.7 below. 
The canonical projection from $[0,1]$ onto $\mm_\infty$ will be denoted by
$\pp$.
By definition, the volume measure $\lambda$ on $\mm_\infty$ is the image
of Lebesgue measure on $[0,1]$ under $\pp$.

In agreement with the presentation that was given in the introduction above, 
it is often useful to view the Brownian map as a quotient space of the
random real tree $\t_{\ov\eg}$. Note that the equivalence relation $\approx$
makes sense on $\t_{\ov\eg}$: If $a,b\in\t_{\ov\eg}$, $a\approx b$
if and only if there exist a representative $s$ of $a$ and a representative $t$
of $b$ in $[0,1]$ such that $s\approx t$. Then, $D$ induces a
pseudo-distance on $\t_{\ov\eg}$: $D(a,b)=D(s,t)$ if $s$, respectively $t$, is any
representative of $a$, resp. $b$. The quotient space $\t_{\ov\eg}\,/\!\approx$, equipped
with $D$, is a metric space, which is canonically identified to
$(\mm_\infty,D)$. This slightly different perspective will be
useful as the genealogical structure of $\t_{\ov\eg}$ plays an important
role in our arguments.  Note that the
equivalence class of any $a\in{\rm Sk}(\t_{\ov\eg})$ for $\approx$
is a singleton, by Lemma \ref{equirel}.

As in Section 1,  the canonical projection from $\t_{\ov\eg}$ onto $\mm_\infty$
is denoted by $\Pi$. Note that $\pp=\Pi \circ p_{\ov\eg}$
(as above $p_{\ov\eg}$ denotes the canonical projection
from $[0,1]$ onto $\t_{\ov\eg}$).
Both projections $\pp$ and $\Pi$ are continuous, when $[0,1]$ is
equipped with the usual topology and $\tree$ with the distance
$d_{\ov\eg}$: For the first one, this is a consequence of (iv), and the result for the
second one follows because the topology of
$\t_{\ov\eg}$ is the quotient topology.
 It will 
be important to carefully distinguish elements of $[0,1]$ from their
equivalence classes in $\mm_\infty$ or in $\t_{\ov\eg}$. We will 
typically use the
letters $s,t$ to denote elements of $[0,1]$, 
$a,b$ for elements of $\t_{\ov\eg}$ and $x,y$ for elements
of $\mm_\infty$. The symbol $\rho$ will stand both for the root of
$\t_{\ov\eg}$ and for the corresponding element in $\mm_\infty$,
which is just the equivalence class of $0$ or of $1$.

If $x\in\mm_\infty$,  property (iii) above implies that
$\ov Z_t=D(\rho,x)$ for every $t\in[0,1]$ such that $\pp(t)=x$. So $\ov Z$
can also be viewed as a random function on $\mm_\infty$. We will
write indifferently $\ov Z_x=\ov Z_a=\ov Z_t$ when $x\in\mm_\infty$, $a\in\t_{\ov\eg}$
and $t\in[0,1]$ are such that $\pp(t)=\Pi(a)=x$.

\subsection{Discrete maps and the Bouttier-Di Francesco-Guitter bijection}

Recall that the integer $p\geq 2$ is fixed throughout this work, and that
$\m^p_n$ denotes the set of all
rooted $2p$-angulations with $n$ faces. 
We will first discuss the Bouttier-Di Francesco-Guitter bijection between 
$\m^p_n$ and the set of all $p$-mobiles with $n$ black vertices.

We use the standard formalism for plane trees as found in section 2.1
of \cite{IM}. 
A plane tree $\tau$ is a finite subset of the set
$${\cal U}=\bigcup_{n=0}^\infty \N^n  $$
of all finite sequences of positive integers (including the empty 
sequence $\varnothing$), which satisfies three obvious conditions:
First $\varnothing\in\tau$, then, for every $v=(u_1,\ldots,u_k)\in\tau$
with $k\geq 1$, the sequence $(u_1,\ldots,u_{k-1})$ (the ``parent'' of $v$)
also belongs to $\tau$, and finally for every $v=(u_1,\ldots,u_k)\in\tau$ there
exists an integer $k_v(\tau)\geq 0$ (the ``number of children'' of $v$) such that
the vertex $vj:=(u_1,\ldots,u_k,j)$ belongs to $\tau$ if
and only if $1\leq j\leq k_v(\tau)$. 
The generation of $v=(u_1,\ldots,u_k)$ is denoted by $|v|=k$. 
The notions of an ancestor and a descendant in the tree $\tau$
are defined in an obvious way.

\begin{center}
\begin{picture}(450,140)
\linethickness{1.5pt}
\put(160,0){\circle{12}}
\put(146,-14){$\varnothing$}
\put(155,3){\thinlines\line(-2,1){56}}

\put(165,3){\thinlines\line(2,1){56}}
\put(220,30){\circle*{12}}
\put(228,17){$2$}
\put(245,55){\circle{12}}
\put(223,33){\line(1,1){18}}
\put(195,55){\circle{12}}
\put(180,42){$21$}
\put(217,33){\line(-1,1){18}}
\put(252,42){$22$}

\put(100,30){\circle*{12}}
\put(90,17){$1$}
\put(125,55){\circle{12}}
\put(132,42){$12$}
\put(103,33){\line(1,1){18}}
\put(75,55){\circle{12}}
\put(61,42){$11$}
\put(97,33){\line(-1,1){18}}
\put(75,85){\circle*{12}}
\put(55,72){$111$}
\put(75,61){\thinlines\line(0,1){20}}
\put(100,110){\circle{12}}
\put(78,88){\line(1,1){18}}
\put(50,110){\circle{12}}
\put(23,97){$1111$}
\put(72,88){\line(-1,1){18}}
\put(106,97){$1112$}

\thinlines \put(300,0){\vector(1,0){140}}
\put(300,0){\circle*{3}}
\thinlines \put(300,0){\vector(0,1){130}}
\put(314,56){\circle*{3}}
\thicklines \put(300,0){\line(1,4){14}}
\thicklines \put(314,56){\line(1,4){14}}
\thicklines \put(328,112){\line(1,0){14}}
\put(328,112){\circle*{3}}
\thicklines \put(342,112){\line(1,-4){14}}
\put(342,112){\circle*{3}}
\thicklines \put(356,56){\line(1,0){14}}
\put(356,56){\circle*{3}}
\thicklines \put(370,56){\line(1,-4){14}}
\put(370,56){\circle*{3}}
\thicklines \put(384,0){\line(1,4){14}}
\put(384,0){\circle*{3}}
\thicklines \put(398,56){\line(1,0){14}}
\put(398,56){\circle*{3}}
\thicklines \put(412,56){\line(1,-4){14}}
\put(412,56){\circle*{3}}
\put(426,0){\circle*{3}}

\put(440,-10){$i$}
\put(282,125){$C^{\tau^\circ}_i$}

\put(422,-10){$pn$}

\put(314,0){\line(0,1){4}}
\put(312,-10){$1$}

\put(293,52){$1$}
\put(300,55){\line(1,0){4}}

\end{picture}

\vspace{8mm}

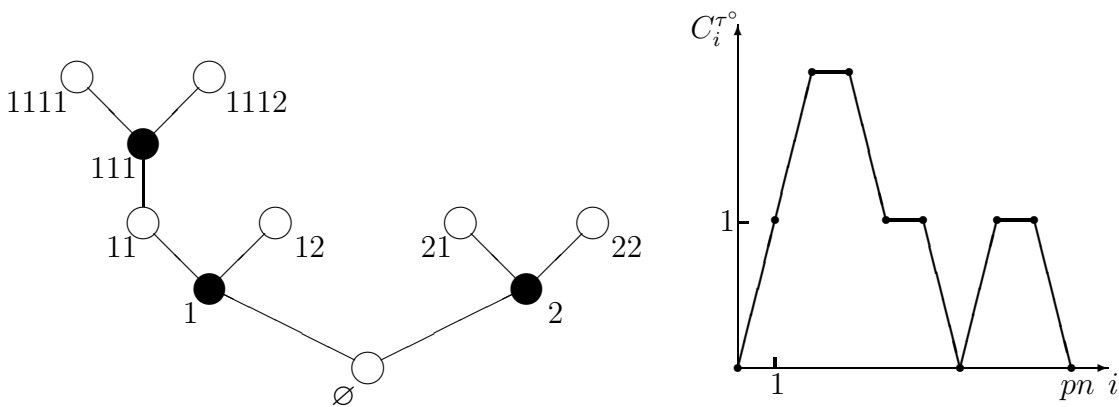
Figure 1. A $3$-tree $\tau$ and the associated 
contour function $C^{\tau^\circ}$ of $\tau^\circ$.

\end{center}

A $p$-tree is a plane tree $\tau$ that satisfies the 
following additional property:
For every $v\in\tau$ such that $|v|$ is odd, $k_v(\tau)=p-1$.

If $\tau$
is a $p$-tree, vertices $v$ of $\tau$ such that $|v|$ is even are called
white vertices, and vertices $v$ of $\tau$ such that $|v|$ is odd are called
black vertices. We denote by $\tau^\circ$ the set of all white vertices
of $\tau$ and by $\tau^\bullet$ the set of all black vertices.
See the left side of Fig.1 for an example of a $3$-tree.

A (rooted) $p$-mobile is a pair $\theta=(\tau,(\ell_v)_{v\in\tau^\circ})$
that consists of a $p$-tree $\tau$ and a collection of integer labels
assigned to the white vertices of
$\tau$, such that the following properties hold:
\begin{description}
\item{(a)} $\ell_\varnothing=1$ and $\ell_v\geq 1$ for each $v\in\tau^\circ$.
\item{(b)} Let $v\in \tau^\bullet$, let $v_{(0)}$ be the parent of $v$ and let
$v_{(j)}=vj$, $1\leq j\leq p-1$, be the children of $v$. Then for every $j\in\{0,1,\ldots,p-1\}$,
$\ell_{v_{(j+1)}}\geq \ell_{v_{(j)}}-1$, where by convention $v_{(p)}=v_{(0)}$.
\end{description}

The left side of Fig.2 gives an example of a $p$-mobile with $p=3$. The numbers appearing 
inside the circles representing white vertices are the labels assigned to
these vertices. Condition (b) above means that if one lists the white vertices
adjacent to a given black
vertex in clockwise order, the labels of these vertices can decrease by at most one
at each step.

\begin{center}
\begin{picture}(400,170)
\linethickness{1.5pt}
\put(140,0){\circle{12}}
\put(137,-4){$1$}
\put(135,3){\line(-4,3){34}}
\put(145,3){\line(4,3){34}}

\put(180,30){\circle*{8}}
\put(205,55){\circle{12}}
\put(183,33){\line(1,1){18}}
\put(155,55){\circle{12}}
\put(152,51){$1$}
\put(177,33){\line(-1,1){18}}
\put(202,51){$2$}

\put(100,30){\circle*{8}}
\put(125,55){\circle{12}}
\put(103,33){\line(1,1){18}}
\put(75,55){\circle{12}}
\put(72,51){$3$}
\put(97,33){\line(-1,1){18}}
\put(122,51){$2$}
\put(75,85){\circle*{8}}
\put(75,61){\thinlines\line(0,1){20}}
\put(100,110){\circle{12}}
\put(78,88){\line(1,1){18}}
\put(50,110){\circle{12}}
\put(47,106){$4$}
\put(72,88){\line(-1,1){18}}
\put(97,106){$3$}
\put(46,114){\line(-1,1){20}}
\put(25,135){\circle*{8}}
\put(104,114){\line(1,1){20}}
\put(125,135){\circle*{8}}
\put(150,160){\circle{12}}
\put(128,138){\line(1,1){18}}
\put(100,160){\circle{12}}
\put(97,156){$2$}
\put(122,138){\line(-1,1){18}}
\put(147,156){$1$}
\put(50,160){\circle{12}}
\put(28,138){\line(1,1){18}}
\put(0,160){\circle{12}}
\put(-3,156){$3$}
\put(22,138){\line(-1,1){18}}
\put(47,156){$2$}

\thinlines \put(250,2){\vector(1,0){180}}
\thinlines \put(250,2){\vector(0,1){160}}

\put(250,35){\circle*{3}}
\put(250,35){\line(1,6){11}}
\put(261,101){\circle*{3}}
\put(261,101){\line(1,3){11}}
\put(272,134){\circle*{3}}
\put(272,134){\line(1,-3){11}}
\put(283,101){\circle*{3}}
\put(283,101){\line(1,-3){11}}
\put(294,68){\circle*{3}}
\put(294,68){\line(1,6){11}}
\put(305,134){\circle*{3}}
\put(305,134){\line(1,-3){11}}
\put(316,101){\circle*{3}}
\put(316,101){\line(1,-3){11}}
\put(327,68){\circle*{3}}
\put(327,68){\line(1,-3){11}}
\put(338,35){\circle*{3}}
\put(338,35){\line(1,6){11}}
\put(349,101){\circle*{3}}
\put(349,101){\line(1,0){11}}
\put(360,101){\circle*{3}}
\put(360,101){\line(1,-3){11}}
\put(371,68){\circle*{3}}
\put(371,68){\line(1,-3){11}}
\put(382,35){\circle*{3}}
\put(382,35){\line(1,0){11}}
\put(393,35){\circle*{3}}
\put(393,35){\line(1,3){11}}
\put(404,68){\circle*{3}}
\put(404,68){\line(1,-3){11}}
\put(415,35){\circle*{3}}

\put(242,32){$1$}
\put(242,65){$2$}
\put(242,98){$3$}
\put(242,131){$4$}
\put(250,68){\line(1,0){3}}
\put(250,101){\line(1,0){3}}
\put(250,134){\line(1,0){3}}

\put(261,2){\line(0,1){3}}
\put(259,-8){$1$}
\put(415,2){\line(0,1){3}}
\put(409,-8){$pn$}

\put(428,-8){$i$}
\put(233,155){$\Lambda^\theta_i$}

\end{picture}

\vspace{4mm}

Figure 2. A $3$-mobile $\theta$ with $5$ black vertices
and the associated spatial contour function.

\end{center}

We will now describe the Bouttier-Di Francesco-Guitter bijection between 
$\m^p_n$ and the set of all $p$-mobiles with $n$ black vertices.
This bijection can be found in Section 2 of \cite{BDG} in the more general
setting of bipartite planar maps.  Note that \cite{BDG} deals with pointed
planar maps rather than with rooted planar maps. However,
the results described below easily follow from \cite{BDG}.

Let $\tau$ be a $p$-tree with $n$ black vertices and let $k=\#\tau -1=pn$. The
depth-first search sequence of $\tau$ is the sequence $w_0,w_1,\ldots,w_{2k}$ of vertices
of $\tau$ which is obtained by induction as follows. First $w_0=\varnothing$, and
then for every $i\in\{0,\ldots,2k-1\}$, $w_{i+1}$ is either the first child of
$w_i$ that has not yet appeared in the sequence $w_0,\ldots,w_i$, or the parent
of $w_i$ if all children of $w_i$ already appear in the sequence $w_0,\ldots,w_i$.
It is easy to verify that $w_{2k}=\varnothing$ and that all vertices of $\tau$
appear in the sequence $w_0,w_1,\ldots,w_{2k}$ (of course some of them 
appear more than once).

Vertices $w_i$ are white when $i$ is even
and black when $i$ is odd.
The contour sequence of $\tau^\circ$ is by definition the sequence 
$v_0,\ldots,v_k$ defined by $v_i=w_{2i}$ for every $i\in\{0,1,\ldots,k\}$.

\begin{center}
\begin{picture}(400,220)
\linethickness{1.5pt}
\put(260,0){\circle{12}}
\put(257,-4){$1$}
\put(255,3){\line(-2,1){52}}
\put(265,3){\line(4,3){34}}

\put(300,30){\circle*{8}}
\put(325,55){\circle{12}}
\put(303,33){\line(1,1){18}}
\put(275,55){\circle{12}}
\put(272,51){$1$}
\put(297,33){\line(-1,1){18}}
\put(322,51){$2$}

\put(344,63){\line(1,1){10}}
\put(344,63){\line(-1,1){10}}

\put(200,30){\circle*{8}}
\put(225,55){\circle{12}}
\put(203,33){\line(1,1){18}}
\put(175,55){\circle{12}}
\put(172,51){$3$}
\put(197,33){\line(-1,1){18}}
\put(222,51){$2$}
\put(175,85){\circle*{8}}
\put(175,61){\thinlines\line(0,1){20}}
\put(200,110){\circle{12}}
\put(178,88){\line(1,1){18}}
\put(150,110){\circle{12}}
\put(147,106){$4$}
\put(172,88){\line(-1,1){18}}
\put(197,106){$3$}
\put(146,114){\line(-1,1){20}}
\put(125,135){\circle*{8}}
\put(204,114){\line(1,1){20}}
\put(225,135){\circle*{8}}
\put(250,160){\circle{12}}
\put(228,138){\line(1,1){18}}
\put(200,160){\circle{12}}
\put(197,156){$2$}
\put(222,138){\line(-1,1){18}}
\put(247,156){$1$}
\put(150,160){\circle{12}}
\put(128,138){\line(1,1){18}}
\put(100,160){\circle{12}}
\put(97,156){$3$}
\put(122,138){\line(-1,1){18}}
\put(147,156){$2$}

\bezier{800}(168,55)(-60,260)(146,165)
\bezier{50}(106,160)(125,170)(143,160)
\bezier{100}(143,110)(100,132)(100,154)
\bezier{50}(156,110)(175,120)(193,110)
\bezier{50}(200,116)(190,135)(200,154)
\bezier{100}(153,165)(200,210)(250,166)
\bezier{50}(206,160)(225,170)(243,160)
\bezier{100}(206,110)(230,85)(225,61)
\bezier{50}(181,55)(200,65)(218,55)
\bezier{80}(225,49)(230,22)(256,5)

\put(300,140){\circle{12}}
\put(297,136){$\partial$}
\bezier{200}(260,6)(250,110)(294,140)
\bezier{200}(256,160)(280,182)(300,146)
\bezier{200}(275,61)(300,100)(300,134)
\bezier{150}(325,49)(320,30)(266,0)
\bezier{200}(265,-3)(400,20)(306,140)

\end{picture}

\vspace{6mm}

\end{center}

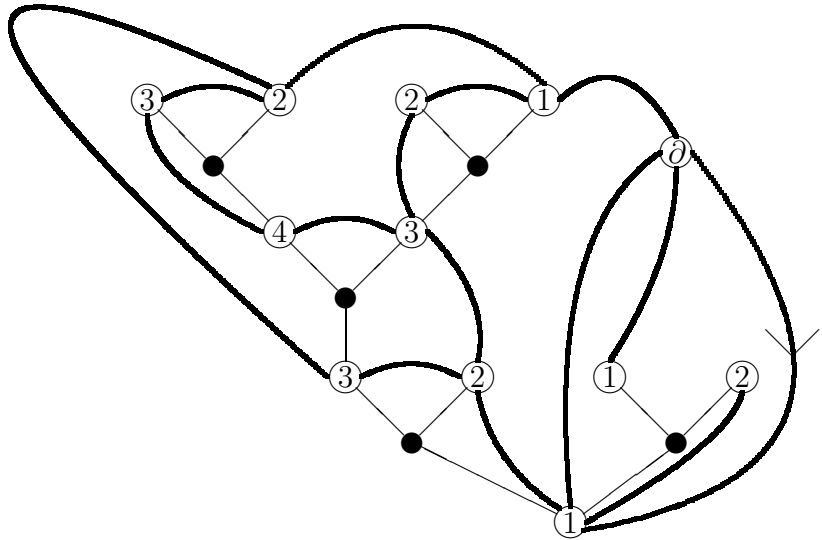
\noindent Figure 3. The Bouttier-Di Francesco-Guitter bijection: A rooted $3$-mobile with $5$ black vertices
and the associated rooted
$6$-angulation with $5$ faces. The root of the map is the edge between the vertex $\partial$
and the root of the tree at the right end of the figure.

\medskip

Now let $\theta=(\tau,(\ell_v)_{v\in\tau^\circ})$ be a $p$-mobile 
with $n$ black vertices. As previously, denote the contour
sequence of $\tau^\circ$ by $v_0,v_1,\ldots,v_{pn}$. Suppose that the tree $\tau$ is drawn 
in the plane as pictured on Fig.3 and add an extra vertex $\partial$.
We associate with $\theta$ a rooted $2p$-angulation $M$
with $n$ faces,
whose set of vertices is
$$V(M)=\tau^\circ \cup\{\partial\}$$
and whose edges are obtained by the following device: For every
$i\in\{0,1,\ldots,pn-1\}$,
\begin{description}
\item{$\bullet$} if $\ell_{v_i}=1$, draw an edge between $v_i$ and $\partial$\ ;
\item{$\bullet$} if $\ell_{v_i}\geq 2$, draw an edge between $v_i$ and $v_j$, where $j$
is the first index in the sequence $i+1,i+2,\ldots,pn$ such that $\ell_{v_j}=\ell_{v_i}-1$
(we then say that $j$ is the
successor of $i$, or sometimes that $v_j$ is a successor of $v_i$  -- note that a given 
vertex $v$ can appear several times in the contour sequence and so 
may have several different successors).

\end{description}

Notice that $v_{pn}=v_0=\varnothing$ and $\ell_\varnothing=1$, and that condition (b) 
in the definition of a $p$-tree
entails that
$\ell_{v_{i+1}}\geq \ell_{v_i}-1$ for every $i\in\{0,1,\ldots,pn-1\}$. 
This ensures that whenever $\ell_{v_i}\geq 2$ there is at least one
vertex among $v_{i+1},v_{i+2},\ldots,v_{pn}$ with label $\ell_{v_i}-1$.
The construction can be made in such a way
that edges do not intersect, except possibly at their endpoints: For every vertex $v$, each index $i$ such that 
$v_i=v$ corresponds to a ``corner'' of $v$, and the associated edge starts from
this corner. We refer to
Section 2 of \cite{BDG} for a more detailed description (here we will only need the fact that edges are
generated in the way described above). The resulting planar
graph $M$ is a $2p$-angulation, which is rooted at the oriented edge 
between $\partial$ and $v_0=\varnothing$, corresponding 
to $i=0$ in the previous construction. Each black
vertex of
$\tau$ is associated with a face of the map $M$. Furthermore the graph 
distance in $M$ between the root vertex $\partial$ and another vertex
$v\in\tau^\circ$ is equal to $\ell_v$. See Fig.3 for the $6$-angulation
associated with the $3$-mobile of Fig.2.

It follows from \cite{BDG} that the preceding construction yields
a bijection between the set $\T^p_n$ of all $p$-mobiles with $n$
black vertices and the set $\m^p_n$. 

The {\it contour function} of $\tau^\circ$ is the discrete sequence
$C^{\tau^\circ}_0,C^{\tau^\circ}_1,\ldots,C^{\tau^\circ}_{pn}$ defined by 
$$C^{\tau^\circ}_i=\frac{1}{2}\,|v_i|\ ,\hbox{ for every }0\leq i\leq pn.$$
See Fig.1 for an example with $p=n=3$.
It is easy to verify that the contour function determines $\tau^\circ$, which
in turn determines the $p$-tree $\tau$ uniquely. We will also
use the {\it spatial contour function} of $\theta=(\tau,(\ell_v)_{v\in\tau^\circ})$, which
is the discrete sequence $(\Lambda^\theta_0,\Lambda^\theta_1,\ldots,\Lambda^\theta_{pn})$ defined by
$$\Lambda^\theta_i=\ell_{v_i}\ ,\hbox{ for every }0\leq i\leq pn.$$
From property (b) of the labels and the definition of the contour sequence, it
is clear that $\Lambda^\theta_{i+1}\geq \Lambda^\theta_i-1$ for every $0\leq i\leq pn-1$
(cf Fig.2). 
The pair $(C^{\tau^\circ},\Lambda^\theta)$ determines $\theta$ uniquely. 

Define an equivalence
relation $\sim_{[\tau]}$ on $\{0,1,\ldots,pn\}$ by setting $i\sim_{[\tau]} j$
if and only if $v_i=v_j$. The quotient space $\{0,1,\ldots,pn\}\,/\!\sim_{[\tau]}$
is then obviously identified with $\tau^\circ$. 
If $i\leq j$, the
relation $i\sim_{[\tau]} j$ implies
$$\inf_{i\leq k\leq j} C^{\tau^\circ}_k=C^{\tau^\circ}_i=C^{\tau^\circ}_j.$$
The converse is not true (except if $p=2$) but 
the conditions $j>i+1$, $C^{\tau^\circ}_i=C^{\tau^\circ}_j$ and
$$C^{\tau^\circ}_k>C^{\tau^\circ}_i\ ,\hbox{ for every }
k\in\,]i,j[\cap \Z$$
imply that $i\sim_{[\tau]} j$. 

\subsection{Convergence towards the Brownian map}

For every integer $n\geq 1$, let $M_n$ be a random rooted $p$-angulation, which is uniformly distributed 
over the set $\m^n_p$, 
and let
$\theta_n=(\tau_n,(\ell^n_v)_{v\in\tau^\circ_n})$ be the random mobile
corresponding to $M_n$ via the Bouttier-Di Francesco-Guitter bijection. Then
$\theta_n$ is uniformly distributed
over the set $\T^p_n$ of all $p$-mobiles with $n$ black vertices. 
We denote by $C^n=(C^n_i)_{0\leq i\leq pn}$ the contour function of
$\tau_n^\circ$ and by
$\Lambda^n=(\Lambda^n_i)_{0\leq i\leq pn}$ the spatial contour function of $\theta_n$.
Recall that the  pair $(C^n,\Lambda^n)$ determines $\theta_n$ and thus $M_n$.

Let $\mm_n$ stand for the vertex set of $M_n$. Thanks to 
the Bouttier-Di Francesco-Guitter bijection we have the
identification
$$\mm_n=\tau^\circ_n\cup\{\partial_n\}$$
where $\partial_n$ denotes the root vertex of $M_n$. The graph distance 
on $\mm_n$ will be denoted by $d_n$. 
In particular, if $a,b\in\tau^\circ_n$, $d_n(a,b)$ denotes the
graph distance between $a$ and $b$ viewed as vertices in
the map $M_n$.

To simplify notation, we write $\sim_{[n]}$ for the equivalence relation $\sim_{[\tau_n]}$
on $\{0,1,\ldots,pn\}$, so that $\tau^\circ_n$
is canonically identified to the quotient $\{0,1,\ldots,pn\}/\!\sim_{[n]}$. 
We also write $p_n$ for the canonical projection from
$\{0,1,\ldots,pn\}$ onto $\tau^\circ_n=\mm_n\backslash\{\partial_n\}$.
To be specific, $p_n(i)=v^n_i$, if $v^n_0,v^n_1,\ldots,v^n_{pn}$ denotes the
contour sequence of $\tau_n^\circ$. If $i,j\in\{0,1,\ldots,pn\}$,
we set $d_n(i,j)=d_n(p_n(i),p_n(j))$. 

The following theorem restates the main result of \cite{IM}
in a form convenient for the present work.
To simplify notation, we set
$$
\lambda_p={1\over 2}\sqrt{\frac{p}{p-1}}\ ,\quad
\kappa_p= \Big(\frac{9}{4p(p-1)}\Big)^{1/4}.
$$

\begin{theorem}
\label{mainIM}
From every sequence of integers converging to $\infty$, we can extract
a subsequence $(n_k)_{k\geq 1}$
such that the following properties hold. On a suitable probability space, 
for every integer $n$ belonging to the
sequence $(n_k)$ we can construct the uniformly distributed random
$p$-angulation $M_n$, in such a way that:
\begin{eqnarray}
\label{basicconv}
&&\hspace{-6mm}\left(\lambda_p\,n^{-1/2}\,C^n_{\lfloor pnt\rfloor },
\kappa_p\,n^{-1/4} \Lambda^n_{\lfloor pnt\rfloor },
\kappa_p\,n^{-1/4}\,d_n({\lfloor pns\rfloor },{\lfloor pnt\rfloor })\!\right)
_{0\leq s\leq 1, 0\leq
t\leq 1}\nonumber\\
\noalign{\smallskip}
&&\qquad\build{\la}_{n\to\infty}^{} \left(\ov\eg_t,\ov Z_t
,D(s,t)\right)_{0\leq s\leq 1,0\leq t\leq
1}, \qquad a.s.
\end{eqnarray}
where the convergence holds uniformly in $s,t\in[0,1]$ along the 
sequence $(n_k)_{k\geq 1}$. In {\rm(\ref{basicconv})},
$(\ov\eg_t,\ov Z_t)_{0\leq t\leq 1}$ has the distribution 
described in subsection 2.4, and
$(D(s,t))_{0\leq s\leq 1,0\leq t\leq 1}$ is a continuous random process
that satisfies properties (i)--(iv) stated in subsection 2.5. 
Furthermore, the pointed compact metric spaces 
$$(\mm_n,\kappa_p\,n^{-1/4}d_n,\partial_n)$$
converge almost surely to $(\mm_\infty,D,\rho)$ in the sense of the
Gromov-Hausdorff convergence.
\end{theorem}

\rems (a) The convergence of the first two components in (\ref{basicconv})
does not require the use of a subsequence: See Theorem 3.3 in \cite{We}.
A subsequence is needed only to get the convergence of the third component
via a compactness argument.

\noindent (b)
The last assertion of the theorem refers to the Gromov-Hausdorff
distance on the space of isometry classes of pointed compact
metric spaces. The definition of this distance is recalled 
in Section 8 below (this definition is not needed in the proof of our main results, and our main tool will be the convergence (\ref{basicconv})). The last assertion is then
a rather simple consequence of the convergence (\ref{basicconv}), and the
fact that the process $D$ satisfies the above-mentioned properties (i)--(iv):
See the proof of Theorem \ref{rerootinv} below for a sketch of the argument.
The key point in the proof of Theorem \ref{mainIM} is to verify property (ii)
for the pseudo-metric $D$. See \cite{IM} for more details.

\smallskip
We will need a simple application of (\ref{basicconv})
to the convergence of ``discrete snakes'' associated with the
$p$-mobiles $(\tau_n,(\ell^n_v)_{v\in\tau^\circ_n})$ . Recall that
 the contour sequence of $\tau^\circ_n$ is denoted 
by $v^n_0,v^n_1,\ldots,v^n_{pn}$. Then, for every 
$i\in\{0,1,\ldots,pn\}$, define the finite sequence $W^n_i=(W^n_i(j),0\leq j\leq C^n_i)$
by requiring that $W^n_i(j)=\ell^n_{u^n_i(j)}$, where $u^n_i(j)\in\tau^\circ_n$ is the ancestor of 
$v^n_i$ at generation $2j$ in the tree $\tau_n$. In particular,
$W^n_i(C^n_i)=\Lambda^n_i$. Then, if (\ref{basicconv}) holds, we have also, along the sequence $(n_k)$,
\begin{equation}
\label{snakeconv}
\sup_{0\leq s\leq 1}\Big(\sup_{r\geq 0} \Big|\kappa_pn^{-1/4} W^n_{\lfloor pnt\rfloor}
(\lfloor\lambda_p^{-1} n^{1/2}r\rfloor \wedge C^n_{\lfloor pnt\rfloor})
-\ov W_s(r\wedge \ov \eg_s)\Big| \Big)
\build{\la}_{n\to\infty}^{} 0 \qquad\hbox{a.s.}
\end{equation}
where the process $(\ov W_s)_{0\leq s\leq 1}$ is defined from the pair
$(\ov\eg_s,\ov Z_s)_{0\leq s\leq 1}$ as explained at the end of subsection 2.4.
The convergence (\ref{snakeconv}) is a consequence of
the convergence of the first two components in 
(\ref{basicconv}). A simple way to verify this is to use the homeomorphism
theorem of 
 \cite{MaMo0}. We leave details 
 to the reader.

\smallskip
In the remaining part of this work (with the important exception of Section 9), we choose a sequence 
$(n_k)$ as in Theorem \ref{mainIM} and we consider only values of $n$
belonging to this sequence. We fix the random maps $M_n$ 
as in the theorem, and we argue on the set of full probability
measure where the convergences (\ref{basicconv}) and (\ref{snakeconv}) hold, and the
process $(D(s,t))_{0\leq s\leq 1,0\leq t\leq 1}$ satisfies the
properties (i)--(iv).

\section{Preliminary results}

In this section, we state and prove a few preliminary facts that
will be used in the subsequent proofs. We start with 
a result which was already mentioned in Section 1. Recall the notation
${\rm Skel}_\infty=\Pi({\rm Sk}(\t_{\ov\eg}))$.

\begin{proposition}
\label{skeleton}
The following properties hold almost surely. The restriction of the projection $\Pi$ to ${\rm Sk}(\t_{\ov\eg})$ is a homeomorphism
from ${\rm Sk}(\tree)$ onto ${\rm Skel}_\infty$,
and the Hausdorff dimension of ${\rm Skel}_\infty$ is less than or equal to $2$.
\end{proposition}

\proof
We already know that the projection $\Pi$ from $\t_{\ov\eg}$ onto $\mm_\infty$
is continuous, and so is its restriction to
${\rm Sk}(\t_{\ov\eg})$. This restriction is also
one-to-one by Lemma \ref{equirel}.  We need to verify that its inverse is continuous
in order to get the first assertion.
To see this, let $(x_k)_{k\geq 1}$ be a sequence in ${\rm Skel}_\infty$
that converges to a point $x_\infty\in {\rm Skel}_\infty$, in the sense of the
metric $D$. For every $k\geq 1$, we have $x_k=\Pi(a_k)$ with $a_k\in{\rm Sk}(\t_{\ov\eg})$
and similarly $x_\infty=\Pi(a_\infty)$ with $a_\infty\in{\rm Sk}(\t_{\ov\eg})$.
Since $\t_{\ov\eg}$ is compact, we may find a subsequence $(a_{k_j})_{j\geq 1}$
that converges to a point $b\in\t_{\ov\eg}$. By the continuity of $\Pi$, we have
then $\Pi(b)=x=\Pi(a_\infty)$. Since $a_\infty\in{\rm Sk}(\t_{\ov\eg})$, 
Lemma \ref{equirel} shows that this is 
possible only if $b=a_\infty$. We conclude that the sequence $(a_k)$ must converge
to $a_\infty$ as required. 

As for the second assertion, we observe that ${\rm Sk}(\t_{\ov\eg})$ is a countable 
union of sets that are isometric to line segments, so that its Hausdorff dimension is $1$. 
On the other hand, the bound (iv) in subsection 2.5 easily implies that the mapping
$\Pi$ from $(\t_{\ov\eg},d_{\ov\eg})$ onto $(\mm_\infty,D)$ is H\"older continuous with exponent
$\frac{1}{2}-\ve$, for every $\ve>0$ (see Lemma 5.1 in \cite{IM} for a
closely related statement). The desired result follows.
\cq

\rem It is natural to conjecture that the Hausdorff dimension of ${\rm Skel}_\infty$
is equal to $2$.

\smallskip

Our next lemma shows that for any
fixed time $U\in\,]0,1[$, $p_{\ov\eg}(U)$ does not belong to 
the skeleton of $\t_{\ov\eg}$.
Although this statement is intuitively clear, it is not so easy to give a 
precise argument.

\begin{lemma}
\label{leaf}
Let $U\in\,]0,1[$. Then $p_{\ov\eg}(U)\notin {\rm Sk}(\t_{\ov\eg})$ and $\pp(U)\notin {\rm Skel}_\infty$, almost surely.
\end{lemma}

\proof To get the first assertion, we need to prove that for every $\ve\in\,]0,U\wedge(1-U)[$,
we have almost surely
$$m_{\ov\eg}(U, U+\ve)<\ov\eg_U\quad\hbox{and}\quad
m_{\ov\eg}(U-\ve, U)<\ov\eg_U.$$
Since the law of the pair $(\eg_t,Z_t)_{0\leq t\leq 1}$ is invariant under the time-reversal
operation $t\to 1-t$, the same holds for the pair $(\ov\eg_t,\ov Z_t)_{0\leq t\leq 1}$, and
so it is enough to prove the first part of the preceding assertion. Equivalently, we must
show that $P(A)=0$, where
$$A=\bigcup_{\ve\in\,]0,1-U[}\{m_{\ov\eg}(U, U+\ve)=\ov\eg_U\}.$$

Now recall the definition of $s_*$ and $\ov\eg$
in subsection 2.4. It follows from the invariance under re-rooting stated in
subsection 2.3 of \cite{LGW} that $s_*$ is independent of $\ov\eg$ and uniformly
distributed over $[0,1]$. In particular, if we fix $\ve_0\in\,]0,1-U[$, we have
$$(1-U-\ve_0)P(A)
=P(\{s_*<1-U-\ve_0\}\cap A).$$
On the other hand, from the definition of $\ov\eg$ in terms of $\eg$, we get
$$(\{s_*<1-U-\ve_0\}\cap A)\subset (A_1 \cup A_2)$$
where
$$A_1=\{s_*<1-U-\ve_0\}\cap \Big(\bigcup_{\ve\in\,]0,\ve_0[} 
\{m_\eg(s_*+U, s_*+U+\ve)=\eg_{s_*+U}\}\Big)$$
and 
$$A_2=\{s_*<1-U-\ve_0\}\cap \{m_\eg(s_*,s_*+U)=\eg_{s_*+U}\}.$$
Let us verify that $P(A_1)=0$. For every rational $q\in\,]0,1[$, set
$$s_*(q)=\sup\{s\in[0,q]: Z_s= \min_{0\leq r\leq q} Z_r\}.$$
Clearly the random variable $s_*(q)$ is measurable with respect to
the $\sigma$-field generated by $(\eg_r,W_r)_{0\leq r\leq q}$. As a consequence,
the random time
$$T_{(q)}=q + (U+s_*(q)-q)_+$$
is a stopping time of the filtration generated by $(\eg_r,W_r)_{0\leq r\leq 1}$. Notice that
the process $(\eg_t)_{0\leq t\leq 1}$ is Markovian with respect to this filtration. By a standard
property of linear Brownian motion we get that, almost surely on the event
$\{T_{(q)}<1\}$, we have for every $\ve>0$,
$$m_\eg(T_{(q)},(T_{(q)}+\ve)\wedge 1)<\eg_{T_{(q)}}.$$
Finally, on the event $\{s_*<1-U-\ve_0\}$, we can pick a rational $q$
such that $s_*<q<s_*+U$, and then we have $s_*(q)=s_*$,
and $T_{(q)}=s_*+U$. From the preceding observations we see that
$A_1$ does not hold almost surely. A similar argument gives $P(A_2)=0$.
It follows that $P(A)=0$, which completes the proof of the first assertion. 
The second one follows since the equivalence class of any
$a\in {\rm Sk}(\t_{\ov\eg})$ for $\approx$ is a singleton, by Lemma \ref{equirel}. \cq

\smallskip
Recall from Section 1 the definition of a geodesic connecting
two points of $\mm_\infty$. Also recall that $D(\rho,x)={\ov Z}_x$, for every
$x\in\mm_\infty$.

\begin{lemma}
\label{skeletongeo}
Almost surely, for every $x\in\mm_\infty$, every geodesic $\omega=(\omega(t))_{0\leq t\leq {\ov Z}_x}$
from $\rho$ to $x$, and every interval $[u,v]$ such that $0\leq u<v\leq {\ov Z}_x$, the range of
$\omega$ over $[u,v]$ intersects $\mm_\infty\backslash {\rm Skel}_\infty$.
\end{lemma}

\proof 
Let $x\in \mm_\infty$ and let $\omega$ be a geodesic from $\rho$ to
$x$. Also, let $u$ and $v$ be reals such that $0\leq u<v<{\ov Z}_x$. We argue by contradiction
and suppose that the set $\{\omega(t): u\leq t\leq v\}$ is contained in ${\rm Skel}_\infty$.
From Proposition \ref{skeleton}, we may write $\omega(r)=\Pi(\psi(r))$, for
every $r\in[u,v]$, where $\psi$ is a continuous mapping 
from $[u,v]$ into ${\rm Sk}(\t_{\ov\eg})$. This mapping $\psi$ is also one-to-one
because $\omega$ is a geodesic. Since ${\rm Sk}(\t_{\ov\eg})$ is a
(non-compact) real tree, it follows that the range of $\psi$ is the 
line segment $\llbracket \psi(u),\psi(v)\rrbracket$, and moreover 
$\psi$ is a homeomorphism from $[u,v]$ onto $\llbracket \psi(u),\psi(v)\rrbracket$.
The fact that ${\ov Z}_{\psi(r)}=r$ for every $r\in[u,v]$ then implies that
the mapping $a\la {\ov Z}_a$ is monotone increasing over $\llbracket \psi(u),\psi(v)\rrbracket$.
However, this is absurd since $\ov Z$ has been constructed by
shifting the process $Z$, and we already observed that 
$Z$ can be interpreted as Brownian motion indexed by the tree $\t_\eg$,
so that it cannot vary monotonically on a line segment of the tree. \cq

\smallskip
Recall our notation $\tree(c)$ for the subtree of descendants of a vertex $c\in\tree$.
Before we state the next lemma, we make a simple remark
concerning the boundary of the set $\Pi(\t_{\ov\eg}(c))$, when $c$
is a vertex of ${\rm Sk}(\tree)$. We claim that
this boundary consists of the point $\Pi(c)$ and the points $\Pi(d)$,
for all $d\in\tree(c)$ such that there exists $d'\in\tree\backslash\tree(c)$
with $d'\approx d$. To see this, first note that $\Pi(c)$ belongs to the 
boundary since the set $\Pi(\llbracket\rho,c\llbracket)$
is contained in $\mm_\infty\backslash\Pi(\tree(c))$ by Lemma \ref{equirel}.
Similarly, for every point $d\in\tree(c)$ such that there exists $d'\in\tree\backslash\tree(c)$
with $d'\approx d$, the point $\Pi(d)$ belongs to the boundary. 
Conversely, let $x$
be a point of the boundary of $\Pi(\t_{\ov\eg}(c))$, and write $x=\Pi(b)$ with 
$b\in\tree(c)$. We may assume that $b\not =c$. Then we have $\Pi(b)=\lim \Pi(a_k)$
where $(a_k)_{k\geq 1}$ is a sequence in $\tree\backslash\tree(c)$. Note that
$(\tree\backslash\tree(c))\cup\{c\}$ is compact as the image under 
$p_{\ov\eg}$ of the union of two closed subintervals of $[0,1]$, and so we may assume that
the sequence $(a_k)$ converges to $a\in (\tree\backslash\tree(c))\cup\{c\}$. 
Then $\Pi(b)=\Pi(a)$ and $a\not = c$ since we assumed that $b\not =c$. 
This gives our claim.

\begin{lemma}
\label{labelsline}
Almost surely, for every $a,b\in\t_{\ov\eg}$, and every continuous curve
$(\omega(t),0\leq t\leq T)$  
in $\mm_\infty$, such that $\omega(0)=\Pi(a)$ and $\omega(T)=\Pi(b)$, we have
\beq
\label{labelsline1}
{\ov Z}_c\geq \inf\{{\ov Z}_{\omega(t)}:0\leq t\leq T\},
\eeq
for every $c\in \llbracket a\tri b,b\rrbracket$.
Furthermore, if equality holds in {\rm(\ref{labelsline1})}, 
then there exists $t\in[0,T]$ such that $\omega(t)= \Pi(c)$.
\end{lemma}

\proof When $c=b$ there is nothing to prove. Suppose first that 
$c\in\,\rrbracket a\tri b,b\llbracket$. We can also assume that
$\Pi(c)\notin\{\omega(t):0\leq t\leq T\}$. Note that
$b\in\t_{\ov\eg}(c)$ but $a\notin \t_{\ov\eg}(c)$. Set
$$t_0=\inf\{t\in[0,T]:\omega(t)\in \Pi(\t_{\ov\eg}(c))\}.$$
Notice that $\Pi(\t_{\ov\eg}(c))$ is closed as the image under
the projection $\pp$ of a compact subinterval of
$[0,1]$. It follows that $\omega(t_0)\in\Pi(\t_{\ov\eg}(c))$. 

If $t_0=0$,
then $\Pi(a)=\omega(0)=\Pi(a')$ for some $a'\in\t_{\ov\eg}(c)$. Thus 
$a\approx a'$, which from the definition of the equivalence relation $\approx$
entails ${\ov Z}_c\geq {\ov Z}_a={\ov Z}_{a'}$
(note that $c\in\llbracket a,a'\rrbracket$, so that, if $a=p_{\ov\eg}(s)$ and $a'=p_{\ov\eg}(t)$, the set
$\{p_{\ov\eg}(r):s\wedge t\leq r\leq s\vee t\}$ contains $c$). There is even a strict inequality
because otherwise we would have $c\approx a$, which is impossible by Lemma \ref{equirel},
since $c$ belongs to the skeleton of the tree $\t_{\ov\eg}$.

If $t_0>0$, then $\omega(t_0)$ belongs to the boundary of
$\Pi(\t_{\ov\eg}(c))$ in $\mm_\infty$. We appply the remark
preceding the statement of the lemma, noting that 
$\omega(t_0)\not =\Pi(c)$
since $\Pi(c)\notin\{\omega(t):0\leq t\leq T\}$. Thus, there exists a point $d\in\t_{\ov\eg}(c)$
such that $\omega(t_0)=\Pi(d)=\Pi(d')$ for some 
$d'\in\t_{\ov\eg}\backslash\t_{\ov\eg}(c)$. As in the case $t_0=0$, this entails
that ${\ov Z}_c\geq {\ov Z}_d$, and there is even a strict inequality.

Finally, if $c=a\tri b$, the bound (\ref{labelsline1}) 
is immediate by a continuity argument. To get the desired conclusion when there 
is equality in (\ref{labelsline1}), we can assume that $a\not =a\tri b$
and $b\not = a\tri b$.
We then argue in a similar way as in the case 
$c\in\,\rrbracket a\tri b,b\llbracket$, but we replace
$\t_{\ov\eg}(c)$ by $p_{\ov\eg}([r_1,r_2])$, where $[r_1,r_2]$ is the unique 
compact subinterval of $[0,1]$ such that $p_{\ov\eg}(r_1)=p_{\ov\eg}(r_2)=a\tri b$,
$b\in p_{\ov\eg}([r_1,r_2])$ and $a\notin p_{\ov\eg}([r_1,r_2])$. We leave details
to the reader. \cq

\section{Two classes of geodesics}

In this section, we will discuss 
several classes of geodesics connecting the root to a point of $\mm_\infty$. This point
will be of the form $\pp(s)$ for $s\in[0,1]$. All the results of this section are valid 
outside a set of zero probability (discarding such a set is needed e.g. to apply Lemma
\ref{equirel} or Lemma \ref{labelsline}). It is important to note that this 
set of zero probability does not depend on $s$. We will omit the words ``almost surely''
in the statements of this section.

\subsection{Simple geodesics}

For every $s\in[0,1]$, we define a mapping $\varphi_s:[0,\ov Z_s]\la [0,1]$ by setting
$$\vf_s(t)=\sup\{r\leq s: {\ov Z}_r\leq t\}\;,\ 0\leq t\leq {\ov Z}_s.$$
Clearly, $\pp(\vf_s(0))=\rho$ and $\pp(\vf_s({\ov Z}_s))=\pp(s)$. Also, 
${\ov Z}_{\vf_s(t)}=t$ and if $0\leq t\leq t'\leq {\ov Z}_s$,
$$D(\vf_s(t),\vf_s(t'))\leq {\ov Z}_{\vf_s(t)} + {\ov Z}_{\vf_s(t')}-2\inf_{[\vf_s(t),\vf_s(t')]} {\ov Z}_r
=t+t'-2t=t'-t.$$
Since $D(\vf_s(t),\vf_s(t'))\geq D(0,\vf_s(t'))-D(0,\vf_s(t))=t'-t$, we must have 
 $$D(\vf_s(t),\vf_s(t'))=t'-t$$ for every $0\leq t\leq t'\leq {\ov Z}_s$. Thus 
 $(\pp(\vf_s(t)),0\leq t\leq {\ov Z}_s)$ is a geodesic from $\rho$
 to the point $\pp(s)$. We will write
 $$\Phi_s(t)=\pp(\vf_s(t))$$
 for every $t\in[0,\ov Z_s]$ and call $\Phi_s$
 a {\it simple} geodesic from
 $\rho$ to $\pp(s)$ (sometimes we abusively say that $\vf_s$
 itself is a simple geodesic). We can also give a dual definition 
 of simple geodesics, by setting, for every $s\in[0,1]$,
 $$\ov \vf_s(t)=\inf\{r\geq s: {\ov Z}_r\leq t\}\;,\ 0\leq t\leq {\ov Z}_s.$$
It is immediate that $\vf_s(t)\approx\ov\vf_s(t)$, so that
$\pp(\ov\vf_s(t))=\pp(\vf_s(t))=\Phi_s(t)$, for every $t\in[0,Z_s]$.
 
 Let $s$ and $s'$ be two distinct points in $]0,1[$ such that $\pp(s)=\pp(s')$.
 If $s\approx s'$, we have
 $${\ov Z}_s={\ov Z}_{s'}=\min_{[s\wedge s',s\vee s']} {\ov Z}_r$$
 and it readily follows that the geodesics $\Phi_s$ and $\Phi_{s'}$
 coincide. On the other hand, if $s\sim s'$ (so that $s\approx s'$ does not
 hold), it is not hard to see that the geodesics $\Phi_s$ and $\Phi_{s'}$
 differ. We will come back to this later.
 
 It is important to note that a simple geodesic from $\rho$ to $\pp(s)$ does not intersect 
 ${\rm Skel}_\infty$, except possibly at its endpoint $\pp(s)$. 
 Indeed, if $0\leq t<\ov Z_s$, it is immediate from the definition that $\vf_s(t)$ is a 
 right-increase time of $\ov Z$, and then Lemma \ref{equirel} implies that
 $\pp(\vf_s(t))\notin {\rm Skel}_\infty$. Much of what follows is devoted to
 proving that a similar property holds for any geodesic connecting the root
 to another point of $\mm_\infty$. 
 
 \subsection{Minimal geodesics}
 
 Let $s\in[0,1]$, and denote by ${\rm G}(s)$ the set of all
 mappings $\ga:[0,{\ov Z}_s]\la[0,s]$ such that $(\pp(\ga(t)),0\leq t\leq {\ov Z}_s)$
 is a geodesic from $\rho$ to $\pp(s)$. Notice that ${\rm G}(s)$
 is not empty since $\vf_s\in {\rm G}(s)$. We then 
 define a mapping $\ga_s:[0,\ov Z_s]\la [0,1]$ by setting
 $$\gamma_s(t)=\inf\{\gamma(t):\gamma\in {\rm G}(s)\}\;,\ 0\leq t\leq {\ov Z}_s.$$
 In contrast with $\vf_s$, $\ga_s$ only depends on $\pp(s)$ and not on
 the particular choice of a representative of $\pp(s)$ in $[0,1]$.
 
 We have $\ga_s(0)=0$ and $\pp(\ga_s({\ov Z}_s))=\pp(s)$. From the
 continuity of ${\ov Z}$ and of the distance $D$, it is also clear that
 \beq
 \label{minitech1}
 {\ov Z}_{\ga_s(r)}=D(0,\ga_s(r))=r
 \eeq
 and
 \beq
 \label{minitech2}
 D(\ga_s(r),s)={\ov Z}_s-r
 \eeq
 for every $r\in[0,{\ov Z}_s]$. 
 
It is easy to verify that the mapping $t\to \ga_s(t)$ is monotone increasing.
 Indeed, let $t\in\,]0,{\ov Z}_s]$ and $\ga\in {\rm G}(s)$. Put
 $$\widetilde\ga(r)=\left\{\begin{array}{ll}
 \gamma(r)&\hbox{if }r\in[t,{\ov Z}_s],\\
 \vf_{\gamma(t)}(r)&\hbox{if }r\in[0,t[.
 \end{array}
 \right.
 $$
 Then, we have also $\widetilde\ga\in {\rm G}(s)$. Furthermore,
 $\widetilde\ga(r)< \ga(t)$ for every $r\in[0,t[$. From the definition of
 $\ga_s$ it now follows that $\ga_s(r)< \ga_s(t)$ for every $0\leq r<t\leq {\ov Z}_s$
 (note that $\ga_s(r)=\ga_s(t)$ is impossible since $ {\ov Z}_{\ga_s(r)}=r\not = t=
  {\ov Z}_{\ga_s(t)}$).
  
  To simplify notation, we set $\Gamma_s(t)=\pp(\ga_s(t))$ for every $t\in[0,\ov Z_s]$. 
  
  \begin{proposition}
 \label{minigeo}
 The curve $(\Gamma_s(t),0\leq t\leq {\ov Z}_s)$
 is a geodesic from  $\rho$ to $\pp(s)$. It is called
 the minimal geodesic from $\rho$ to $\pp(s)$. (Sometimes we also
 say that $\ga_s$ is a minimal geodesic.)
 \end{proposition}
 
 The proof of Proposition \ref{minigeo} is based on the following lemma.
 We denote by  ${\rm G}^*(s)$ the set of all $\ga\in {\rm G}(s)$
 such that $\ga(t)$ is the smallest representative of $\pp(\gamma(t))$,
 for every $t\in[0,{\ov Z}_s]$. Obviously we may replace $ {\rm G}(s)$
 with $ {\rm G}^*(s)$ in the definition of $\gamma_s$.

 \begin{lemma}
 \label{minilemma}
 If $\ga$ and $\ga'$ belong to ${\rm G}^*(s)$, then $\gamma\wedge\ga'$
 also belongs to ${\rm G}^*(s)$.
 \end{lemma}
 
 \proof We first claim that the set $\{t\in\,]0,{\ov Z}_s[:\ga(t)<\ga'(t)\}$
 is open. To see this, let $t_1$ be an element of this set. We argue by
 contradiction and assume that there
 exists a sequence $(r_k)_{k\geq 1}$ in $]0,{\ov Z}_s[$ such that
 $r_k\to t_1$ and $\ga(r_k)\geq \ga'(r_k)$ for every $k\geq 1$. Via
 a compactness argument, we may further assume that $\ga(r_k)\la u$
 and $\ga'(r_k)\la u'$ as $k\to\infty$, where $0\leq u'\leq u\leq 1$. By the continuity
 of geodesics, we have $\pp(u)=\pp(\ga(t_1))$ and $\pp(u')=\pp(\ga'(t_1))$. 
 In particular,
 $${\ov Z}_u={\ov Z}_{\ga(t_1)}=t_1={\ov Z}_{\ga'(t_1)}={\ov Z}_{u'}.$$
 Since $\ga'(t_1)$ is the smallest representative of $\pp(\ga'(t_1))=\pp(u')$,
 we have $\ga'(t_1)\leq u'$, and we also know that $\ga(t_1)<\ga'(t_1)$
 by assumption. Thus we have both $\ga(t_1)<\ga'(t_1)\leq u'\leq u$
 and $\pp(u)=\pp(\ga(t_1))$. Note that $\ga(t_1)\approx u$ is impossible,
 because the previous observations 
 (and the fact that ${\ov Z}_{\ga(t_1)}={\ov Z}_{\ga'(t_1)}$) would imply that we have
 also $\gamma(t_1)\approx \gamma'(t_1)$ and then $\gamma(t_1)= \gamma'(t_1)$,
 contradicting our initial assumption. Thus, we have $\ga(t_1)\sim u$, and the 
 inequalities $\ga(t_1)<\ga'(t_1)\leq u$ imply that $p_{\ov\eg}(\ga'(t_1))$
 is a descendant of $p_{\ov\eg}(\ga(t_1))$ in the tree $\t_{\ov\eg}$. Equivalently,
 $p_{\ov\eg}(\ga(t_1))\in\, \rrbracket\rho,p_{\ov\eg}(\ga'(t_1))\llbracket$, and
 in particular $p_{\ov\eg}(\ga(t_1))\in {\rm Sk}(\t_{\ov\eg})$.
 
 Consider first the case when $p_{\ov\eg}(\ga(t_1))\in \,\rrbracket\rho,p_{\ov\eg}(\ga'(t_1))
 \tri p_{\ov\eg}(s)\llbracket$. The fact that $\ga'(t_1)$ is the smallest
 representative of $\pp(\ga'(t_1))$ entails that
 $$u'=\liminf_{k\to\infty} \ga'(r_k)\geq \ga'(t_1).$$
 On the other hand, the fact that $\pp(\ga(r_k))$ converges to $\pp(\gamma(t_1))$
 entails that 
 $$u=\limsup_{k\to\infty} \ga(r_k) \leq v$$
 where $v$ is the greatest  representative of $p_{\ov\eg}(\ga(t_1))$ in $[0,s]$
 (at this point we use in a crucial way the fact that $\ga$ takes values in $[0,s]$). Since
 $p_{\ov\eg}(\ga(t_1))\in \llbracket\rho,p_{\ov\eg}(\ga'(t_1))
 \tri p_{\ov\eg}(s)\llbracket$, 
 simple considerations about the genealogy of the tree $\t_{\ov\eg}$ show
 that we must have $v < \ga'(t_1)$. It follows that $u\leq v<\ga'(t_1)\leq u'$, which
 is a contradiction.
 
 Let us turn to the case when $p_{\ov\eg}(\ga(t_1))\in \llbracket p_{\ov\eg}(\ga'(t_1))
 \tri p_{\ov\eg}(s),p_{\ov\eg}(\ga'(t_1))\llbracket$. In that case,
 we apply Lemma \ref{labelsline} to the continuous curve
 $(\pp(\ga'({\ov Z}_s-t)),0\leq t\leq {\ov Z}_s-t_1)$. The assumptions of this lemma
 are satisfied with $a=p_{\ov\eg}(s)$, $b=p_{\ov\eg}(\ga'(t_1))$ and $c=p_{\ov\eg}(\ga(t_1))$. Since 
 $${\ov Z}_{\ga(t_1)}=t_1=\inf\{{\ov Z}_{\ga'(t)}:t_1\leq t\leq {\ov Z}_s\},$$
 we get that $\pp(\ga(t_1))=\pp(\ga(t^*))$ for some $t^*\in[t_1,{\ov Z}_s]$.
 Necessarily $t^*=t_1$ and $\ga(t_1)=\ga'(t_1)$, which is again a contradiction. This
 completes the proof of the claim.
 
 We can now complete the proof of Lemma \ref{minilemma}. 
  Let us fix $r\in\,]0,{\ov Z}_s[$ and verify that, for every $t\in[r,{\ov Z}_s[$,
 \beq
 \label{minitech3}
 D(\ga\wedge\ga'(r),\ga\wedge\ga'(t))=t-r.
 \eeq
 Notice that this equality holds for $t={\ov Z}_s$ since both $\ga({\ov Z}_s)$
 and $\ga'({\ov Z}_s)$ are equal to the smallest representative of $\pp(s)$. 
 If $\ga(r)=\ga'(r)$, (\ref{minitech3}) holds because $\pp(\ga(t))$
 and $\pp(\ga'(t))$ are both geodesics. We can thus assume
 that $\ga(r)\not =\ga'(r)$, and by symmetry we may restrict our attention
 to the case $\ga(r)<\ga'(r)$.
 
 We again argue by contradiction and assume that (\ref{minitech3})
 fails for some $t\in\,]r,{\ov Z}_s[$. We then set
 $$t_0=\inf\{t\in\,]r,{\ov Z}_s[:D(\ga\wedge\ga'(r),\ga\wedge\ga'(t))\not=t-r\} < {\ov Z}_s.$$
 We consider the following three cases:
 
 \smallskip
 $\bullet$ If $\ga(t_0)=\ga'(t_0)$, then
 $$D(\ga\wedge\ga'(r),\ga\wedge\ga'(t_0))=t_0-r$$
 and for every $t\in[t_0,{\ov Z}_s]$,
 $$D(\ga\wedge\ga'(t_0),\ga\wedge\ga'(t))=t-t_0.$$
 Thus 
 $$D(\ga\wedge\ga'(r),\ga\wedge\ga'(t))\leq t-r$$
 for every $t\in[t_0,{\ov Z}_s]$. The reverse inequality is also clear by writing
 $$D(\ga\wedge\ga'(r),\ga\wedge\ga'(t))
 \geq D(0,s)-D(0,\ga\wedge\ga'(r))-D(\ga\wedge\ga'(t),s)$$
 and we get a contradiction with the definition of $t_0$.
 
 \smallskip
 $\bullet$ If $\ga(t_0)<\ga'(t_0)$, then the claim gives $\ve>0$ such that 
 $\ga(t)<\ga'(t)$ for every $t\in[t_0,t_0+\ve[$, and thus
 $$D(\ga\wedge\ga'(r),\ga\wedge\ga'(t))=D(\ga(r),\ga(t))=t-r.$$
 This again contradicts the definition of $t_0$. 
 
 \smallskip
 $\bullet$. If $\ga(t_0)>\ga'(t_0)$, then $t_0>r$ (recall that $\ga(r)<\ga'(r)$) and the claim gives  
 $\ve>0$ such that $\ga(t)>\ga'(t)$ for every $t\in\,]t_0-\ve,t_0+\ve[$. Then, for every
 $t,t'\in\,]t_0-\ve,t_0+\ve[$, with $t<t'$,
 $$D(\ga\wedge\ga'(t),\ga\wedge\ga'(t'))=D(\ga'(t),\ga'(t'))=t'-t.$$
 Fix $t\in\,]t_0-\ve,t_0[$. By the definition of $t_0$,
 $$D(\ga\wedge\ga'(r),\ga\wedge\ga'(t))=t-r$$
 and so
 $$D(\ga\wedge\ga'(r),\ga\wedge\ga'(t'))\leq t'-r$$
 for every $t'\in[t_0,t_0+\ve[$. As in the first case above, the last inequality must be
 an equality, which again contradicts the definition of $t_0$.
 
 This completes the proof of Lemma \ref{minilemma}. \cq
 
 \medskip
 \noi{\bf Proof of Proposition \ref{minigeo}:} We need to verify that
 \beq
 \label{minitech4}
 D(\ga_s(r),\ga_s(t))=t-r
 \eeq
 for every $0\leq r< t\leq {\ov Z}_s$. If $r=0$, or if $t={\ov Z}_s$, (\ref{minitech4})
 follows from (\ref{minitech1}) and (\ref{minitech2}). Let
 $\ve>0$ and let us fix $r$ and $t$
 with $0<r<t<{\ov Z}_s$. By the definition of $\ga_s$, we can find 
 $\ga,\ga'\in {\rm G}^*(s)$ such that
 $$\ga_s(r)\leq \ga(r)\leq \ga_s(r)+\ve$$
 and 
 $$\ga_s(t)\leq \ga'(t)\leq\ga_s(t)+\ve.$$
 It follows that
 \ba
 &&\ga_s(r)\leq \ga\wedge\ga'(r)\leq \ga_s(r)+\ve,\\
 &&\ga_s(t)\leq \ga\wedge\ga'(t)\leq \ga_s(t)+\ve.
 \ea
 Choosing $\ve$ small, we see that $D(\ga_s(r),\ga_s(t))$
 can be made arbitrarily close to the quantity $D(\ga\wedge\ga'(r),\ga\wedge\ga'(t))=t-r$,
 by Lemma \ref{minilemma}. This completes the proof of (\ref{minitech4}). \cq
 
 \smallskip
 In view of our applications, it will also be important to
 consider a notion that is dual to the notion of a minimal geodesic. 
 For every $s\in[0,1]$, we  denote by $\ov{\rm G}(s)$ the set of all
 mappings $\ga:[0,{\ov Z}_s]\la[s,1]$ such that $(\pp(\ga(t)),0\leq t\leq {\ov Z}_s)$
 is a geodesic from $\rho$ to $\pp(s)$. We then set 
 $$\ov\gamma_s(t)=\sup\{\gamma(t):\gamma\in \ov{\rm G}(s)\}\;,\;
 \ov\Gamma_s(t)=\pp(\ov\ga_s(t))\;,$$
 for every $t\in[0,\ov Z_s]$. 
 The very same arguments as in the proof of Proposition \ref{minigeo}
 show that $(\ov\Gamma_s(t),0\leq t\leq \ov Z_s)$ is a geodesic from
 $\rho$ to $\pp(s)$, which is called the maximal geodesic. 
 
 The next two sections
 are devoted to a number of technical results about geodesics.
 Although we concentrate on minimal geodesics, it should be noted that
 symmetric arguments yield corresponding results for maximal geodesics.
 
 \section{The main estimate}
 
 In this section we fix $U\!\!\in\,]0,1[$. Recall that $\gamma_U$ is
 the minimal geodesic that was introduced in the previous section. 
 We denote by $\r(\ga_U)$ the range of $\ga_U$:
 $$\r(\ga_U)=\{\ga_U(t):0\leq t\leq {\ov Z}_U\}.$$
 We also fix $u_0\in\,]0,U/2[$. 
 
  \begin{lemma}
 \label{keylemma}
 Let $\delta\!\in\,]0,u_0/2[$ and $\eta >0$.
For every $u\in\,]0,U[$,
let  $B_{\eta,u}$ denote the event
 $$B_{\eta,u}=\Big\{{\ov\eg}_u>\eta\;,\;\inf_{\eta/2\leq t\leq\ov\eg_u} \ov W_u(t) > \eta\Big\}.$$
 There exists a constant $C=C(p,U,u_0,\eta,\delta)$ such that,
 for every $u\in[u_0,U[$, $v\in\,]u,U]$ and $\al >0$,
 \beq
 \label{keybound}
 P\Big[\{\r(\ga_U)\cap]u,v]\not =\varnothing\}\cap B_{\eta,u}\cap \Big\{
 \inf_{u-\delta\leq t\leq u} {\ov\eg}_t>{\ov\eg}_u-\al\Big\}\Big]
 \leq C\,\al\,P[\r(\ga_U)\cap]u,v]\not =\varnothing].
 \eeq
 \end{lemma}
 
 The remaining part of this
 section is devoted to the proof of Lemma \ref{keylemma}. 
 To simplify notation we will write $\ga=\ga_U$ and $\Gamma=\Gamma_U$ in this proof.
 The underlying idea is to apply the Markov property
 in reverse time to the process $({\ov\eg}_r,\ov W_r)_{0\leq r\leq 1}$ at time $u$.
 Notice that the event $B_{\eta,u}$ only depends on $(\ov\eg_u,\ov W_u)$. Assuming
 that the event $\{\r(\ga)\cap]u,v]\not =\varnothing\}$ is also measurable with
 respect to the $\sigma$-field generated by $(\ov\eg_r,\ov W_r)_{u\leq r\leq 1}$,
 we would need to bound the conditional probability of the event
 $$\Big\{ \inf_{u-\delta\leq t\leq u} {\ov\eg}_t>{\ov\eg}_u-\al\Big\}$$
 given that $\sigma$-field. The latter conditional probability is bounded
 above by $C\al$ as would be the case if ${\ov\eg}$ were replaced by
 a linear Brownian path.
 
 The above line of reasoning is not easy to implement, mainly because 
 the required measurability property of the event $\{\r(\ga)\cap]u,v]\not =\varnothing\}$
 seems difficult to establish. Instead, we will use a discrete
 version of the preceding ideas, and we will rely on the convergence
 (\ref{basicconv}) to derive the bound of Lemma \ref{keylemma}
 from our discrete estimates.
 
 The proof of Lemma \ref{keylemma} depends on several intermediate lemmas.
 Before stating the first of these lemmas, we start with a few simple remarks.
 
 We already noticed that the minimal geodesic $(\ga(r),0\leq r\leq {\ov Z}_U)$
 is monotone increasing. Moreover, for every $r\in\,]0,{\ov Z}_U]$ we have $\pp(\ga(r-))=
 \Gamma(r-)=\Gamma(r)=\pp(\ga(r))$ 
 by continuity, and the definition of the minimal geodesic shows
 that we must have $\ga(r-)=\ga(r)$. The mapping $r\to \ga(r)$ is thus left-continuous
 on $]0,U]$. 
  
 Then let $r_1,r_2,\ldots$ be a sequence dense in $]0,\infty[$. We have
 for every $u<v$,
 $$\{\r(\ga)\cap]u,v]\not =\varnothing\}=\lim_{k\to\infty}\uparrow
 \Big(\bigcup_{i=1}^k \{u<\ga(r_i)\leq v\}\Big).$$
 Here and below we make the convention that $\ga(r)=\infty$ if $r>{\ov Z}_U$. 
 
 Thus, if we can prove that a bound analogous to (\ref{keybound}) holds
 when the event $\{\r(\ga)\cap]u,v]\not =\varnothing\}$ is replaced by
 an event of the type
 $$\bigcup_{i=1}^k \{u<\ga(r_i)\leq v\}\},$$
 with a constant $C$ independent of $k$, then a monotone passage to the limit will give us the desired result.
 
 Our first technical lemma will relate events of the form $\{\ga(r)\leq u\}$
 to other events that are more suitable for the discrete approximations that
 we will use to derive our estimates. For every $u\in\,]0,1[$, we denote by
 ${\rm Fr}([0,u])$ the set of all $v\in[u,1]$ such that $\pp(v)=\pp(v')$ for
 some $v'\in[0,u]$. Note that ${\rm Fr}([0,u])$ is a closed subset of $[0,1]$. 
 Furthermore, the ancestral line of 
 $p_{\ov\eg}(u)$ is contained in $p_{\ov\eg}({\rm Fr}([0,u]))$.
 
 In the following three lemmas,
 we fix $u\in\,]0,U[$ and $r>0$. 
 
 \begin{lemma}
 \label{approxevent}
For every $\ve>0$, let $A_\ve(r,u)$ stand for the event
 on which there exist an integer $q\geq 1$ and points $s_0,s_1,\ldots,s_q\in[u,U]$
 such that $s_0=U,\,s_q\in{\rm Fr}([0,u])$ and, for every $i\in\{0,1,\ldots,q-1\}$,
 \begin{eqnarray}
 \label{approxt1}
 &&{\ov Z}_{s_i}\geq r,\\
 \label{approxt2}
 &&D(s_i,s_{i+1})\leq \ve,
 \end{eqnarray}
 and, for every $i\in\{0,1,\ldots,q-2\}$,
 \beq
 \label{approxt3}
 D(s_i,s_{i+1})\leq {\ov Z}_{s_i}-{\ov Z}_{s_{i+1}}+\frac{\ve}{q}.
 \eeq
 Then,
 $$\{\ga(r)\leq u\}\subset \bigcap_{\ve>0} A_{\ve}(r,u)\qquad {a.s.}$$
 \end{lemma}
 
 \rem The inclusion in the conclusion of the lemma can in fact be replaced by
 an equality. See Lemma \ref{approxcc} below.
   
 \smallskip
 \proof Let $\ve>0$ and let us verify that $\{\ga(r)\leq u\}\subset A_\ve(r,u)$ a.s. 
 We set 
 $$r_0=\sup\{t\in[0,{\ov Z}_U]: \ga(t)\leq u\}.$$
  If $r_0={\ov Z}_U$, then
 $\ga({\ov Z}_U)\leq u$ by left-continuity, and since $\pp(U)=\pp(\ga({\ov Z}_U))$
 we get that $U\in{\rm Fr}([0,u])$. In that case, we simply take $q=1$
 and $s_0=s_1=U$, noting that ${\ov Z}_U\geq r$ on the event $\{\ga(r)\leq u\}$
 (recall that we made the convention that $\ga(r)=\infty$ if $r>{\ov Z}_U$).
 
 So we can assume that $r_0<{\ov Z}_U$. 
 We observe that $\ga(r_0+)\in{\rm Fr}([0,u])$. Indeed
 $\ga(r_0)\leq u$, $\ga(r_0+)\geq u$ and $\pp(\ga(r_0))=\pp(\ga(r_0+))$
 by continuity. 
 
 We then set
 $$s_0=U,\,s_1=\ga({\ov Z}_U-\ve),\ldots,\,s_{q-1}=\ga({\ov Z}_U-(q-1)\ve)$$
 where $q$ is the first integer such that ${\ov Z}_U-q\ve\leq r_0$. We finally set
 $s_q=\ga(r_0+)$. Then $s_q\geq u$ and $s_i\in[u,U]$ for every $i\in\{0,1,\ldots,q\}$
 because $\ga$ is nondecreasing.  Conditions (\ref{approxt2}) and (\ref{approxt3}) clearly
 hold: In (\ref{approxt3}) we have even $D(s_i,s_{i+1})= {\ov Z}_{s_i}-{\ov Z}_{s_{i+1}}$.
 Moreover the condition $\ga(r)\leq u$ implies $r\leq r_0$ and thus ${\ov Z}_{s_i}\geq {\ov Z}_{\ga(r_0)}
 \geq {\ov Z}_{\ga(r)}=r$ for every $i\in\{0,1,\ldots,q\}$. \cq
 
 \smallskip
 We now want to use the convergence (\ref{basicconv}) in order to get 
 discrete approximations of the sets $A_\ve(r,u)$. For every integer $n$,
 and for $0\leq k\leq k'\leq pn$, we denote by $\g^n(k,k')$ the 
 $\sigma$-field generated by the variables $(C^n_i,\Lambda^n_i)$, $k\leq i\leq k'$.
 
 Fix $k\in\{0,1,\ldots,pn\}$.
 If $i,j\in\{k,k+1,\ldots,pn\}$, the distance $d_n(i,j)$ is in general not
 measurable with respect to the $\sigma$-field $\g^n(k,pn)$. However,
 we may define a ``modified distance'' $\wt d^k_n(i,j)$ for which this will
 be true: To this end, we restrict our attention to edges of the map
 $\mm_n$ that are generated between steps $k$ and $pn-1$
 of the Bouttier-Di Francesco-Guitter bijection of subsection 2.6. To give a more precise definition, recall from
 subsection 2.6 the notion of a successor (with respect to the mobile $(\tau_n,
 (\ell^n_v)_{v\in\tau^\circ_n})$), and
 say by convention that the successor of any integer
 $i\in\{0,1,\ldots,pn\}$ such that the corresponding vertex in $\tau^\circ_n$
 has label $1$ is $pn+1$. Then,
 for $i,j\in\{k,k+1,\ldots,pn\}$, $\wt d^k_n(i,j)$ is the minimal integer
 $\ell\geq 0$ for which there exists a sequence $i_0,j_0,i_1,j_1,\ldots,i_\ell,j_\ell$
 of integers in $\{k,k+1,\ldots,pn,pn+1\}$, such that:
 \begin{description}
 \item{$\bullet$}
  $i_0=i$, $j_\ell=j$;
  \item{$\bullet$} for every $0\leq m\leq \ell$, either $i_m\not =pn+1$ and $j_m\not =pn+1$,
  or $i_m=j_m=pn+1$;
  \item{$\bullet$}
 $i_m\sim_{[n]} j_m$ for every $0\leq m\leq \ell$ such that $i_m\not =pn+1$;
 \item{$\bullet$} for every $1\leq m\leq \ell$,
 either $i_m$ is the successor of $j_{m-1}$ or $j_{m-1}$ is the successor of
 $i_m$. 
 \end{description}
 By convention $\wt d^k_n(i,j)=\infty$ if $i\notin \{k,k+1,\ldots,pn\}$
 or $j\notin \{k,k+1,\ldots,pn\}$.
 It is then easy to verify that $\wt d^k_n(i,j)$ is measurable
 with respect to $\g^n(k,pn)$. Furthermore $\wt d^k_n(i,j)\geq d_n(i,j)$.
 
 Recall our notation $p_n$ for the canonical projection from 
 $\{0,1,\ldots,pn\}$ onto $\tau^\circ_n$ ($p_n(i)=v^n_i$ if $v^n_0,v^n_1,\ldots,v^n_{pn}$ denotes the 
 contour sequence of the tree $\tau^\circ_n$). We denote by ${\rm Fr}_n([0,k])$ the set of
 all $j\in\{k,k+1,\ldots,pn\}$ for which at least one of the following 
 two conditions holds:
\begin{description}
\item{(a)} there exists $i\in\{0,1,\ldots,k\}$
 such that $d_n(p_n(i),p_n(j))=0$ or $1$;
 \item{(b)} $\Lambda^n_i\geq \Lambda^n_j$ for every $k\leq i\leq j$.
 \end{description}
One can then verify that, for every $j\in\{k,k+1,\ldots,pn\}$, the event 
 $\{j\in {\rm Fr}_n([0,k])\}$ is measurable with respect to the 
 $\sigma$-field $\g^n(k,pn)$. In particular, the properties
 of the contour function imply that the event $\{p_n(j)\in\{p_n(0),\ldots,p_n(k)\}\}$
 is measurable for the $\sigma$-field generated by $(C^n_i,k\leq i\leq pn)$. 
 The event $\{\exists i\in\{0,\ldots,k\}: d_n(p_n(i),p_n(j))=1\}$ is not 
  $\g^n(k,pn)$-measurable in general, but its union with the event where (b) holds
  is $\g^n(k,pn)$-measurable.

To simplify notation,
 we write $\wt d_n=\wt d^{\lfloor pnu\rfloor}_n$.
 We can now define our discrete approximation of the set
 $A_\ve(r,u)$. We let $A^n_\ve(r,u)$ denote the event on which there
 exists an integer $q_n\geq 1$ and points $s^n_0,\ldots,s^n_{q_n}$
 in $\{\lfloor pnu\rfloor,\ldots,\lfloor pnU\rfloor\}$ such that 
 $s^n_0=\lfloor pnU\rfloor$, $s^n_{q_n}
 \in {\rm Fr}_n([0,\lfloor pnu\rfloor])$ and for every $i\in\{0,1,\ldots,q_n-1\}$,
 \begin{eqnarray}
 \label{approxtd1}
 &&\Lambda^n_{s^n_i}\geq \kappa_p^{-1}\,(r-\ve)n^{1/4},\\
 \label{approxtd2}
 &&\wt d_n(s^n_i,s^n_{i+1})\leq 2\kappa_p^{-1}\,\ve\,n^{1/4},
 \end{eqnarray}
 and, for every $i\in\{0,1,\ldots,q_n-2\}$,
 \beq
 \label{approxtd3}
\wt d_n(s^n_i,s^n_{i+1})\leq \Lambda^n_{s^n_i}-\Lambda^n_{s^n_{i+1}}+\frac{2\kappa_p^{-1}\ve}{q_n}
\,n^{1/4}.
 \eeq
 It is important to observe that the event $A^n_\ve(r,u)$ is measurable with respect
 to the $\sigma$-field $\g^n(\lfloor pnu\rfloor,pn)$.

\begin{lemma}
\label{approxdc}
For every $\ve\in\,]0,r/10[$,
$$A_\ve(r,u)\subset \liminf_{n\to\infty} A^n_\ve(r,u)\qquad\hbox{a.s.}$$
\end{lemma}

\proof On the event $A_\ve(r,u)$, we can choose $s_0,s_1,\ldots,s_q\in[u,U]$ as
specified in the statement of Lemma \ref{approxevent}, so that 
in particular properties (\ref{approxt1}),(\ref{approxt2}),(\ref{approxt3}) hold.
Since $s_q\in{\rm Fr}([0,u])$, we can find $s'_q\in]0,u]$
such that $\pp(s_q)=\pp(s'_q)$.
We set $\ov s^n_i=\lfloor pns_i\rfloor$ for every $i\in\{0,1,\ldots,q-1\}$,
and $\ov s^n_q=\lfloor pn s'_q\rfloor -1 <\lfloor pn u\rfloor$.
Then $\ov s^n_0=\lfloor pnU\rfloor$ and  the convergence (\ref{basicconv}) ensures that a.s. for 
all $n$ sufficiently large, for
every $i\in\{0,1,\ldots,q-1\}$,
\begin{eqnarray}
 \label{approxtd1'}
 &&\Lambda^n_{\ov s^n_i}\geq \kappa_p^{-1}\,(r-\ve)n^{1/4},\\
 \label{approxtd2'}
 &&d_n(\ov s^n_i,\ov s^n_{i+1})\leq 2\kappa_p^{-1}\ve\,n^{1/4},
 \end{eqnarray}
 and, for every $i\in\{0,1,\ldots,q-2\}$,
 \beq
 \label{approxtd3'}
d_n(\ov s^n_i,\ov s^n_{i+1})\leq \Lambda^n_{\ov s^n_i}-\Lambda^n_{\ov s^n_{i+1}}+
\frac{2\kappa_p^{-1}\ve}{q}
\,n^{1/4}.
 \eeq

Consider the first index $i_0\in\{0,1,\ldots,q-1\}$ such that $d_n(\ov s^n_{i_0},\ov s^n_{i_0+1})
< \wt d_n(\ov s^n_{i_0},\ov s^n_{i_0+1})$. This index exists since 
$\wt d_n(\ov s^n_{q-1},\ov s^n_q)=\infty$, because $\ov s^n_q<\lfloor pnu\rfloor$. Then choose a 
discrete geodesic $\omega_n$ from $p_n(\ov s^n_{i_0})$ to $p_n(\ov s^n_{i_0+1})$ 
in the map $\mm_n$. By (\ref{approxtd1'}), (\ref{approxtd2'}) and our assumption
$\ve<r/10$, it is clear that this geodesic does not visit the root of $\mm_n$. 
Furthermore, at least one
point on the geodesic $\omega_n$ must belong to $\{p_n(0),\ldots,p_n(\lfloor pnu \rfloor)\}$
(otherwise we would have $d_n(\ov s^n_{i_0},\ov s^n_{i_0+1})
=\wt d_n(\ov s^n_{i_0},\ov s^n_{i_0+1})$).

We let $b_n$ be the first point on the geodesic $\omega_n$ that belongs 
to $\{p_n(0),\ldots,p_n(\lfloor pnu \rfloor)\}$. We put $s^n_{i_0+1}=\ov s^n_{i_0}$
if $b_n=p_n(\ov s^n_{i_0})$, and otherwise, we choose $s^n_{i_0+1}
\in\{\lfloor pnu \rfloor,\ldots,pn\}$ such that $p_n(s^n_{i_0+1})$ is the point preceding
$b_n$ on the geodesic $\omega_n$. In both cases it is clear that
$s^n_{i_0+1}\in {\rm Fr}_n([0,\lfloor pnu \rfloor])$.

We also set
$q_n=i_0+1$ and $s^n_i=\ov s^n_i$ for every $i\in\{0,1,\ldots,i_0\}$. Then
we have $\wt d_n(s^n_{i_0},s^n_{i_0+1})=d_n(s^n_{i_0},s^n_{i_0+1})$. Indeed,
this equality is trivial if $b_n=p_n(\ov s^n_{i_0})$, and otherwise
there exists a geodesic from $p_n(s^n_{i_0})$ to $p_n(s^n_{i_0+1})$ such that for
every point $a$ of this geodesic all representatives of $a$
lie in $\{\lfloor pnu \rfloor,\ldots,pn\}$. It then follows from
(\ref{approxtd1'}),(\ref{approxtd2'}),(\ref{approxtd3'}) that the points $s^n_0,s^n_1,\ldots,s^n_{q_n}$
satisfy (\ref{approxtd1}),(\ref{approxtd2}),(\ref{approxtd3}), so that 
$A^n_\ve(r,u)$ holds for $n$ large. This completes the proof. \cq

\begin{lemma}
\label{approxcc}
For every $\ve>0$, set
$$\wt A_\ve(r,u)=\limsup_{n\to\infty} A^n_\ve(r,u).$$
Then,
$$\{\ga(r)\leq u\}= \bigcap_{\ve>0} A_\ve(r,u)= \bigcap_{\ve>0} \wt A_\ve(r,u)\;,\qquad\hbox{a.s.}$$
\end{lemma}

\proof
From Lemma \ref{approxdc}, we have $A_\ve(r,u)\subset\wt A_\ve(r,u)$ a.s. if $\ve < r/10$. Then 
it follows from Lemma \ref{approxevent} that
$$\{\ga(r)\leq u\}\subset \bigcap_{\ve>0} A_\ve(r,u) \subset \bigcap_{\ve>0} \wt A_\ve(r,u)\;,\qquad\hbox{a.s.}$$
Thus we only need to prove that
$$\bigcap_{\ve>0} \wt A_\ve(r,u)\subset \{\ga(r)\leq u\}\;,\qquad\hbox{a.s.}$$

From now on we assume that the event $\wt A_\ve(r,u)$ holds for every
$\ve$ belonging to a (fixed) sequence decreasing to $0$. Fix one value of 
$\ve$ in this sequence. By the definition of the set $\wt A_\ve(r,u)$
we can find, for every $n$ belonging to a (random) sequence converging
to $\infty$, a discrete path $\ga^\ve_{(n)}=(\ga^\ve_{(n)}(i),0\leq i\leq L^\ve_{(n)})$
taking values in $\{0,1,\ldots,pn\}$, such that the following properties hold:
\begin{description}
\item{$\bullet$} $\ga^\ve_{(n)}$ starts from $\lfloor pnU\rfloor$
and ends at a point $z_n\in {\rm Fr}_n([ 0,\lfloor pnu\rfloor])$;
\item{$\bullet$} the length $L^\ve_{(n)}$ is bounded above by $\Lambda^n_{\lfloor pnU\rfloor}
-\Lambda^n_{z_n}+ 6 \kappa_p^{-1}\ve\,n^{1/4}$;
\item{$\bullet$} any point on the path $\ga^\ve_{(n)}$ lies within $d_n$-distance
at most $2\kappa_p^{-1}\ve\,n^{1/4}$ from a point $z$ of $\{\lfloor pnu\rfloor,\ldots,\lfloor pnU\rfloor\}$
such that $\Lambda^n_z \geq \kappa_p^{-1}(r-\ve)n^{1/4}$;
\item{$\bullet$} $d_n(\ga^\ve_{(n)}(i),\ga^\ve_{(n)}(j))\leq |j-i|$ for every $i,j\in\{0,1,\ldots,L^\ve_{(n)}\}$;
\item{$\bullet$} if $y$ and $y'$ are two points that come in that order on the path $\ga^\ve_{(n)}$,
\beq
\label{approxtc1}
d_n(y,y')\leq \Lambda^n_{y}-\Lambda^n_{y'}+10\kappa_p^{-1}\ve\,n^{1/4}.
\eeq
\end{description}
The path $\ga^\ve_{(n)}$ is constructed by choosing the points 
$s^n_0,s^n_1,\ldots,s^n_{q_n}$ as in the definition of $A^n_\ve(r,u)$
and then concatenating discrete geodesics (relative to $d_n$) between $s^n_{i-1}$ and $s^n_i$,
for $1\leq i\leq q_n$. The preceding properties of $\ga^\ve_{(n)}$
follow from the properties stated in the definition of $A^n_{\ve}(r,u)$.
In particular, (\ref{approxtc1}) follows from (\ref{approxtd2}) and (\ref{approxtd3}).

Using (\ref{basicconv}) and extracting a diagonal subsequence, we can assume that along a sequence
$(n'_j)_{j\geq 1}$ of values of $n$ converging to infinity,
$\kappa_p n^{-1/4}L^\ve_{(n)}$ converges to $L^\ve\leq {\ov Z}_U+6\ve$, 
$(pn)^{-1}z_n$ converges to $z_\infty$ and
for every rational $a\in[0,L^\ve[$,
$(pn)^{-1}\ga^\ve_{(n)}(\lfloor \kappa_p^{-1}n^{1/4}a\rfloor)$
converges to a number $\ga^\ve(a)\in[0,1]$. Then, for every
rationals $a,b\in[0,L^\ve[$, the convergence (\ref{basicconv}) gives
$$D(\ga^\ve(a),\ga^\ve(b))=
\lim_{j\to\infty} \kappa_p(n'_j)^{-1/4}\,d_{n'_j}(\ga^\ve_{(n'_j)}(\lfloor \kappa_p^{-1}(n'_j)^{1/4}a\rfloor),
\ga^\ve_{(n'_j)}(\lfloor \kappa_p^{-1}(n'_j)^{1/4}b\rfloor))
\leq |a-b|.$$
Set $\omega^\ve(a)=\pp(\ga^\ve(a))$ for every rational $a\in[0,L^\ve[$. By
the preceding bound, the mapping $a\la \omega^\ve(a)$, from
$[0,L^\ve[\cap \Q$ into $\mm_\infty$, is $1$-Lipschitz.
It can thus be
extended to a $1$-Lipschitz path from $[0,L^\ve]$ into $\mm_\infty$, which we
still denote by $\omega^\ve$. Clearly, $\omega^\ve(0)=\pp(U)$ and
$\omega^\ve(L^\ve)=\pp(z_\infty)$.

Using the fact that $z_n\in {\rm Fr}_n([ 0,\lfloor pnu\rfloor])$, it is not hard to
verify that $z_\infty\in {\rm Fr}([0,u])$, and thus $\omega^\ve(L^\ve)\in
\pp({\rm Fr}([0,u]))$. From the bound $L^\ve_{(n)}\leq \Lambda^n_{\lfloor pnU\rfloor}
-\Lambda^n_{z_n}+ 6 \kappa_p^{-1}\ve\,n^{1/4}$, we also get
$$L^\ve \leq {\ov Z}_U-{\ov Z}_{\omega^\ve(L^\ve)}+6\ve.$$
Moreover (\ref{approxtc1}) gives
$$D(\omega^\ve(t),\omega^\ve(t'))\leq {\ov Z}_{\omega^\ve(t)}-{\ov Z}_{\omega^\ve(t')} +10\ve$$
for every $0\leq t\leq t'\leq L^\ve$. Finally, any point of the path $\omega^\ve$
lies within $D$-distance at most $2\ve$ from a point $x$ of $\pp([u,U])$
such that ${\ov Z}_x\geq r-\ve$. In particular, ${\ov Z}_{\omega^\ve(t)}\geq r-3\ve$
for every $t\in[0,L^\ve]$.

The preceding construction can be made for every $\ve$ belonging to a
sequence decreasing to $0$. Again via a compactness argument, we can
extract a subsequence of values of $\ve$ along which $L^\ve$
converges to $L$ and the $1$-Lipschitz paths $(\omega^\ve(t\wedge L^\ve),0\leq t\leq L)$
converge uniformly to a limiting path $(\omega^0(t),0\leq t\leq L)$. The path 
$\omega^0$ satisfies the following properties:
\begin{description}
\item{$\bullet$} $\omega^0$ is $1$-Lipschitz, starts from $\pp(U)$ and ends at 
a point $\omega^0(L)\in \pp({\rm Fr}([0,u]))$;
\item{$\bullet$} for every $0\leq t\leq t'\leq L$,
\beq
\label{approxtc2}
D(\omega^0(t),\omega^0(t'))\leq {\ov Z}_{\omega^0(t)}-{\ov Z}_{\omega^0(t')};
\eeq
\item{$\bullet$}
for every $t\in[0,L]$, $\omega^0(t)\in\pp([0,U])$ and ${\ov Z}_{\omega^0(t)}\geq r$;
\item{$\bullet$}
$L\leq {\ov Z}_U-{\ov Z}_{\omega^0(L)}$.
\end{description}
By the last property and the triangle inequality,
$$L\leq D(\rho,\pp(U))-D(\rho,\omega^0(L))\leq D(\pp(U),\omega^0(L)).$$
Since $\omega^0$ is $1$-Lipschitz, we also have the reverse inequality
$D(\pp(U),\omega^0(L))\leq L$, and thus
$$D(\pp(U),\omega^0(L))= L,$$
which implies that $\omega^0$ is a geodesic from $\pp(U)$ to $\omega^0(L)$. For
similar reasons, the inequality in (\ref{approxtc2}) must be an equality, and in particular,
$$D(\rho,\pp(U))=D(\rho,\omega^0(L))+ D(\omega^0(L),\pp(U)).$$

Since $\omega^0(L)\in \pp({\rm Fr}([0,u]))$, we can choose $t_1\in[0,u]$ such that $\pp(t_1)=\omega^0(L)$.
Thanks to the property in the last display, we can concatenate $\omega^0$ with the (time-reversed) simple
geodesic $\Phi_{t_1}$, in order to get a geodesic from $\pp(U)$
to $\rho$. Write $\wt\omega$ for the time-reversal of the latter geodesic.
Then $\wt \omega$ is a geodesic from $\rho$ to $\pp(U)$, which takes
values in $\pp([0,U])$. Hence we can find a mapping $\wt\ga:[0,{\ov Z}_U]\la [0,U]$
such that $\wt \omega(t)=\pp(\wt\ga(t))$ for every $t\in[0,{\ov Z}_U]$. We can impose
$\wt\ga({\ov Z}_{\omega^0(L)})=t_1\leq u$. However, by the definition of $\ga=\ga_U$, we have
$$\ga(t)\leq \wt \ga(t),\qquad t\in[0,{\ov Z}_U]$$
and thus in particular $\ga({\ov Z}_{\omega^0(L)})\leq \wt \ga({\ov Z}_{\omega^0(L)})\leq u$. 
On the other hand, we know that ${\ov Z}_{\omega^0(L)}\geq r$, and by the
monotonicity of $\ga$ we conclude that $\ga(r)\leq u$. This completes the
proof of Lemma \ref{approxcc}. \cq

\smallskip
Before we proceed to the proof of Lemma \ref{keylemma}, we need
one more estimate. Recall from subsection 2.7 the notation $(W^n_i)_{0\leq i\leq pn}$
for the ``discrete snake'' associated with the $p$-mobile
$(\tau_n,(\ell^n_v)_{v\in\tau^\circ_n})$. 

\begin{lemma}
\label{excursion}
Let $\delta\in\,]0,u_0/2[$ and $\eta,\eta'>0$. There exists a constant $K=K(p,U,u_0,\eta,\eta',\delta)$
such that, for every $u\in[u_0,U[$ and 
$\al >0 $, we have for all $n$ sufficiently large, 
$$P\Big(\inf_{\lfloor pn(u-\delta)\rfloor\leq i\leq \lfloor pnu\rfloor} C^n_i
\geq C^n_{\lfloor pnu \rfloor} -\al n^{1/2} \,\Big|\,\g^n(\lfloor pnu \rfloor,pn)\Big)
\leq K\,\al$$
on the event 
$$\Big\{C^n_{\lfloor pnu \rfloor}\geq \eta n^{1/2},
\inf_{\eta n^{1/2}/2\leq j \leq C^n_{\lfloor pnu \rfloor}} W^n_{\lfloor pnu \rfloor}(j)
\geq \eta' n^{1/4}\Big\}.$$
\end{lemma}

The proof of Lemma \ref{excursion} is postponed to the Appendix.

\smallskip
\noi{\bf Proof of Lemma \ref{keylemma}:} As we already noticed after the statement
of Lemma \ref{keylemma}, it is enough to establish the bound
\beq 
\label{keyt1}
P\Big[\Big(\bigcup_{i=1}^k \{u<\ga(r_i)\leq v\}\Big) \cap B_{\eta,u}\cap \Big\{ \inf_{u-\delta\leq t\leq u} {\ov\eg}_t
> {\ov\eg}_u-\al\Big\} \Big]
\leq C\,\al\,P\Big[\bigcup_{i=1}^k \{u<\ga(r_i)\leq v\}\Big]
\eeq
where $r_1,\ldots,r_k$ are fixed positive numbers, and the constant $C$ only depends
on $p,U,u_0,\eta$ and $\delta$.

By Lemma \ref{approxcc}, we have for every 
$i\in\{1,\ldots,k\}$,
$$\ind{\{\ga(r_i)\leq u\}} = \lim_{\ve \to 0} \ind{A_\ve(r_i,u)}
=\lim_{\ve\to 0} \ind{\wt A_\ve(r_i,u)}\qquad\hbox{a.s.}$$
and consequently
$$\ind{\{u<\ga(r_i)\leq v\}} = \lim_{\ve \to 0}\Big( \ind{A_\ve(r_i,v)}
- \ind{\wt A_\ve(r_i,v)}\Big)^+\qquad\hbox{a.s.}$$
So we have
$${\bf 1}\Big\{\bigcup_{i=1}^k\{u<\ga(r_i)\leq v\}\Big\} = \lim_{\ve \to 0}\Big( \sup_{1\leq i\leq k}
\Big(\ind{A_\ve(r_i,v)}
- \ind{\wt A_\ve(r_i,u)}\Big)^+\Big)\qquad\hbox{a.s.}$$
and a symmetric argument gives
$${\bf 1}\Big\{\bigcup_{i=1}^k\{u<\ga(r_i)\leq v\}\Big\} = \lim_{\ve \to 0}\Big(\sup_{1\leq i\leq k}
\Big(\ind{\wt A_\ve(r_i,v)}
- \ind{A_\ve(r_i,v)}\Big)^+\Big)\qquad\hbox{a.s.}$$
We will verify that, for every $\ve>0$,
\begin{eqnarray}
\label{keyt2}
&&E\Big[\sup_{1\leq i\leq k}
\Big(\ind{A_\ve(r_i,v)}
- \ind{\wt A_\ve(r_i,u)}\Big)^+ {\bf 1}\Big\{\inf_{u-\delta\leq t\leq u}{\ov\eg}_t>{\ov\eg}_u-\al\Big\}\,
\ind{B_{\eta,u}}\Big]
\nonumber\\
&&\quad \leq C\al\,E\Big[\sup_{1\leq i\leq k}
\Big(\ind{\wt A_\ve(r_i,v)}
- \ind{A_\ve(r_i,u)}\Big)^+\Big],
\end{eqnarray}
where $C=\lambda_p^{-1}K(p,U,u_0,\lambda_p^{-1}\eta,\kappa_p^{-1}\eta,\delta)$ 
with the notation of Lemma \ref{excursion}.
By passing to the limit $\ve\to 0$, we see that (\ref{keyt1}) follows from (\ref{keyt2}).

By Lemmas \ref{approxdc} and \ref{approxcc}, we have a.s. for $i\in\{1,\ldots,k\}$,
$$
\Big(\ind{A_\ve(r_i,v)}
- \ind{\wt A_\ve(r_i,u)}\Big)^+\leq\Big(\liminf_{n\to\infty}
\ind{A^n_\ve(r_i,v)} -\limsup_{n\to\infty} \ind{A^n_\ve(r_i,u)}\Big)^+
\leq \liminf_{n\to\infty}
\Big(\ind{A^n_\ve(r_i,v)} -\ind{A^n_\ve(r_i,u)}\Big)^+.
$$
Therefore,
$$\sup_{1\leq i\leq k}
\Big(\ind{A_\ve(r_i,v)}
- \ind{\wt A_\ve(r_i,u)}\Big)^+ 
\leq \liminf_{n\to\infty}\Big(\sup_{1\leq i\leq k}
\Big(\ind{A^n_\ve(r_i,v)} -\ind{A^n_\ve(r_i,u)}\Big)^+\Big),\qquad{\rm a.s.}
$$
Then, using (\ref{basicconv}), (\ref{snakeconv}) and Fatou's lemma, the left-hand side of (\ref{keyt2})
is bounded above by
$$\liminf_{n\to\infty}
E\Big[\sup_{1\leq i\leq k}
\Big(\ind{A^n_\ve(r_i,v)}
- \ind{ A^n_\ve(r_i,u)}\Big)^+ 
{\bf 1}\Big\{\inf_{\lfloor pn(u-\delta)\rfloor\leq i\leq \lfloor pnu\rfloor} C^n_i
\geq C^n_{\lfloor pnu \rfloor} -\lambda_p^{-1}\al n^{1/2}\Big\}\,\ind{B^n_{\eta,u}}\Big],$$
where 
$$B^n_{\eta,u}
=\Big\{C^n_{\lfloor pnu \rfloor}\geq \lambda_p^{-1}\eta n^{1/2},
\inf_{\lambda_p^{-1}\eta n^{1/2}/2\leq j \leq C^n_{\lfloor pnu \rfloor}} W^n_{\lfloor pnu \rfloor}(j)
\geq \kappa_p^{-1} \eta n^{1/4}\Big\}.$$
Now note that both $B^n_{\eta,u}$ and the variable
$$\sup_{1\leq i\leq k}
\Big(\ind{A^n_\ve(r_i,v)}
- \ind{A^n_\ve(r_i,u)}\Big)^+ $$
are measurable with respect to $\g^n(\lfloor pnu \rfloor,pn)$. Hence, we can apply Lemma \ref{excursion}
and we get that the left-hand side of (\ref{keyt2}) is bounded above by
$$\lambda_p^{-1}K\,\al\,\liminf_{n\to\infty}
E\Big[\sup_{1\leq i\leq k}
\Big(\ind{A^n_\ve(r_i,v)}
- \ind{A^n_\ve(r_i,u)}\Big)^+ \Big]$$
where $K=K(p,U,u_0,\lambda_p^{-1}\eta,\kappa_p^{-1}\eta,\delta)$.
So in order to get (\ref{keyt2}), it only remains to verify that the latter liminf is bounded above by
$$E\Big[\sup_{1\leq i\leq k}
\Big(\ind{\wt A_\ve(r_i,v)}
- \ind{A_\ve(r_i,u)}\Big)^+\Big].$$
To this end, we just have to note that
\ba
\limsup_{n\to\infty}\Big(
\sup_{1\leq i\leq k}
\Big(\ind{A^n_\ve(r_i,v)}
- \ind{A^n_\ve(r_i,u)}\Big)^+ \Big)&=&
\sup_{1\leq i\leq k} \Big(\limsup_{n\to\infty}
\Big(\ind{A^n_\ve(r_i,v)}
- \ind{A^n_\ve(r_i,u)}\Big)^+ \Big)\\
&\leq& \sup_{1\leq i\leq k} \Big(\limsup_{n\to\infty}
\ind{A^n_\ve(r_i,v)}
- \liminf_{n\to\infty} \ind{A^n_\ve(r_i,u)}\Big)^+ \\
&\leq& \sup_{1\leq i\leq k} \Big(\ind{\wt A_{\ve}(r_i,v)}
-\ind{A_\ve(r_i,u)}\Big)^+.
\ea
This completes the proof of Lemma \ref{keylemma}. \cq

\section{Uniform estimates for the volume of balls}

Our next goal is to get uniform bounds for the volume of small balls in $\mm_\infty$. 
The main ingredient of the proof will be a bound on the moments 
of the quantities
$$\ov{\cal J}([0,\ve])$$
where $\ve >0$, and $\ov{\cal J}$ is the occupation measure of $\ov Z$, which
is the random measure on $[0,\infty[$
defined by
$$\ov{\cal J}(A)=\int_0^1 1_A(\ov Z_s)\,ds$$
for every Borel subset $A$ of $[0,\infty[$. 

We will rely on certain estimates under the infinite excursion
measure $\ov\N_0$ of the conditioned Brownian snake.
Let us briefly recall the definition of $\ov\N_0$. More details can be found
in \cite{LGW}.

Let $C([0,\infty[,\W_0)$ stand for the space of all continuous functions from $\R_+$
into $\W_0$, and denote the canonical process on this space by $(W_s)_{s\geq 0}$.
Also let $\zeta_s$ denote the lifetime of $W_s$, for every $s\geq 0$. For $\omega\in
C([0,\infty[,\W_0)$, set $\sigma(\omega)=\sup\{s\geq 0:\zeta_s(\omega)>0\}$.

Denote by 
$\ov\N_0^{(1)}$ the probability measure on $C([0,\infty[,\W_0)$ which is the law
of the process $(\ov W_{s\wedge 1})_{s\geq 0}$. Note that $\sigma=1$, $\ov\N^{(1)}_0$ a.s.
We can now use scaling transformations to define a probability measure $\ov\N^{(a)}_0$,
for every $a>0$, in the following way. For every $\mu>0$, we set
$$\begin{array}{ll}
\zeta^{(\mu)}_s=\mu^2\,\zeta_{s/\mu^4}\;,&s\geq 0,\\
W^{(\mu)}_s(t)=\mu\,W_{s/\mu^4}(t/\mu^2)\;,\qquad&s\geq 0,\,0\leq t\leq \zeta^{(\mu)}_s\,.
\end{array}
$$
By definition, $\ov\N^{(a)}_0$ is the law of $(W^{(a^{1/4})}_s)_{s\geq 0}$ 
under $\N^{(1)}_0$. 
Clearly, $\sigma=a$, $\ov\N^{(a)}_0$ a.s.

We may now set
\beq
\label{uniball2}
\ov\N_0 = \int_0^\infty \frac{da}{2\sqrt{2\pi a^5}}\,\ov\N^{(a)}_0.
\eeq

\begin{lemma} 
\label{uniball}
For every integer $k\geq 2$, there exists a constant $c_k<\infty$
such that, for every $\ve >0$,
$$\ov\N_0\Big(\Big(\int_0^\sigma ds\,\ind{\{\wh W_s\leq \ve\}}\Big)^k\Big)=c_k\,\ve^{4k-6}.$$
Consequently, for every integer $k\geq 1$ and every $\delta\in\,]0,1]$, there exists a 
constant $c_{k,\delta}<\infty$ such that, for every $\ve >0$,
\beq
\label{uniball1}
E[\ov\j([0,\ve])^k]\leq c_{k,\delta}\,\ve^{4k-\delta}.
\eeq
\end{lemma}

\rem It should be possible to replace $\ve^{4k-\delta}$ by $\ve^{4k}$ in (\ref{uniball1}): See
p.669 of \cite{IM} for the case $k=1$.

\smallskip
\proof Recall the notation $W^{(\mu)}$ introduced before the lemma.
Then the law of $(W^{(\mu)}_{s})_{s\geq 0}$ under $\ov\N_0$
is $\mu^6\ov\N_0$. This is easily seen from the decomposition formula
(\ref{uniball2}).
It follows that, for every $\ve>0$,
\ba
\ov\N_0\Big(\Big(\int_0^\sigma ds\,\ind{\{\wh W_s\leq \ve\}}\Big)^k\Big)
&=&\ve^{-6}\,\ov\N_0\Big(\Big(\int_0^{\sigma^{(\ve)}} ds\,\ind{\{\wh W^{(\ve)}_s\leq \ve\}}\Big)^k\Big)\\
&=&\ve^{-6}\,\ov\N_0\Big(\Big(\int_0^{\ve^4\sigma} ds\,\ind{\{\wh W_{s/\ve^4}\leq 1\}}\Big)^k\Big)\\
&=&\ve^{4k-6}\,\ov\N_0\Big(\Big(\int_0^\sigma ds\,\ind{\{\wh W_s\leq 1\}}\Big)^k\Big).
\ea
So, setting
$$c_k=\ov\N_0\Big(\Big(\int_0^\sigma ds\,\ind{\{\wh W_s\leq 1\}}\Big)^k\Big),$$
we need to verify that $c_k<\infty$ if $k\geq 2$.

To this end we will use Theorem 5.1 in \cite{LGW}. We need to introduce some notation.
A marked tree is a pair $\theta=(\tau,(h_v)_{v\in\tau})$ where $\tau$ is a plane tree
(in the formalism of subsection 2.6), and, for every vertex $v\in\tau$, $h_v$ is a nonnegative
real number, which is interpreted as the length of the branch associated with $v$. We then say that
$\tau$ is the discrete skeleton of $\theta$. The plane tree $\tau$
is called binary if every vertex of $\tau$ has either $0$ or $2$ children, and vertices without
children are called leaves.

The
uniform measure on the set of all (binary) marked trees with $k$ leaves is
defined by
$$\int \Lambda_p(d\theta)\,F(\theta)=\sum_{\tau \in {\bf A}^{\rm bin}_k}
\int_{\R_+^\tau} \prod_{v\in\tau} dh_v\,F(\tau,(h_v)_{v\in\tau})$$
where ${\bf A}_k^{\rm bin}$ is the set of all binary plane trees with $k$
leaves. Let $r\geq 0$ and let $\theta=(\tau,(h_v)_{v\in\tau})$ be a marked tree.
We can combine the branching structure of $\theta$ with spatial displacements 
given by $9$-dimensional Bessel processes, to get random variables 
$\ov V_a,a\in\tau$, which are defined under a probability measure $Q^\theta_r$
and constructed as follows. First consider a $9$-dimensional Bessel process 
$R^\varnothing=(R^{\varnothing}_t)_{t\geq 0}$
started from $r$, and set $\ov V_\varnothing = R^\varnothing_{h_\varnothing}$. Then,
conditionally given $R^\varnothing$, consider two other independent
$9$-dimensional Bessel processes $R^1$ and $R^2$, both started from
$\ov V_\varnothing$, and set $\ov V_1=R^1_{h_1}$ and $\ov V_2=R^2_{h_2}$.
We can easily continue the construction by induction. See Section 5 of \cite{LGW}
for more
details.

As a direct consequence of Theorem 5.1 in \cite{LGW}, we have $c_k=2^{k-1}k!\,\wt c_k$, where
$$\wt c_k=\int \Lambda_k(d\theta)\,Q^\theta_0\Big(\prod_{a\in I(\theta)} (\ov V_a)^4
\prod_{a\in L(\theta)} ((\ov V_a) ^{-4}\,\ind{\{\ov V_a\leq 1\}})\Big).$$
Here $L(\theta)$ stands for the set of all leaves of the discrete
skeleton $\tau$ of $\theta$, and $I(\theta)$ is the set of all other vertices of $\tau$. If
$k\geq 2$, we can decompose the binary tree $\theta$ at its first node to obtain
$$\wt c_k=E^{(9)}_0\Big[\int_0^\infty dt\,R_t^4\sum_{j=1}^{k-1} \wt c_j(R_t)\wt c_{k-j}(R_t)\Big]$$
where $R$ denotes a $9$-dimensional Bessel process that starts from $r$ under the 
probability measure $P^{(9)}_r$, and for every $j\geq 1$ and $r\geq 0$,
$$\wt c_j(r)=\int \Lambda_j(d\theta)\,Q^\theta_r\Big(\prod_{a\in I(\theta)} (\ov V_a)^4
\prod_{a\in L(\theta)} ((\ov V_a) ^{-4}\,\ind{\{\ov V_a\leq 1\}})\Big).$$
Let us prove by induction that for every integer $j\geq 1$ there exists
a constant $M_j$ such that, for every $r\geq 0$,
\beq
\label{uniball3}
\wt c_j(r)\leq\left\{\begin{array}{ll}
M_j(r^{-2}\wedge r^{-7})\qquad&{\rm if}\ j=1,\\
M_j(1\wedge r^{-7})&{\rm if}\ j\geq 2.
\end{array}
\right.
\eeq
If $j=1$,
$$\wt c_1(r)=E^{(9)}_r\Big[\int_0^\infty dt\,R_t^{-4}\,\ind{\{R_t\leq 1\}}\Big]
=\int_{\R^9}dz\,G_9(y_r,z)|z|^{-4}\ind{\{|z|\leq 1\}}$$
where $G_9(y,z)=\alpha |z-y|^{-7}$ is the Green function of 
$9$-dimensional Brownian motion, and $y_r$ denotes an
arbitrary point in $\R^9$ such that $|y_r|=r$. From the preceding explicit formula,
straightforward estimates give the bound (\ref{uniball3}) when $j=1$. Let 
$\ell\geq 2$ and assume that the bound (\ref{uniball3}) holds for 
$j=1,\ldots,\ell-1$. Using again a decomposition at the first node, we get
\ba
\wt c_\ell(r)&=&E^{(9)}_r\Big[\int_0^\infty dt\,R_t^4\sum_{j=1}^{\ell-1} \wt c_j(R_t)
\wt c_{\ell-j}(R_t)\Big]\\
&\leq& \sum_{j=1}^{\ell-1} M_jM_{\ell-j}\,E^{(9)}_r\Big[\int_0^\infty dt\,R_t^4 \,
(R_t^{-2}\wedge R_t^{-7})^2\Big]\\
&=&\wt M_\ell\, E^{(9)}_r\Big[\int_0^\infty dt\,(1\wedge R_t^{-10})\Big]
\ea
for some constant $\wt M_\ell$. We have then
$$E^{(9)}_r\Big[\int_0^\infty dt\,(1\wedge R_t^{-10})\Big]
=\int_{\R^9} dz\,G_9(y_r,z)(1\wedge |z|^{-10})
=\alpha \int_{\R^9} dz\,|z-y_r|^{-7}(1\wedge |z|^{-10}).$$
Simple estimates show that the last integral is bounded above by
a constant times $1\wedge r^{-7}$, which gives (\ref{uniball3}) for $j=\ell$.

If $k\geq 2$, then $\wt c_k=\wt c_k(0)<\infty$, which completes the proof of the first
assertion of Lemma \ref{uniball}.

The bound (\ref{uniball1}) now comes as an easy consequence. Using the decomposition
(\ref{uniball2}) and a scaling argument, we get 
$$\ov\N^{(1)}_0\Big(\Big(\int_0^1 ds\,\ind{\{\wh W_s\leq \ve\}}\Big)^k\Big)\leq c'_k\,\ve^{4k-6}.$$
with another constant $c'_k$. Since $\ov\j([0,\ve])$ has the distribution 
of $\int_0^1 ds\,\ind{\{\wh W_s\leq \ve\}}$ under $\ov \N^{(1)}_0$, 
it follows that
$$E[\ov\j([0,\ve])^k]\leq c'_k\,\ve^{4k-6}.$$
This holds for every integer $k\geq 2$, so that a simple application
of H\"older's inequality yields (\ref{uniball1}). \cq

\smallskip
If $x\in\mm_\infty$ and $\ve >0$, we denote by $B_D(x,\ve)$ the closed ball
of radius $\ve$ centered at $x$ in the metric space $(\mm_\infty,D)$. Recall our notation
$\lambda$ for the volume measure on $\mm_\infty$
(cf subsection 2.5).

\begin{corollary}
\label{unibound}
Let $\delta\in\,]0,1]$, and
$$S=\sup_{\ve >0} \Big( \sup_{x\in\mm_\infty } \frac{\lm(B_D(x,\ve))}{\ve^{4-\delta}}\Big).$$
Then $E[S^k]<\infty$ for every integer $k\geq 1$. 
\end{corollary}

\proof Let $k\geq 1$ and $\ve>0$. The same argument as in the proof of Lemma 6.2 in \cite{IM}
gives the bound
$$E\Big[\int_{\mm_\infty} \lm(dx)\,\lm(B_D(x,\ve))^k\Big]
\leq E[\lambda(B_D(\rho,\ve))^k]= E[\ov\j([0,\ve])^k] $$
(the inequality is in fact an equality, see Theorem \ref{rerootinv} below).
Using Lemma \ref{uniball}, we thus get
$$E\Big[\int_{\mm_\infty} \lm(dx)\,\lm(B_D(x,\ve))^k\Big]
\leq c_{k,1}\,\ve^{4k-1}.$$
Let $r>0$. If we suppose
that there exists a point $x_0\in\mm_\infty$ such that $\lm(B_D(x_0,\ve))\geq r\,\ve^{4-\delta}$,
then we have for every $x\in B_D(x_0,\ve)$, 
$$\lm(B_D(x,2\ve))\geq \lm(B_D(x_0,\ve))\geq r\,\ve^{4-\delta},$$
and thus
$$\int_{\mm_\infty} \lm(dx)\,\lm(B_D(x,2\ve))^k\geq (r\,\ve^{4-\delta})^{k+1}.$$
It follows that
\begin{eqnarray}
\label{unibound1}
P(\exists x_0\in\mm_\infty : \lm(B_D(x_0,\ve))\geq r\,\ve^{4-\delta})
&\leq& (r\,\ve^{4-\delta})^{-k-1}\,E\Big[\int_{\mm_\infty} \lm(dx)\,\lm(B_D(x,2\ve))^k\Big]\nonumber\\
&\leq& c_{k,1} \,(r\,\ve^{4-\delta})^{-k-1}\,(2\ve)^{4k-1}\nonumber\\
&=& 2^{4k-1} c_{k,1}\,r^{-k-1}\,\ve^{(k+1)\delta-5}.
\end{eqnarray}

We apply this bound with $\ve=2^{-j}$, $j\in\Z_+$, and $k$ chosen so large that 
$(k+1)\delta\geq 6$. Setting 
$$S_j= \sup_{x\in\mm_\infty } \frac{\lm(B_D(x,2^{-j}))}{(2^{-j})^{4-\delta}},$$
we get from (\ref{unibound1}) that, for every $r>0$,
$$P(S_j>r)\leq \ov c_k\,r^{-k-1}\,2^{-j}.$$
where $\ov c_k=2^{4k-1} c_{k,1}$. It follows that, for every $r>1$,
$$P(S>r)\leq \sum_{j=0}^\infty P(S_j>2^{-(4-\delta)}r)\leq 2^{(4-\delta)(k+1)}\ov c_k\,r^{-k-1}.$$
Since this holds for every (large enough) integer $k$, the result of the corollary readily follows. \cq

\smallskip

We will now state and prove the key proposition that 
motivated the technical lemmas of this section and of the previous one. 
We denote by $LI$ the set of all 
left-increase times of
${\ov\eg}$, and by $RI$ the set of all its right-increase times. Note that $p_{\ov\eg}(LI)
=p_{\ov\eg}(RI)={\rm Sk}(\tree)\cup\{\rho\}$. If $U\in\,]0,1[$ is fixed, then
Lemma \ref{leaf} shows that $U\notin LI\cup RI$ a.s.

Recall from subsection 4.2 the notation $\ga_U$ for the minimal
geodesic associated with $U\in[0,1]$, and $\ov\ga_U$ for the
maximal geodesic. Also recall the notation
$\r(\ga_U)$ for the range of $\ga_U$. As usual $\ov{\r(\ga_U)}$
denotes the closure of $\r(\ga_U)$. 

\begin{proposition}
\label{leftinc}
For every fixed $U\in\,]0,1[$, we have $\ov{\r(\ga_U)}\cap LI =\varnothing$ 
and $\ov{\r(\ov\ga_U)}\cap RI =\varnothing$ almost surely.
\end{proposition}

\proof Fix $k\in\N$,
and denote by $N_k$ the number of those intervals of the form 
$]i2^{-k},(i+1)2^{-k}]$, $0\leq i\leq 2^k-1$ that intersect $\r(\ga_U)$. 

We also fix $u_0\in\,]0,U/2[$, $\eta>0$ and $\delta\in\,]0,u_0/2[$. Let $\beta\in\,]0,\frac{1}{16}[$,
and choose $\xi\in\,]0,\frac{1}{2}[$ such that
\beq
\label{left1}
\xi>1-(\frac{1}{4}-\beta)(3-\beta).
\eeq
We denote by $N^*_k$ the number of those intervals 
$]i2^{-k},(i+1)2^{-k}]$ that intersect $\r(\ga_U)$
and satisfy the following additional properties:
\begin{description}
\item{\rm(a)} $i2^{-k}\geq u_0$;
\item{\rm(b)} $\displaystyle{\inf_{i2^{-k}-\delta\leq t\leq i2^{-k}}{\ov\eg}_t>{\ov\eg}_{i2^{-k}}-(2^{-k})^\xi}$;
\item{\rm(c)} ${\ov\eg}_{i2^{-k}}>\eta$ and
$\displaystyle{\inf_{\eta/2\leq t\leq {\ov\eg}_{i2^{-k}}} \ov W_{i2^{-k}}(t) > \eta}$.
\end{description}
Applying Lemma \ref{keylemma} with $u=i2^{-k}$ and $v=(i+1)2^{-k}\wedge U$, for every
integer $i$ such that $u_0\leq i2^{-k}<U$,
we get the bound
\beq
\label{left2}
E[N^*_k]\leq C (2^{-k})^\xi\,E[N_k]
\eeq
with a 
constant $C$ that depends only on $p,U,u_0,\eta$ and $\delta$. 

We now need to bound $E[N_k]$. For every $\ve>0$, set
$$\Omega(\ve)=\sup_{r,t\in[0,1],\, |r-t|\leq\ve} D(r,t).$$
Recall the bound
\begin{equation}
\label{boundD2}
D(r,t)\leq {\ov Z}_r+{\ov Z}_t - 2\inf_{r\wedge t\leq u\leq r\vee t} {\ov Z}_u
\end{equation}
and also recall the construction of $\ov Z$ by shifting the process $Z^\eg$,
which conditionally given $\eg$ has the distribution of the centered
Gaussian process with covariance function $m_\eg(r,t)$. Using the
bound (\ref{boundD2})
together with standard chaining arguments, one easily gets the
estimate
\beq
\label{left3}
P(\Omega(\ve)>\ve^{\frac{1}{4}-\beta} )=o(\ve)
\eeq
as $\ve\to 0$ (compare with Lemma 5.1 in \cite{IM}). 

Set $\ve_k=2^{-k(\frac{1}{4}-\beta)}$ to simplify notation
and write 
$$\r_{\ve_k}(\ga_U)=\{s\in[0,1]: D(r,s)\leq \ve_k\hbox{ for some }r\in\r(\ga_U)\}$$ 
for the tubular neighborhood
of radius $\ve_k$ of $\r(\ga_U)$, with respect to the pseudo-metric $D$. 
Then, on the event $\{\Omega(2^{-k})\leq \ve_k\}$, any interval 
of the form $]i2^{-k},(i+1)2^{-k}]$ that intersects $\r(\ga_U)$
is contained in $\r_{\ve_k}(\ga_U)$, and so we have on the same event
$$\lm(\pp({\r_{\ve_k}(\ga_U)})) \geq N_k\,2^{-k}.$$
Thus
$$E[N_k]
\leq 2^k\,E[\lm(\pp({\r_{\ve_k}(\ga_U)}))]+2^k\,P[\Omega(2^{-k})>\ve_k].$$
By (\ref{left3}), the second term of the sum is $o(1)$ as $k\to\infty$. 
In order to bound the first term,
we cut the geodesic $\ga_U$ into slices of length $\ve_k$: Precisely
we observe that $\pp(\r_{\ve_k}(\ga_U))$ is contained in the union of the balls
$B_D(\pp(\ga_U(i\ve_k)),2\ve_k)$ for $i\in\{0,1,\ldots,\lfloor \ve_k^{-1}\ov Z_U\rfloor\}$.
It follows that
\ba
E[\lm(\pp({\r_{\ve_k}(\ga_U)}))]
&\leq& E\Big[(\ve_k^{-1}{\ov Z}_U+1)\,\sup_{x\in\mm_\infty}
\lm(B_D(x,2\ve_k))\Big]\\
&\leq&\ve_k^{-1}\, E[({\ov Z}_U+1)^2]^{1/2}\,
E\Big[\sup_{x\in\mm_\infty}
\lm(B_D(x,2\ve_k))^2\Big]^{1/2}\\
&=&O(\ve_k^{3-\beta})
\ea
as $k\to\infty$, by Corollary \ref{unibound}. Finally, we have
$$E[N_k]=O(2^{k(1-(\frac{1}{4}-\beta)(3-\beta))})$$
and, using (\ref{left1}) and (\ref{left2}), it follows that
$$E[N^*_k]\build{\la}_{k\to\infty}^{} 0.$$
Thus $P(N^*_k\geq 1)\la 0$ as $k\to \infty$. 

To complete the proof, suppose that there exists
a point $t_0\in\ov{\r(\ga_U)}\cap LI$, and suppose in
addition that we have the following properties:
\begin{description}
\item{\rm(a)'} $t_0>u_0$;
\item{\rm(b)'} ${\ov\eg}_t\geq {\ov\eg}_{t_0}$ for every $t\in\,]t_0-2\delta,t_0[$;
\item{\rm(c)'} ${\ov\eg}_{t_0}>\eta$ and $\displaystyle{
\inf_{\eta/2 \leq r\leq {\ov\eg}_{t_0}} \ov W_{t_0}(r)>\eta}$.
\end{description}
By Lemma \ref{leaf}, we must have $t_0<U$ (outside a 
set of zero probability depending on $U$).
For every $k\geq 1$, let $i(k)$ be the unique index such that
$t_0\in\,]i(k)2^{-k},(i(k)+1)2^{-k}]$. If $k$ is large enough,
both indices $i=i(k)$ and $i=i(k)+1$ will satisfy the properties (a)--(c) listed above
(for (b), we use (b)' together with the fact that the function 
$t\la{\ov\eg}_t$ is H\"older continuous with 
exponent $\frac{1}{2}-\ve$, for every $\ve>0$). Moreover,
$\r(\ga_U)$ intersects $]i(k)2^{-k},(i(k)+1)2^{-k}]$, or possibly 
$](i(k)+1)2^{-k},(i(k)+2)2^{-k}]$ in the case when 
$t_0=(i(k)+1)2^{-k}$, and so we conclude that
$N^*_k\geq 1$ for $k$ large enough. Since $P(N^*_k\geq 1)$ tends
to $0$ as $k\to\infty$, we get that, with probability one, there
exists no point $t_0 \in\ov{\r(\ga_U)}\cap LI$ that satisfies the 
additional properties (a)'--(c)'. We apply this result to rationals
$u_0,\eta$ and $\delta$ that can be made arbitrarily small
and we obtain the first assertion of the proposition. 

The case  of the maximal geodesic is treated in exactly the same manner.
All lemmas of Section 5 can be extended to
cover this case, with some minor changes due to the lack of symmetry in
the construction of the map $M_n$ from the pair $(C^n,\Lambda^n)$. 
We leave details to the reader. \cq

\section{The main results}

For a fixed $U\!\in\,]0,1[$, we have defined 
the simple geodesic $\Phi_U$, the minimal geodesic 
$\Gamma_U$
and the maximal geodesic $\ov\Gamma_U$. The next proposition will
imply that these three geodesics coincide a.s.

\begin{proposition}
\label{keyprop}
Let $U\!\in\,]0,1[$. 
\begin{description}
\item{\rm(i)}
Almost surely, for
every mapping $\vf:[0,{\ov Z}_U]\la[0,U]$ such that $(\pp(\vf(t)),0\leq t\leq {\ov Z}_U)$
is a geodesic from $\rho$ to $\pp(U)$, we have
$$\pp(\vf(t))=\Gamma_U(t)$$
for every $t\in[0,{\ov Z}_U]$. In particular, $\Phi_U=\Gamma_U$ a.s.
\item{\rm(ii)} Almost surely, for
every mapping $\vf:[0,{\ov Z}_U]\la[U,1]$ such that $(\pp(\vf(t)),0\leq t\leq {\ov Z}_U)$
is a geodesic from $\rho$ to $\pp(U)$, we have
$$\pp(\vf(t))=\ov\Gamma_U(t)$$
for every $t\in[0,{\ov Z}_U]$.
In particular,
$\Phi_U=\ov\Gamma_U$ a.s.
\end{description}
\end{proposition}

\proof Let us prove part (i) of the proposition. By the definition of the minimal geodesic, we have 
$\ga_U(t)\leq \vf(t)$
for every $t\in[0,{\ov Z}_U]$. We argue by contradiction and assume that there
exists $r_0\in\,]0,{\ov Z}_U[$ such that $\ga_U(r_0)<\vf(r_0)$ and $\pp(\ga_U(r_0))\not =\pp(\vf(r_0))$.
Then we can find $s\in\,]\ga_U(r_0),\vf(r_0)[$ such that
$${\ov Z}_s<r_0={\ov Z}_{\ga_U(r_0)}={\ov Z}_{\vf(r_0)}\,,$$
because otherwise we would have $\ga_U(r_0)\approx \vf(r_0)$. Since the mapping $t\la\ga_U(t)$
is nondecreasing and left-continuous, there exists $r_1\in[r_0,{\ov Z}_U]$ such that
$$\begin{array}{ll}
\ga_U(t)>s\qquad &{\rm if}\ t\in\,]r_1,{\ov Z}_U],\\
\ga_U(t)\leq s\qquad &{\rm if}\ t\in[0,r_1].
\end{array}
$$
Note that $r_1={\ov Z}_U$ does not occur a.s., because this would imply
$\ga_U({\ov Z}_U)\leq s< U$ together with $\pp(U)=\pp(\ga_U({\ov Z}_U))$ and 
${\ov Z}_s<{\ov Z}_U$: This is only possible if
$U\sim \ga_U({\ov Z}_U)$, which contradicts Lemma \ref{leaf}.

Thus $r_1<\ov Z_U$, and the fact that $\ov Z_s<r_0\leq r_1=\ov Z_{\ga(r_1)}=\ov Z_{\ga(r_1+)}$
implies that we have the strict inequalities $\ga_U(r_1)< s <\ga_U(r_1+)$. 
From these inequalities it also follows that $\ga_U(r_1)\approx \ga_U(r_1+)$ does not hold.
However we have $\pp(\ga_U(r_1))=\pp(\ga_U(r_1+))$ by continuity. Hence, we must  have $\ga_U(r_1)\sim \ga_U(r_1+)$, and it follows that $p_{\ov\eg}(\ga_U(r_1))=p_{{\ov\eg}}(\ga_U(r_1+))$
is an ancestor of $p_{\ov\eg}(s)$ in $\t_{\ov{\eg}}$. This means in particular
that $\ga_U(r_1+)$ is a left-increase point of ${\ov\eg}$, which is impossible
by Proposition \ref{leftinc}. This contradiction completes the proof of the first statement 
in (i). By applying this statement to $\vf=\vf_U$, we get $\Phi_U=\Gamma_U$ a.s.

Part (ii) of the proposition is proved in a similar way. In particular, the identity 
$\Phi_U=\ov\Gamma_U$ a.s. is obtained by taking $\vf=\ov\vf_U$.
 \cq

\begin{proposition}
\label{keyprop2}
Part {\rm(i)} of Proposition {\rm \ref{keyprop}} holds simultaneously for all
$U\in\,]0,1[\backslash LI$ outside a single set of zero probability. Similarly, part {\rm(ii)}
of Proposition {\rm \ref{keyprop}} holds simultaneously for all
$U\in\,]0,1[\backslash RI$ outside a single set of zero probability. 
\end{proposition}

\proof The assertions of
Proposition \ref{keyprop} hold simultaneously for every rational number $U\in\,]0,1[$
outside a single set of zero probability. From now on, we argue outside this set,
and we deal only with part (i) of Proposition \ref{keyprop}.

 Set
$$H=\{U\in\,]0,1]:\pp(U)=\pp(U')\hbox{ for some }U'<U\}$$
so that $LI\subset H$ in particular. If $U\in H\backslash LI$, then 
there exists $U'<U$ such that $U'\approx U$, and 
$U'\in\,]0,1[\backslash H$. Note that ${\ov Z}_r\geq {\ov Z}_{U'}={\ov Z}_U$ for every $r\in[U',U]$. 
It is then obvious that if $\vf$ is as in the first statement of Proposition \ref{keyprop},
we have $\vf(r)\in[0,U']$ for every $r\in[0,{\ov Z}_U[$. Putting
$\vf'(r)=\vf(r)$ if $r\in[0,{\ov Z}_U[$ and $\vf'({\ov Z}_U)=U'$, we get a function $\vf':[0,{\ov Z}_{U'}]\la[0,U']$
such that $(\pp(\vf'(t)),0\leq t\leq {\ov Z}_{U'})$ is a geodesic from $\rho$
to $\pp(U')$. So if we can prove that $\pp(\vf'(t))=\Gamma_{U'}(t)$
for every $t\in[0,{\ov Z}_{U'}]$, we will get the desired conclusion for $\vf$
since $\pp(\vf(t))=\pp(\vf'(t))$
and $\Gamma_U(t)=\Gamma_{U'}(t)$ for every $t\in[0,{\ov Z}_U]$.

The previous observations show that we may restrict our attention to
the case $U\in\,]0,1[\backslash H$, which we consider from now on.
Let $\ve\in\,]0,{\ov Z}_U[$. Then we can find $s\in\,]0,U[$ such that the following properties hold:
\begin{description}
\item{(a)} ${\ov\eg}_t>{\ov\eg}_s$ for every $t\in\,]s,U]$, and thus $p_{\ov\eg}(s)\in\,\rrbracket \rho,p_{\ov\eg}(U)\llbracket$;
\item{(b)} ${\ov Z}_t\geq {\ov Z}_U-\ve$ for every $t\in[s,U]$;
\item{(c)} ${\displaystyle\min_{a\in\llbracket p_{\ov\eg}(s),p_{\ov\eg}(U)\rrbracket} {\ov Z}_a< {\ov Z}_s.}$
\end{description}
The reason why we can impose condition (c) comes from the fact that the
mapping $a\to \ov Z_a$ cannot be monotone on a line segment of the tree
$\tree$. 

Let $\vf$ be as in the first statement of Proposition \ref{keyprop}. Since $U\notin H$ and 
$\vf$ takes values in $[0,U]$, it is clear that $\vf({\ov Z}_U)=U$. Set
$$t^*=\sup\{t\in[0,{\ov Z}_U]:\vf(t)\leq s\}.$$
We can choose a sequence $(r_k)_{k\geq 1}$ with $0\leq r_k\leq t^*$,
$r_k\la t^*$ as $k\to\infty$ and $\vf(r_k)\leq s$ for every $k\geq 1$. By compactness,
we may assume that $\vf(r_k)\la v\in[0,s]$, and
we have $\pp(\vf(t^*))=\pp(v)$ by continuity. In particular, we must have
$t^*<{\ov Z}_U$, because otherwise it would follow that
$\pp(U)=\pp(\vf({\ov Z}_U))=\pp(v)$, which is impossible 
since $v\leq s<U$ and $U\notin H$. 

Then let $(t_k)_{k\geq 1}$ be a sequence such that $t^*<t_k<{\ov Z}_U$
and $t_k\la t^*$ as $k\to\infty$. By compactness again, we may assume that
$\vf(t_k)\la u\in[s,U]$, and we have $\pp(u)=\pp(\vf(t^*))=\pp(v)$.
We claim that $u>s$. Indeed, if $u=s$, the curve $(\pp(\vf(t)),t^*\leq t\leq {\ov Z}_U)$
starts from $\pp(s)$ and ends at $\pp(U)$, and then (c),
together with the fact that ${\ov Z}_{\pp(\vf(t))}\geq {\ov Z}_{\pp(\vf(t^*))}={\ov Z}_s$
for every $t\in[t^*,{\ov Z}_U]$, gives a contradiction with Lemma \ref{labelsline}.

\bigskip
\begin{center}
\ifx\JPicScale\undefined\def\JPicScale{1}\fi
\unitlength \JPicScale mm
\begin{picture}(88,78)(0,0)
\linethickness{0.3mm}
\multiput(8,3)(0.16,0.12){500}{\line(1,0){0.16}}
\linethickness{0.3mm}
\multiput(43,53)(0.12,-0.48){42}{\line(0,-1){0.48}}
\linethickness{0.3mm}
\multiput(45,45)(0.12,0.21){67}{\line(0,1){0.21}}
\linethickness{0.3mm}
\multiput(15.75,40.25)(0.12,-0.33){75}{\line(0,-1){0.33}}
\linethickness{0.3mm}
\multiput(20,29.38)(0.12,0.21){58}{\line(0,1){0.21}}
\linethickness{0.3mm}
\multiput(12,29)(0.2,-0.12){50}{\line(1,0){0.2}}
\linethickness{0.3mm}
\put(36,24){\line(1,0){29}}
\linethickness{0.3mm}
\multiput(52,24)(0.16,0.12){50}{\line(1,0){0.16}}
\linethickness{0.3mm}
\multiput(45,24)(0.12,-0.15){67}{\line(0,-1){0.15}}
\linethickness{0.3mm}
\multiput(72,51)(0.12,0.46){50}{\line(0,1){0.46}}
\linethickness{0.3mm}
\multiput(64,70)(0.12,-0.13){83}{\line(0,-1){0.13}}
\linethickness{0.3mm}
\multiput(69,65)(0.12,0.36){25}{\line(0,1){0.36}}
\linethickness{0.3mm}
\multiput(48,64)(0.12,-1.5){8}{\line(0,-1){1.5}}
\linethickness{0.3mm}
\multiput(62.38,59)(0.24,-0.12){42}{\line(1,0){0.24}}
\linethickness{0.3mm}
\multiput(63.62,55)(0.63,0.12){8}{\line(1,0){0.63}}
\put(51,52){\makebox(0,0)[cc]{}}

\put(5,6){\makebox(0,0)[cc]{$\rho$}}

\put(36,59.62){\makebox(0,0)[cc]{}}

\put(66,78){\makebox(0,0)[cc]{$p_{\ov\eg}(u)$}}

\put(61,63){\makebox(0,0)[cc]{$p_{\ov\eg}(r)$}}

\put(77.25,59.62){\makebox(0,0)[cc]{}}

\put(88,67){\makebox(0,0)[cc]{$p_{\ov\eg}(U)$}}

\put(62.25,38.38){\makebox(0,0)[cc]{$p_{\ov\eg}(s)$}}

\put(44,68){\makebox(0,0)[cc]{$p_{\ov\eg}(v)$}}

\linethickness{0.3mm}
\put(61,43){\circle*{2}}

\linethickness{0.3mm}
\put(63,59){\circle*{2}}

\linethickness{0.3mm}
\put(72,74){\circle*{2}}

\linethickness{0.3mm}
\put(87,63){\circle*{2}}

\linethickness{0.3mm}
\put(9,4){\circle*{2}}

\linethickness{0.3mm}
\put(48,64){\circle*{2}}

\end{picture}

\end{center}
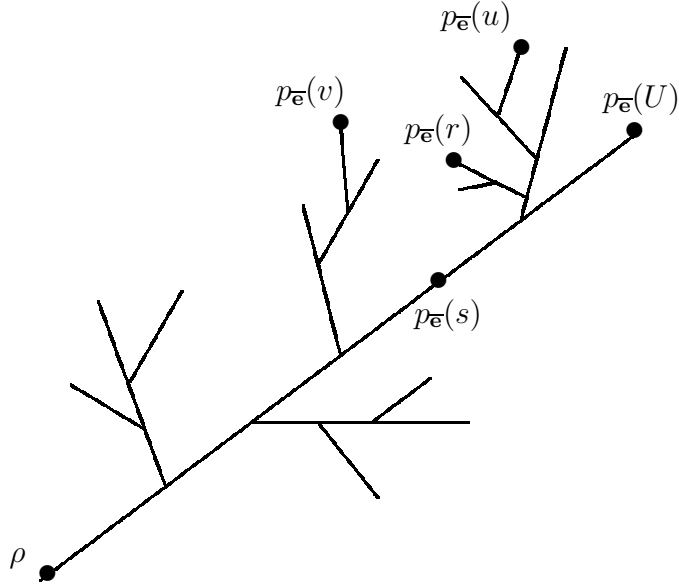
\noindent Figure 4. Illustration of the proof of Proposition \ref{keyprop2}: We have 
$u<s<v$, $u\approx v$ and for every rational $r\in]s,u[$, $\vf(t)=\Phi_r(t)=\Phi_U(t)$
for every $t\leq \ov Z_u=\ov Z_v$.

\medskip

To complete the proof, choose a rational $r\in\,]s,u[$. Since $v<r<u$
and $u\approx v$ ($u\sim v$ is impossible by (a) since $v\leq s<u$), the definition of
the simple geodesic shows that $\Phi_r({\ov Z}_v)=\pp(v)$. 
On the other hand, we have
${\ov Z}_v=t^*$ since $\pp(v)=\pp(\vf(t^*))$, and by setting
$$\wt\vf(t)=\left\{\begin{array}{ll}
\vf_r(t)\qquad&{\rm if}\ t^*\leq t\leq {\ov Z}_r,\\
\vf(t)&{\rm if}\ 0\leq t<t^*,
\end{array}
\right.
$$
we get a map $\wt\vf:[0,{\ov Z}_r]\la [0,r]$ such that
$(\pp(\wt\vf(t)),0\leq t\leq {\ov Z}_r)$ is a geodesic from $\rho$ to $\pp(r)$. 
Since we assumed that the assertions of Proposition \ref{keyprop} hold
for every rational, we obtain that $\pp(\vf(t))=\Phi_r(t)$ for every
$t\in[0,t^*]$. From (b), we get
$t^*={\ov Z}_u\geq {\ov Z}_U-\ve$ and also $\Phi_r(t)=\Phi_U(t)$
for $t<{\ov Z}_U-\ve$. Therefore $\pp(\vf(t))=\Phi_U(t)$
for every $t\in[0,{\ov Z}_U-\ve[$. Since $\ve$ can be taken arbitrarily
small, we have $\pp(\vf(t))=\Phi_U(t)$ for every $t\in[0,\ov Z_U]$. By taking
$\vf=\ga_U$, we also see that $\Gamma_U=\Phi_U$, which completes the proof. \cq

\smallskip

For every $U\!\in\,]0,1[$, we set
$$\l_U=p_{\ov\eg}^{-1}(\llbracket \rho,p_{\ov\eg}(U)\rrbracket)$$
and
$$\r_U=\Big\{s\in[0,U]:{\ov Z}_s=\min_{s\leq r\leq U}{\ov Z}_r\Big\}
\cup \Big\{s\in[U,1]:{\ov Z}_s=\min_{U\leq r\leq s}{\ov Z}_r\Big\}.$$
Then $\l_U$ and $\r_U$ are closed, and $\pp(\l_U)\cap\pp(\r_U)=\{\rho,\pp(U)\}$. 
Indeed, if $x\in\pp(\r_U)\backslash\{\rho,\pp(U)\}$, there exist $s\in\,]0,U[$
and $s'\in\,]U,1[$ such that $s\approx s'$ and  $\pp(s)=\pp(s')=x$.
By Lemma \ref{equirel}, $x$ does not belong to ${\rm Skel}_\infty$, and a fortiori 
$x\notin \pp(\l_U)\backslash\{\rho,\pp(U)\}$. Also notice that 
$\pp(\r_U)$ is the range of the simple geodesic
$\Phi_U$.

We also set
\ba
&&O^1_U=[0,U]\backslash (\r_U\cup \l_U),\\
&&O^2_U=[U,1]\backslash (\r_U\cup \l_U).
\ea
Then $O^1_U$ and $O^2_U$ are open and disjoint, and it is not hard to verify
that $\pp^{-1}(\pp(O^1_U))=O^1_U$ and $\pp^{-1}(\pp(O^2_U))=O^2_U$.
From the definition of the quotient topology, it follows that
$\pp(O^1_U)$ and $\pp(O^2_U)$ are two disjoint open subsets of $\mm_\infty$.

\begin{lemma}
\label{ancestline}
Almost surely, for every $U\!\in\,]0,1[$ and every geodesic $(\omega(t),0\leq t\leq {\ov Z}_U)$
from $\rho$ to $\pp(U)$, the set $\{\omega(t):0<t<{\ov Z}_U\}$ does not intersect $\pp(\l_U)$.
\end{lemma}

\proof We argue by contradiction and assume that there exists $t_0\in\,]0,{\ov Z}_U[$
such that $\omega(t_0)\in\pp(\l_U)$. Then $\omega(t_0)\in\pp(\l_U)\backslash\pp(\r_U)$.
Since $\pp(\r_U)$ is closed, a continuity 
argument shows that $\omega(t)\notin\pp(\r_U)$ if $t>t_0$ is sufficiently close 
to $t_0$. Also Lemma \ref{skeletongeo} implies that we can find values of $t>t_0$
arbitrarily close to $t_0$ such that $\omega(t)\notin\pp(\l_U)$. It follows that
we can find a connected component $]t_1,t_2[$ of the open set
$$\{t\in\,]t_0,{\ov Z}_U[:\omega(t)\notin \pp(\r_U\cup\l_U)\}$$
such that $\omega(t_1)\in\pp(\l_U)$. Since $\omega(]t_1,t_2[)$
is connected, we have either $\omega(]t_1,t_2[)\subset \pp(O^1_U)$
or $\omega(]t_1,t_2[)\subset \pp(O^2_U)$. We assume for definiteness that
$\omega(]t_1,t_2[)\subset \pp(O^1_U)$.

Let $a_1$ be the unique element of $\tree$ such that $\Pi(a_1)=\omega(t_1)$,
and let $r_1<U$ denote the smallest element of $[0,1]$ such that $p_{\ov\eg}(r_1)=a_1$. 
Recall that $\t_{\ov\eg}(a_1)$ denotes the subtree of all descendants of $a_1$.
As a simple consequence of Lemma \ref{labelsline}, we have 
$\omega([t_1,{\ov Z}_U])\subset \Pi(\t_{\ov\eg}(a_1))$.
Indeed, suppose that there exists $t'_1\in\,]t_1,\ov Z_U[$ such that
$\omega(t'_1)\notin \Pi(\t_{\ov\eg}(a_1))$, and let 
$a'_1\in\tree\backslash\tree(a_1)$ be such that $\omega(t'_1)=\Pi(a'_1)$.
Then $a_1\in \llbracket a'_1\tri p_{\ov\eg}(U),p_{\ov\eg}(U)\rrbracket$,
and by applying Lemma \ref{labelsline} to the curve $(\omega(t),t'_1\leq t\leq \ov Z_U)$
we get $\ov Z_{\omega(t_1)}\geq \min\{\ov Z_{\omega(t)}: t'_1\leq t\leq \ov Z_U\}$, which is
a contradiction.

Fix any $t^*\in\,]t_1,t_2[$, and choose $a^*\in\tree(a_1)$ such that 
$\Pi(a^*)=\omega(t^*)$. Set $a_0=a^*\tri p_{\ov\eg}(U)$. Again a simple
application of Lemma \ref{labelsline} shows that $a_0\not = a_1$,
and therefore $a_0\in\,\rrbracket a_1,p_{\ov\eg}(U)\rrbracket$. 
If $r_0$ is the smallest  representative of $a_0$ in $[0,1]$, we have thus $r_1<r_0$.
Set
$$T=\inf\{t\geq t_1: \omega(t)\in \Pi(\tree(a_0))\}.$$
We have $T>t_1$ because $\Pi(\tree(a_0))$ is closed in $\mm_\infty$, and
$\omega(t_1)=\Pi(a_1)\notin \Pi(\tree(a_0))$.
Moreover $T\leq t^*$ because $a^*\in \Pi(\tree(a_0))$. Thus we have $t_1<T\leq t^*<t_2$. 
Furthermore, we can write $\omega(T)=\Pi(b)$ for some $b\in \tree(a_0)$. 
Note that $\omega(T)$ belongs to the boundary of $\Pi(\tree(a_0))$, and $\omega(T)\not =\Pi(a_0)$
(because $\{\omega(t):t_1<t<t_2\}$ does not intersect $\pp(\l_U)$).
As noticed in Section 3, this implies the existence of $b'\in\tree\backslash \tree(a_0)$
such that $b\approx b'$.

\begin{center}
\ifx\JPicScale\undefined\def\JPicScale{1}\fi
\unitlength \JPicScale mm
\begin{picture}(119.38,83.12)(0,0)
\linethickness{0.3mm}
\multiput(10,5)(0.18,0.12){583}{\line(1,0){0.18}}
\put(7.5,8.12){\makebox(0,0)[cc]{$\rho$}}

\put(5,10){\makebox(0,0)[cc]{}}

\linethickness{0.3mm}

\linethickness{0.3mm}
\put(115.31,75.31){\circle*{0.62}}

\put(117.5,78.12){\makebox(0,0)[cc]{$p_{\ov\eg}(U)$}}

\put(119.38,77.5){\makebox(0,0)[cc]{}}

\put(87.5,83.12){\makebox(0,0)[cc]{$p_{\ov\eg}(O^1_U)$}}

\put(106.88,41.25){\makebox(0,0)[cc]{$p_{\ov\eg}(O^2_U)$}}

\linethickness{0.3mm}
\put(70.62,45.62){\line(0,1){26.26}}
\linethickness{0.3mm}
\multiput(59.38,76.88)(0.12,-0.14){94}{\line(0,-1){0.14}}
\linethickness{0.3mm}
\multiput(63.75,71.88)(0.12,0.34){26}{\line(0,1){0.34}}
\linethickness{0.3mm}
\multiput(63.12,80.62)(0.12,-0.23){16}{\line(0,-1){0.23}}
\linethickness{0.3mm}
\multiput(70.62,58.12)(0.12,0.13){94}{\line(0,1){0.13}}
\linethickness{0.3mm}
\multiput(75.62,70.62)(0.13,-0.62){10}{\line(0,-1){0.62}}
\linethickness{0.3mm}
\multiput(61.88,63.75)(0.12,-0.15){73}{\line(0,-1){0.15}}
\linethickness{0.3mm}
\multiput(63.75,56.88)(0.75,-0.13){5}{\line(1,0){0.75}}
\linethickness{0.3mm}
\put(70.62,45.09){\line(0,1){0.44}}

\put(70.62,45.09){\line(0,1){0.44}}

\linethickness{0.3mm}
\put(52.41,33.21){\line(0,1){0.45}}
\put(51.97,33.66){\line(1,0){0.44}}
\put(51.97,33.21){\line(0,1){0.45}}
\put(51.97,33.21){\line(1,0){0.44}}

\linethickness{0.3mm}
\multiput(51.88,50)(0.12,-0.17){68}{\line(0,-1){0.17}}
\linethickness{0.3mm}
\multiput(56.25,43.75)(0.12,0.35){16}{\line(0,1){0.35}}
\linethickness{0.3mm}
\linethickness{0.3mm}
\multiput(45.62,50)(0.28,-0.12){31}{\line(1,0){0.28}}
\linethickness{0.3mm}
\multiput(58.75,61.88)(0.3,-0.12){21}{\line(1,0){0.3}}
\put(53.75,30.62){\makebox(0,0)[cc]{$a_1$}}

\put(74,43){\makebox(0,0)[cc]{$a_0$}}

\put(56.88,65){\makebox(0,0)[cc]{$b$}}

\put(44.38,53.75){\makebox(0,0)[cc]{$b'$}}

\put(58.75,82.5){\makebox(0,0)[cc]{$a_*$}}

\linethickness{0.3mm}
\put(62.59,80.71){\line(0,1){0.45}}
\put(62.59,81.16){\line(1,0){0.44}}
\put(63.03,80.71){\line(0,1){0.45}}
\put(62.59,80.71){\line(1,0){0.44}}

\linethickness{0.3mm}

\linethickness{0.3mm}
\put(58.84,61.97){\line(0,1){0.44}}
\put(58.84,62.41){\line(1,0){0.45}}
\put(59.29,61.97){\line(0,1){0.44}}
\put(58.84,61.97){\line(1,0){0.45}}

\linethickness{0.3mm}
\put(45.31,50.31){\circle*{0.62}}

\linethickness{0.5mm}
\multiput(50,35)(0.19,-0.12){10}{\line(1,0){0.19}}
\linethickness{0.5mm}
\multiput(46.88,38.12)(0.12,-0.16){16}{\line(0,-1){0.16}}
\linethickness{0.5mm}
\multiput(45,42.5)(0.12,-0.31){10}{\line(0,-1){0.31}}
\linethickness{0.5mm}
\put(45,44.38){\line(0,1){3.12}}
\linethickness{0.5mm}
\multiput(46.25,51.88)(0.12,0.12){21}{\line(1,0){0.12}}
\linethickness{0.5mm}
\multiput(50,55.62)(0.15,0.12){21}{\line(1,0){0.15}}
\linethickness{0.5mm}
\multiput(54.38,59.38)(0.19,0.12){16}{\line(1,0){0.19}}
\linethickness{0.5mm}
\multiput(60,63.75)(0.12,0.2){16}{\line(0,1){0.2}}
\linethickness{0.5mm}
\multiput(61.88,71.25)(0.12,-0.5){5}{\line(0,-1){0.5}}
\linethickness{0.3mm}
\put(48.12,41.25){\line(1,0){9.38}}
\linethickness{0.3mm}
\multiput(50,44.38)(0.12,-0.2){16}{\line(0,-1){0.2}}
\linethickness{0.3mm}
\multiput(25.62,38.12)(0.12,-0.13){109}{\line(0,-1){0.13}}
\linethickness{0.3mm}
\multiput(32.5,31.25)(0.12,0.24){31}{\line(0,1){0.24}}
\linethickness{0.3mm}
\multiput(33.12,37.5)(0.11,-0.23){11}{\line(0,-1){0.23}}
\linethickness{0.3mm}
\multiput(25,30.62)(0.48,-0.12){21}{\line(1,0){0.48}}
\linethickness{0.3mm}
\multiput(28.12,32.5)(0.13,-0.31){10}{\line(0,-1){0.31}}
\linethickness{0.3mm}
\multiput(24.38,23.75)(0.12,-0.13){42}{\line(0,-1){0.13}}
\linethickness{0.3mm}
\multiput(26.88,21.25)(0.12,0.44){10}{\line(0,1){0.44}}
\linethickness{0.3mm}
\multiput(35,21.25)(1.81,-0.12){10}{\line(1,0){1.81}}
\linethickness{0.3mm}
\multiput(45,20.62)(0.22,0.12){26}{\line(1,0){0.22}}
\linethickness{0.3mm}
\multiput(49.38,20)(0.12,-0.25){10}{\line(0,-1){0.25}}
\linethickness{0.3mm}
\multiput(63.12,40.62)(0.48,-0.12){26}{\line(1,0){0.48}}
\linethickness{0.3mm}
\multiput(68.75,39.38)(0.13,-1){5}{\line(0,-1){1}}
\linethickness{0.3mm}
\multiput(69.38,36.25)(0.25,-0.12){10}{\line(1,0){0.25}}
\linethickness{0.3mm}
\put(84.38,54.38){\line(1,0){18.74}}
\linethickness{0.3mm}
\multiput(94.38,54.38)(0.12,-0.19){36}{\line(0,-1){0.19}}
\linethickness{0.3mm}
\multiput(97.5,50)(1.06,-0.12){10}{\line(1,0){1.06}}
\linethickness{0.3mm}
\multiput(103.12,49.38)(0.12,0.12){31}{\line(1,0){0.12}}
\linethickness{0.3mm}
\multiput(98.12,54.38)(0.12,0.16){31}{\line(0,1){0.16}}
\linethickness{0.3mm}
\multiput(96.25,62.5)(0.12,0.94){10}{\line(0,1){0.94}}
\linethickness{0.3mm}
\multiput(91.88,71.88)(0.12,-0.12){42}{\line(1,0){0.12}}
\linethickness{0.3mm}
\multiput(96.88,68.75)(0.2,0.12){16}{\line(1,0){0.2}}
\linethickness{0.3mm}
\put(9.69,4.69){\circle*{0.62}}

\linethickness{0.3mm}
\put(71,46){\circle*{2}}

\linethickness{0.3mm}
\put(45,50){\circle*{2}}

\linethickness{0.3mm}
\put(53,34){\circle*{2}}

\linethickness{0.3mm}
\put(10,5){\circle*{2}}

\linethickness{0.3mm}
\put(59,62){\circle*{2}}

\linethickness{0.3mm}
\put(63,81){\circle*{2}}

\linethickness{0.3mm}
\put(116,75){\circle*{2}}

\end{picture}
\end{center}

\vspace{-6mm}
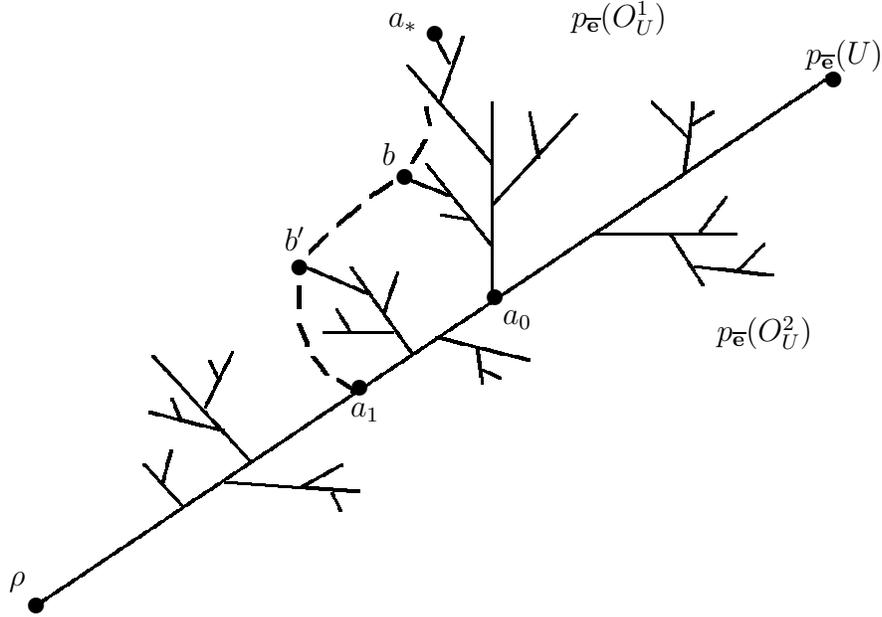
\noindent Figure 4. Illustration of the proof of Lemma \ref{ancestline}: The geodesic $\omega$
visits $a_1$, then stays in the (image under $\Pi$ of the) subtree $\tree(a_1)$. It hits the smaller subtree
$\tree(a_0)$ at time $T$ at the point $\Pi(b)=\Pi(b')$. 

\medskip

Let  $s$ be an arbitrary representative of $b$ in $[0,1]$, so that $\pp(s)=\Pi(b)=\omega(T)$. Since
$b\in\tree(a_0)$ and $b\not =a_0$,
we have $r_0<s$. On the other hand, for every $t\in\,]t_1,T[$, we
know that $\omega(t)\in\mm_\infty\backslash \Pi(\tree(a_0))$, and
$\omega(t)\in \pp(O^1_U)$. It follows that, for every $t\in\,]t_1,T[$, we can find 
$\psi(t)\in[0,r_0]$ such that $\omega(t)=\pp(\psi(t))$. 

We also set $\psi(T)=s$, and for every $t\in[0,t_1]$, we put $\psi(t)=\vf_{r_1}(t)$. The mapping
$(\psi(t),0\leq t\leq T)$ is such that $(\pp(\psi(t)),0\leq t\leq T)$
is a geodesic from $\rho$ to $\omega(T)$, and
$\psi(t)\in[0,\psi(T)]$ for every $t\in[0,T]$. Since $\psi(T)=s\notin LI$
(by Lemma \ref{equirel}, since the equivalence class of $s$ for $\approx$ is not a singleton), 
the first statement of
Proposition \ref{keyprop2} entails that $(\pp(\psi(t)),0\leq t\leq T)$
coincides with the simple geodesic $\Phi_s$.
This is a contradiction, since 
$\pp(\psi(t_1))=\omega(t_1)\in{\rm Skel}_\infty$ and a simple geodesic cannot visit a point
of ${\rm Skel}_\infty$, except possibly at its endpoint.

In the case when $\omega(]t_1,t_2[)\subset \pp(O^2_U)$, the proof is exactly similar, but we
now use the second statement of Proposition \ref{keyprop2}.
 \cq

\smallskip
We now come to one of our main results.

\begin{theorem}
\label{main1}
Almost surely, for every $x\in\mm_\infty\backslash{\rm Skel}_\infty$, there is
a unique geodesic from $\rho$ to $x$, which is the simple geodesic $\Phi_s$
for an arbitrary choice of $s\in[0,1]$ such that $\pp(s)=x$. 
\end{theorem}

\proof
Let $x\in\mm_\infty\backslash {\rm Skel}_\infty$ and let $U\in\,]0,1[$
be such that $\pp(U)=x$ (we exclude the trivial case $x=\rho$). Notice that, 
for every $r\in[0,{\ov Z}_U]$, there is exactly
one point $y$ in $\pp(\r_U)$ such that ${\ov Z}_y=r$, and this point is $y=\Phi_U(r)$.
So in order to prove the theorem, we need only prove that any geodesic from $\rho$
to $x$ takes values in $\pp(\r_U)$. 

So let $(\omega(t),0\leq t\leq {\ov Z}_U)$ be a geodesic from $\rho$
to $x=\pp(U)$. From Lemma \ref{ancestline}, we already
know that $(\omega(t),0<t<{\ov Z}_U)$
does not intersect $\pp(\l_U)$. We argue by contradiction and
suppose that the range of $\omega$ intersects $\mm_\infty\backslash \pp(\r_U)$.
Then let $]t_1,t_2[$ be a connected component of the open set
$$\{s\in\,]0,{\ov Z}_U[:\omega(s)\in \mm_\infty\backslash (\pp(\r_U)\cup \pp(\l_U))\}.$$
Then, $\omega(t_1)\in\pp(\r_U)$, and $\omega(t_2)\in\pp(\r_U)$. 
By the remark of the beginning of the proof, we have
$\omega(t_1)=\Phi_U(t_1)$ and $\omega(t_2)=\Phi_U(t_2)$.
Also, as in the preceding
proof, we have $\omega(]t_1,t_2[)\subset \pp(O^1_U)$ or 
$\omega(]t_1,t_2[)\subset \pp(O^2_U)$, and
we assume for definiteness that $\omega(]t_1,t_2[)\subset \pp(O^1_U)$. For every 
$t\in\,]t_1,t_2[$, let $\psi(t)$ be the smallest representative of $\omega(t)$. Also, let
$\psi(t_2)$ be the largest representative of $\omega(t_2)$. Finally set 
$\psi(t)=\vf_U(t)$ for $t\in[0,t_1]$. Then $\psi(t)\leq U$ for every $t\in[0,t_2[$, whereas
$\psi(t_2)\geq U$. So we see that $\psi$ takes values in $[0,\psi(t_2)]$ and
$(\pp(\psi(t)),0\leq t\leq t_2)$ is a geodesic from $\rho$ to $\pp(\psi(t_2))$. In order to apply
Proposition \ref{keyprop2}, we still need to verify that $\psi(t_2)\notin LI$.
This is clear if $t_2<{\ov Z}_U$ (because the equivalence class of $\psi(t_2)$ for 
$\approx$ has at least two representatives) and if $t_2={\ov Z}_U$
this follows from our assumption $x\notin {\rm Skel}_\infty$. 
Proposition \ref{keyprop2} now shows that 
$\pp(\psi(t))$
must coincide with $\Phi_{\psi(t_2)}(t)$ for every $t\in[0,t_2]$. This is a contradiction
since the range of the simple geodesic $\Phi_{\psi(t_2)}$ is contained in
$\pp(\r_U)$. \cq

\begin{corollary}
\label{typicalgeo}
Almost surely, for $\lambda$-almost every $x\in\mm_\infty$, there is a unique
geodesic from $\rho$ to $x$.
\end{corollary}

\proof By Lemma \ref{leaf}, $\lambda({\rm Skel}_\infty)=0$ a.s., so that the statement
follows from Theorem \ref{main1}. \cq

\smallskip

We then consider the case of vertices belonging
to ${\rm Skel}_\infty$. Such a vertex $x$ can be written uniquely as 
$x=\Pi(a)$ with $a\in{\rm Sk}(\t_{\ov\eg})$, and we denote the multiplicity
of $a$ in $\t_{\ov\eg}$ by $m(x)$. Clearly, $m(x)$ is also the
number of connected components of ${\rm Skel}_\infty\backslash\{x\}$.

\begin{theorem}
\label{main2}
Almost surely, for every $x\in {\rm Skel}_\infty$, there are exactly 
$m(x)$ distinct geodesics from $\rho$ to $x$. These geodesics are the simple geodesics
$\Phi_s$ for all $s\in[0,1]$ such that $\pp(s)=x$.
\end{theorem}

\proof Let $x\in{\rm Skel}_\infty$, and let $(\omega(t),0\leq t\leq {\ov Z}_x)$
be a geodesic from $\rho$ to $x$. By Lemma \ref{skeletongeo}, we can find a sequence
$(r_k)_{k\geq 1}$ in $[0,{\ov Z}_x]$ such that $r_k\la {\ov Z}_x$ as $k\to\infty$
and $\omega(r_k)\in \mm_\infty\backslash{\rm Skel}_\infty$ for every $k$. 
Choose $s_k\in[0,1]$ such that $\pp(s_k)=\omega(r_k)$. 
By extracting a subsequence if necessary, we may assume that 
$s_k\la s\in[0,1]$ as $k\to\infty$, and we have then $\pp(s)=x$.
On the other hand, for every $k\geq 1$,
Theorem \ref{main1} and the fact that $\omega(r_k)\in\mm_\infty
\backslash{\rm Skel}_\infty$ imply that we have 
$\omega(t)=\Phi_{s_k}(t)$ for every $t\in[0,r_k]$.
Passing to the limit $k\to\infty$ we get
$\omega(t)=\Phi_s(t)$ for every 
$t\in[0,{\ov Z}_x[$.

Conversely, for every $s\in[0,1]$ such that $\pp(s)=x$, $\Phi_s$
is a geodesic from $\rho$ to $x$, and the geodesics obtained in this way are distinct.
Indeed, consider the case when $m(x)=2$. Then there exist two reals $s_1,s_2$
such that $0<s_1<s_2<1$,  $p_{\ov\eg}(s_1)=p_{\ov\eg}(s_2)$ and $x=\pp(s_1)=\pp(s_2)$.
By Lemma \ref{equirel}, we must have
$$\min_{s_1\leq s\leq s_2} \ov Z_s <\ov Z_{s_1}=\ov Z_{s_2}.$$
It follows that, if $r>0$ is sufficiently small, 
$\vf_{s_2}({\ov Z}_x-r)$ is a point of $]s_1,s_2[$ which is not equivalent to any
point of $[0,s_1]$. Thus the geodesics $\Phi_{s_1}$
and $\Phi_{s_2}$ are distinct. The same argument applies when
$m(x)=3$. This completes the proof. \cq

\smallskip
A consequence of the previous results is the fact that if $x$ and $x'$
are two points of $\mm_\infty$ distinct from the root, and if 
$\omega$, respectively $\omega'$, is a geodesic going from $\rho$ to $x$, resp.
from $\rho$ to $x'$, then $\omega$ and $\omega'$
must coincide
over a small time interval. We state this property in a slightly more precise
form.

\begin{corollary}
\label{coincid}
Almost surely, for every $\eta>0$, there exists $\alpha\in\,]0,\eta[$ such that the
following holds. Let $x,x'\in\mm_\infty$ such that
$D(\rho,x)\geq\eta$ and $D(\rho,x')\geq \eta$, and let $\omega$, respectively $\omega'$,
be a geodesic from $\rho$ to $x$, resp. from $\rho$ to $x'$. Then,
$\omega(t)=\omega'(t)$ for every $t\in[0,\alpha]$. 
\end{corollary}

\proof
If $\eta>0$ is given, we can choose $\ve>0$ such that $\ov Z_s<\eta$
for every $s\in[0,\ve]\cup[1-\ve,1]$. We then set
$$\alpha=\inf_{s\in[\ve,1-\ve]}\ov Z_s.$$
Notice that $\alpha>0$ since $\ov Z_s>0$ for every $s\in]0,1[$. 

Let $x,x',\omega,\omega'$ be as in the corollary. By the previous two
theorems, we can write $\omega=\Phi_r$ and $\omega'=\Phi_{r'}$,
where $r,r'\in[0,1]$ are such that $x=\pp(r)$ and $x'=\pp(r')$.
In particular, $\ov Z_r=D(\rho,x)\geq \eta$ and $\ov Z_{r'}=D(\rho,x')\geq \eta$.
Our choice of $\ve$ now ensures that both $r$ and $r'$
belong to $[\ve,1-\ve]$. But then, from the definition of the functions
 $\vf_r$ and $\vf_{r'}$, it is immediate that $\vf_r(t)=\vf_{r'}(t)$
for every $t\in[0,\alpha[$. \cq

\smallskip
\rem Our results show that all geodesics starting from the root are
simple geodesics. It follows that a point $x$ of $\mm_\infty$
which is a relative interior point of a geodesic from the root must
be of the form $x=\pp(s)$, for some $s\in[0,1[$ which is a right-increase
time of $\ov Z$. One easily checks that the set of all such points $x$
has Hausdorff dimension $1$, and so is a very small subset of
$\mm_\infty$, which has Hausdorff dimension $4$. This may  be
compared to the work of Zamfirescu \cite{Za}, who proved, in the sense
of Baire's category, that on ``most'' convex surfaces, a typical point is
not an interior point of any geodesic segment.

\section{Gromov-Hausdorff distances and the invariance under re-rooting}

The previous sections were devoted to the study of geodesics 
connecting the root of $\mm_\infty$ to another point. 
Similar results hold if we replace the root $\rho=\pp(0)$ by
a random point chosen according to the volume measure 
$\lambda$. This follows from the invariance of the distribution
of the Brownian map under uniform re-rooting, which will
be discussed below.

Let us first recall some basic definitions concerning pointed metric spaces
and the Gromov-Hausdorff convergence. Let $k\geq 0$ be an integer.
A $k$-pointed compact metric space is a triplet $(E,d,(a_1,a_2,\ldots,a_k))$, where $(E,d)$ is
a compact metric space and $(a_1,a_2,\ldots,a_k)$ is a 
$k$-tuple of (not necessarily distinct) distinguished points of $E$. 
A $0$-pointed compact metric space is just a compact metric space, and
when $k=1$ we say pointed rather than $1$-pointed. We denote  
 the space of all isometry classes of $k$-pointed compact metric spaces by $\K^{(k)}$
(and in agreement with the introduction above, we set $\K=\K^{(0}$).
The Gromov-Hausdorff distance on $\K^{(k)}$ is then defined by
$$d_{GH}((E,d,(a_1,\ldots,a_k)),(E',d',(a'_1,\ldots,a'_k)))\!=\!\inf_{(F,\delta),\phi,\phi'} \!\Big\{\delta_H(\phi(E),\phi'(E'))\vee \max_{1\leq i\leq k}\delta(\phi(a_i),\phi'(a'_i))\!\Big\}$$
where the infimum is over all choices of the metric space $(F,\delta)$ and the isometric 
embeddings $\phi:E\la F$ and $\phi':E'\la F$, and $\delta_H$ stands for the usual
Hausdorff distance between compact subsets of $F$. Then $d_{GH}$
is a distance on $\K^{(k)}$, and the space $(\K^{(k)},d_{GH})$ is Polish.
We equip $\K^{(k)}$ with the associated Borel $\sigma$-field.

 An equivalent definition of the
Gromov-Hausdorff distance can be given in terms of correspondences. 
Let $(E,d,(a_1,\ldots,a_k))$ and 
$(E',d',(a'_1,\ldots,a'_k))$
be two $k$-pointed compact metric spaces.
A correspondence ${\cal C}$ between $(E,d)$ and $(E',d')$ is a subset of the Cartesian
product $E\times E'$ such that for every $x\in E$ there exists $y\in E'$ such that
$(x,y)\in {\cal C}$, and conversely. The distortion of $\cal C$ is defined by
$${\rm dist}({\cal C})=\sup\{|d(x,y)-d'(x',y')|:(x,x')\in{\cal C}\hbox{ and }(y,y')\in{\cal C}\}.$$
Then $d_{GH}((E,d,(a_1,\ldots,a_k)),(E',d',(a'_1,\ldots,a'_k)))$ equals half the infimum of the quantities 
${\rm dist}({\cal C})$ when ${\cal C}$ varies over correspondences between
$E$ and $E'$ that contain all the pairs $(a_i,a'_i)$ for $1\leq i\leq k$. See Theorem 7.3.25 in \cite{BBI}
for a proof in the case $k=0$, which is easily extended.

As was already stated in Theorem \ref{mainIM}, we have the almost sure convergence in $(\K^{(1)},d_{GH})$
\begin{equation}
\label{GHconv}
(\mm_n,\kappa_p n^{-1/4}\,d_n,\partial_n)\build{\la}_{n\to\infty}^{} (\mm_\infty,D,\rho).
\end{equation}

The next theorem gives a precise form of the invariance of the distribution
of the Brownian map under uniform re-rooting

\begin{theorem}
\label{rerootinv}
Let $F$ be a nonnegative measurable function on $\K^{(1)}$. Then
$$E\Big[\int \lambda(dx)\,F(\mm_\infty,D,x)\Big]=E[F(\mm_\infty,D,\rho)].$$
More generally, for every integer $k\geq 2$ and every 
nonnegative measurable function $F$ on $\K^{(k)}$,
\begin{eqnarray*}
&&E\Big[\int \lambda(dx_1)\ldots \lambda(dx_k)\,
F(\mm_\infty,D,(x_1,\ldots,x_k))\Big]\\
&&\qquad=
E\Big[\int \lambda(dx_1)\ldots \lambda(dx_{k-1})\,F(\mm_\infty,D,(\rho,x_1,\ldots,x_{k-1}))
\Big].
\end{eqnarray*}

\end{theorem}

\proof We deal only with the first statement of the theorem. 
The second one can be proved in a similar manner.
Let $U$ be a random variable that is uniformly distributed over $[0,1]$ and independent
of the sequence of random maps $(M_n)_{n\geq 1}$
(and in particular of $(\mm_\infty,D)$). The first assertion of the theorem will
follow from the fact that the law of $(\mm_\infty, D, \pp(U))$ coincides with
the law of $(\mm_\infty,D,\rho)$.

To establish this fact, we will use a similar invariance property 
for the discrete maps $M_n$ and then pass to the limit $n\to\infty$
using (\ref{GHconv}). We start with some simple combinatorial 
considerations. We let $\m^{n,*}_p$ denote the set of all pairs
$(M,e)$ where $M$ is a rooted $2p$-angulation with $n$ faces ($M\in\m^n_p$)
and $e$ is a distinguished oriented edge of $M$. If $(M,e)\in\m^{n,*}_p$,
we define $\Gamma(M,e)\in\m ^n_p$ by saying
that $\Gamma(M,e)$ is the ``same'' map as $M$ but re-rooted 
at the oriented edge $e$. Then  it is easy to verify that the image under
$\Gamma$ of the uniform distribution over $\m^{n,*}_p$ is the uniform 
distribution over $\m^n_p$.

Now consider the random $2p$-angulation $M_n$. The 
Bouttier-Di Francesco-Guitter bijection of subsection 2.6 allows us to associate with
each $i\in\{0,1,\ldots,pn-1\}$ a unique non-oriented edge $e_i$ of $M_n$, and in this way we get
a bijection between $\{0,1,\ldots,pn-1\}$ and the set of 
all (non-oriented) edges of $M_n$. Set $U_n=\lfloor pnU \rfloor$,
in such a way that $U_n$ is uniformly distributed over $\{0,1,\ldots,pn-1\}$.
Choose one of the two possible orientations of the edge $e_{U_n}$ with equal probabilities,
independently of $U$ and of the sequence of maps $(M_n)_{n\geq 1}$, and
denote the resulting oriented edge by $\wt e_{U_n}$.
Also set $\wt M_n=\Gamma(M_n,\wt e_{U_n})$. Note that
the vertex set and the edges of $\wt M_n$ are the same as those of $M_n$, but the root is different.
By the combinatorial considerations
of the beginning of the proof, $\wt M_n$ has the same distribution as $M_n$.  Writing 
$\wt\partial _n$ for the root vertex of $\wt M_n$, we have thus
\begin{equation}
\label{reroot1}
(\mm_n,\kappa_pn^{-1/4}\,d_n,\wt\partial_n)
\build{=}_{}^{\rm(d)} 
(\mm_n,\kappa_pn^{-1/4}\,d_n,\partial_n).
\end{equation}
Note that both sides of the equality (\ref{reroot1}) are viewed as random variables
with values in $\K^{(1)}$.

In view of (\ref{GHconv}) and (\ref{reroot1}), the proof of the theorem
will be complete if we can verify that
\begin{equation}
\label{reroot2}
(\mm_n,\kappa_p n^{-1/4}\,d_n,\wt\partial_n)\build{\la}_{n\to\infty}^{\rm a.s.} (\mm_\infty,D,\pp(U)).
\end{equation}
To this end, let us first recall the proof of (\ref{GHconv}) in \cite{IM}. We define a
correspondence ${\cal C}_n$ between $(\mm_n,\kappa_p n^{-1/4}\,d_n)$
and $(\mm_\infty,D)$ by declaring that a pair $(a,b)\in\mm_n\times\mm_\infty$
belongs to ${\cal C}_n$ if and only if $a=\partial_n$ and $b=\rho$, or
if there exists $t\in[0,1]$ such that $a=p_n(\lfloor pnt \rfloor)$ and $b=\pp(t)$. Then it easily
follows from the convergence (\ref{basicconv}) that ${\rm dist}({\cal C}_n)$
converges to $0$ a.s. as $n\to\infty$ (see \cite{IM}, but note that the presentation
there is slightly different because \cite{IM} does not deal with pointed spaces -- the
same argument however goes through without change). We may slightly enlarge the 
correspondence ${\cal C}_n$ by including also those pairs $(a,b)\in\mm_n\times\mm_\infty$
for which there exists $t\in[0,1]$ such that $d_n(a,p_n(\lfloor pnt \rfloor))=1$ and $b=\pp(t)$.
If ${\cal C}'_n$ denotes the enlarged correspondence,
the distortion ${\rm dist}({\cal C}'_n)$ still
converges to $0$ a.s. as $n\to\infty$.
However, from our construction of the root edge of $\wt M_n$,
it is clear that we have $d_n(\wt\partial_n,p_n(\lfloor pnU\rfloor))\leq 1$. Therefore
$(\wt\partial _n,\pp(U))\in{\cal C}'_n$, and the fact that 
${\rm dist}({\cal C}'_n)$ converges to $0$ also yields our claim (\ref{reroot2}).
This completes the proof. \cq

\medskip
As a by-product of the preceding proof, we get the following variant of (\ref{GHconv}),
which will be useful in the next section.

\begin{proposition}
\label{conv2points}
We can construct for every $n$ 
a random vertex $V_n$ that is uniformly distributed over $\mm_n$, in
such a way that
\begin{equation}
\label{GHconv2}
(\mm_n,\kappa_p n^{-1/4}\,d_n,(\partial_n,V_n))\build{\la}_{n\to\infty}^{\rm a.s.} (\mm_\infty,D,(\rho,
\pp(U))),
\end{equation}
where the random variable $U$ is uniformly distributed over $[0,1]$ and independent
of the sequence of random maps $(M_n)_{n\geq 1}$, and the convergence holds in the space $(\K^{(2)},d_{GH})$.
\end{proposition}

\proof The argument for (\ref{reroot2}) also gives the convergence
\begin{equation}
\label{reroot3}
(\mm_n,\kappa_p n^{-1/4}\,d_n,(\partial_n,\wt\partial_n))\build{\la}_{n\to\infty}^{\rm a.s.} (\mm_\infty,D,
(\rho,\pp(U))),
\end{equation}
which holds in the space $(\K^{(2)},d_{GH})$. This does not immediately 
give the desired result, because $\wt\partial_n$ is not uniformly distributed over $\mm_n$. 
However, we may set $N_n=\lfloor ((p-1)n+2)U\rfloor$ and let $V_n$ denote the 
$N_n$-th vertex of $\tau^\circ_n$ in lexicographical order, with the convention
that $V_n=\partial_n$ if $N_n=0$. Then, $V_n$ has the property stated in
the proposition, and furthermore
$$n^{-1/4} d_n(V_n,\wt\partial_n)\build{\la}_{n\to\infty}^{\rm a.s.} 0.$$
See the proof of Proposition 3.5 in \cite{IM} for a justification of this convergence.
The result of the proposition now follows from (\ref{reroot3}). \cq

\smallskip
As was mentioned above, Theorem \ref{rerootinv} allows one to get information
about geodesics starting from a point chosen uniformly in $\mm_\infty$
rather than from the root $\rho$. We illustrate this in the following 
corollary. If $x,x'\in \mm_\infty$, we denote by ${\rm Geo}(x\to x')$ the set
of all geodesics from $x$ to $x'$.

\begin{corollary}
\label{newroot}
The following properties hold almost surely.
\begin{description}
\item{\rm (i)} $
\lambda\otimes\lambda\{(x,x'):\# {\rm Geo}(x\to x') >1\}=0.$
\item{\rm (ii)} $\lambda(\{x\in\mm_\infty: \exists x'\in\mm_\infty, \# {\rm Geo}(x\to x')\geq 4\})=0.$
\item{\rm (iii)} $\lambda(\{x\in\mm_\infty: \exists x'\in\mm_\infty, \# {\rm Geo}(x\to x')=3\})=1.$
\end{description}
\end{corollary}

\rem Property (i) of the corollary should be compared to Theorem 3 in \cite{Mi2}.

\smallskip
\proof We start by proving (i). Let $U$ and $U'$ be two independent random variables that are
uniformly distributed over $[0,1]$, and such that the pair $(U,U')$
is independent of the sequence $(M_n)_{n\geq 1}$. The statement in (i)
is equivalent to proving that $\# {\rm Geo}(\pp(U)\to \pp(U'))=1$ a.s.
However, as a simple consequence of Theorem \ref{rerootinv}, we have
$$(\mm_\infty,D,(\rho,\pp(U)))\build{=}_{}^{\rm(d)} 
(\mm_\infty,D,(\pp(U),\pp(U'))),$$
where both sides of the equality are random variables with values in $\K^{(2)}$.
Note that the set of all $2$-pointed spaces $(E,d,(x,x'))$ such that there
exists a unique geodesic from $x$ to $x'$ is a measurable subset of $\K^{(2)}$.
Indeed,
if $\delta>0$ is fixed, the set of all $(E,d,(x,x'))$ such that there exist
two geodesics $\ga$ and $\ga'$ from $x$ to $x'$ with $\max\{d(\ga(t),\ga'(t)\}\geq\delta$
is closed in $\K^{(2)}$. Hence, using the previous identity in distribution, 
we see that property (i) reduces to $\# {\rm Geo}(\rho\to \pp(U))=1$ a.s.
But this follows from Corollary \ref{typicalgeo}.

The argument for (ii) and (iii) is similar. We use the identity in distribution
$$(\mm_\infty,D,\rho)\build{=}_{}^{\rm(d)} 
(\mm_\infty,D,\pp(U)),$$
and the fact that, for every integer $k\geq 2$,
the set of all pointed spaces $(E,d,x)$ such that there exist a
point $x'\in E$ and at least $k$ distinct geodesics from $x$ to $x'$
is a measurable subset of $\K^{(1)}$. Details are left to the reader. \cq

\section{Asymptotics for large planar maps}

In this section we discuss the applications of our results to
large planar maps, and in particular we prove
Propositions \ref{uniquenessdiscrete}, \ref{multiplediscrete} 
and \ref{maximaldiscrete}
that were stated in Section 1. Note that, in contrast with the previous
sections, we do not limit ourselves to values of $n$ belonging
to a suitable sequence converging to $\infty$. 

\medskip
\noindent{\bf Proof of Proposition \ref{uniquenessdiscrete}:}
We prove the stronger form of Proposition \ref{uniquenessdiscrete}
where ${\rm Geo}_n(\partial_n,a)$ is replaced by 
$\ov{\rm Geo}_n(\partial_n,a)$. Let 
us fix $\delta>0$ and, for every $n\geq 1$, let $V_n$ be a vertex that
is uniformly distributed over $\mm_n$. The desired result is equivalent to
saying that
$$\lim_{n\to\infty}
P\Big(\exists \ga,\ga'\in \ov{\rm Geo}_n(\partial_n,V_n): d(\ga,\ga')\geq \delta n^{1/4}\Big)
=0.$$
We argue by contradiction and suppose that there exists $\eta>0$ and
a subsequence $\ov n_k\uparrow\infty$, such that, for every $k$,
$$P\Big(\exists \ga,\ga'\in \ov{\rm Geo}_{\ov n_k}(\partial_{\ov n_k},V_{\ov n_k}): 
d(\ga,\ga')\geq \delta {\ov n_k}^{1/4}\Big)\geq \eta.$$
From this sequence $(\ov n_k)$ we can now extract another subsequence 
$(n_k)$ such that the properties stated in Theorem \ref{mainIM} hold. 
By Proposition \ref{conv2points}, we may assume that the random vertices $V_{n_k}$
are constructed in such a way that
\begin{equation}
\label{GHconv2bis}
(\mm_{n_k},\kappa_p n_k^{-1/4}\,d_{n_k},(\partial_{n_k},V_{n_k}))
\build{\la}_{{k}\to\infty}^{\rm a.s.} (\mm_\infty,D,(\rho,
\pp(U))),
\end{equation}
where $U$ is uniformly distributed over $[0,1]$ and independent of 
the sequence $(M_n)_{n\geq 1}$.
Consider then the event
$$F=\limsup_{k\to\infty} 
\{\exists \ga,\ga'\in \ov{\rm Geo}_{n_k}(\partial_{n_k},V_{n_k}): d(\ga,\ga')\geq \delta n_k^{1/4}\}.$$
By our assumptions, $P(F)\geq \eta>0$. 

From now on, we argue on the event $F$. We can then find a (random)
subsequence $(m_k)$ of the sequence $(n_k)$, such that, for every $k$, there exist 
$\ga_k,\ga'_k\in \ov{\rm Geo}_{m_k}(\partial_{m_k},V_{m_k})$ such that
$d(\ga_{k},\ga'_{k})\geq \delta m_k^{1/4}$. On the other hand,
by (\ref{GHconv2bis}), we can find for every $k$
a correspondence ${\cal C}_k$ between the $2$-pointed spaces
$(\mm_{m_k},\kappa_p m_k^{-1/4}\,d_{m_k},(\partial_{m_k},V_{m_k}))$
and $(\mm_\infty,D,(\rho,
\pp(U)))$, in such a way that ${\rm dist}({\cal C}_k)$ tends to $0$
as $k\to\infty$. This implies in particular that $\kappa_p m_k^{-1/4} d_{m_k}(\partial_{m_k},V_{m_k})$
converges to $D(\rho,\pp(U))=\ov Z_U$. For every rational
$r\in[0,\ov Z_U[$, and for $k$ sufficiently large, we can find 
$x_k(r),x'_k(r)\in \mm_\infty$ such that $(\ga_k(\lfloor \kappa_p^{-1}m_k^{1/4}r\rfloor),x_k(r))\in{\cal C}_k$
and $(\ga'_k(\lfloor \kappa_p^{-1}m_k^{1/4}r\rfloor),x'_k(r))\in{\cal C}_k$.
Via a compactness argument and extracting a diagonal subsequence, we can assume that
$x_k(r)$, respectively $x'_k(r)$, converges as $k\to\infty$, for every rational $r\in[0,\ov Z_U[$, and we denote the
limit by $\ga_\infty(r)$, resp. $\ga'_\infty(r)$. The fact that each $\ga_k$ is an approximate geodesic, together with the
property ${\rm dist}({\cal C}_k)\to 0$, implies that $D(\ga_\infty(r),\ga_\infty(r'))=|r-r'|$
for every rationals $r,r'\in[0,\ov Z_U[$. It easily follows that $\ga_\infty$ can be extended to
a geodesic from $\rho$ to $\pp(U)$ in $\mm_\infty$. Similarly, we can extend $\ga'_\infty$
to a geodesic from $\rho$ to $\pp(U)$ in $\mm_\infty$. However,
$$\max_{0\leq r\leq \ov Z_U} D(\ga_\infty(r),\ga'_\infty(r))
\geq \liminf_{k\to\infty} \kappa_p m_k^{-1/4}\,d(\ga_k,\ga'_k)\geq \kappa_p\delta >0$$
and so we find that there are two distinct geodesics from $\rho$ to $\pp(U)$
in $\mm_\infty$. This holds on the event $F$, which has positive probability, thus contradicting
Corollary \ref{typicalgeo}. \cq

\medskip
\noi{\bf Proof of Proposition \ref{multiplediscrete}:} This is very
similar to the proof of Proposition \ref{uniquenessdiscrete}. We argue by contradiction
and assume that there is an increasing  sequence $(\ov n_k)_{k\geq 1}$ such that,
for every $k\geq 1$,
$$
P(\exists a\in \mm_{\ov n_k}: \ov{\rm Mult}_{\ov n_k,\delta}(\partial_{\ov n_k},a)\geq 4) \geq \eta$$
for some $\eta >0$. We can then extract a subsequence $(n_k)_{k\geq 1}$ and construct the
random maps $M_{n_k}$ in such a way that we have the almost sure convergence
\begin{equation}
\label{convGH1}
(\mm_{n_k},\kappa_p n_k^{-1/4}\,d_{n_k},\partial_{n_k})
\build{\la}_{{k}\to\infty}^{\rm a.s.} (\mm_\infty,D,\rho)
\end{equation}
in $\K^{(1)}$. The event
$$F=\limsup_{k\to\infty} \{\exists a\in \mm_{n_k}: \ov{\rm Mult}_{n_k,\delta}(\partial_{n_k},a)\geq 4\}$$
has probability bounded below by $\eta$. Furthermore, on the event $F$, we can use the
convergence (\ref{convGH1}) and compactness arguments to construct a
point $x$ of $\mm_\infty$ which is connected to $\rho$ by at least $4$ distinct 
geodesics, thus contradicting Theorem \ref{numbergeodesics}. We leave details to
the reader. \cq

\rem As we already mentioned in the introduction, the statements of Propositions \ref{uniquenessdiscrete} and \ref{multiplediscrete}
remain valid if $\partial_n$
is replaced by a vertex chosen uniformly at random in $\mm_n$. Let us outline
the argument in the case of Proposition \ref{uniquenessdiscrete}. We consider two random vertices
$V_n$ and $V'_n$ that are distributed independently uniformly over the vertex set of $\mm_n$.
From the argument in the proof of Proposition \ref{conv2points}, we can construct $V_n$ and $V'_n$,
and the random maps $M_n$, in such a way that along a suitable sequence
$(n_k)_{k\geq 1}$ we have the almost sure convergence
\begin{equation}
\label{GHconv2aux}
(\mm_n,\kappa_p n^{-1/4}\,d_n,(V_n,V'_n))\build{\la}_{n\to\infty}^{\rm a.s.} (\mm_\infty,D,(\pp(U),
\pp(U'))),
\end{equation}
where $U$ and $U'$ are independent and uniformly distributed over $[0,1]$, and 
$(U,U')$ is independent of the sequence $(M_n)_{n\geq 1}$. Then a simple adaptation of the
proof of Proposition \ref{uniquenessdiscrete}, using Corollary \ref{newroot} (i), shows that
$P(\ov{\rm Mult}_{n,\delta}(V_n,V'_n)\geq 2)$ tends to $0$ as $n\to\infty$
along the sequence $(n_k)_{k\geq 1}$. Since from any sequence $(n'_k)$ we can 
extract a subsequence $(n_k)$ such that the previous conclusion holds, the
desired result follows. 

\medskip
\noindent{\bf Proof of Proposition \ref{maximaldiscrete}:} Recall our notation
$v^n_0,v^n_1,\ldots,v^n_{pn}$ for the contour sequence of $\tau^\circ_n$. Let
$\al>0$. We denote by ${\cal H}_n(\al)$ the set of all triplets $(j_1,j_2,j_3)$
of integers such that $0\leq j_1\leq j_2\leq j_3\leq pn$ and the following properties hold:
\begin{description}
\item{$\bullet$} $v^n_{j_1}=v^n_{j_2}=v^n_{j_3}$\,;
\item{$\bullet$} ${\displaystyle
\max\Big(\min_{j_1 \leq j\leq j_2} \Lambda^n_j,\min_{j_2 \leq j\leq j_3} \Lambda^n_j\Big)
\leq \Lambda^n_{j_1} - \al n^{1/4}.}$
\end{description}
Suppose that  $(j_1,j_2,j_3)\in{\cal H}_n(\al)$. Then we have $C^n_{j_1}=C^n_{j_2}=C^n_{j_3}$,
$\Lambda^n_{j_1}=\Lambda^n_{j_2}=\Lambda^n_{j_3}$
and $C^n_j\geq C^n_{j_1}$ if $j_1\leq j\leq j_3$. From the properties of the spatial
contour function, it is also clear that $\Lambda^n$ takes all integer values between
$\Lambda^n_{j_1}$ and $\Lambda^n_{j_1} - \al n^{1/4}$ on each of the
intervals $[j_1,j_2]$, $[j_2,j_3]$ and $[j_3,pn]$.

We denote the event $\{{\cal H}_n(\al)\not =\emptyset\}$ by $\Delta_n(\al)$.

\begin{lemma}
\label{maxitech}
We have
$$\lim_{\al\to 0} \Big(\liminf_{n\to\infty} P(\Delta_n(\al))\Big)=1.$$
\end{lemma}

We postpone the proof of Lemma \ref{maxitech}
 to the Appendix.

Let $\delta>0$. For $i=1$ or $2$, we say that the event $\Delta^{(i)}_n(\al,\delta)$
holds if, for every triplet $(j_1,j_2,j_3)\in{\cal H}_n(\al)$, for every
$j\in[j_i,j_{i+1}]$ such that
$$\Lambda^n_j=\min_{j_i\leq \ell\leq j} \Lambda^n_\ell
\quad\hbox{and}\quad \Lambda^n_{j_1}-\frac{2\al}{3}n^{1/4}\leq \Lambda^n_j
\leq \Lambda^n_{j_1}-\frac{\al}{3}n^{1/4}$$
we have 
\begin{equation}
\label{maxitech1}
d_n(j,j')\geq \delta\,n^{1/4}\;,\quad\hbox{for every }j'\in\{ j_{i+1},\ldots,pn\}.
\end{equation}
We claim that, for every $\al>0$,
\begin{equation}
\label{maxitech2}
\sup_{n\geq 1} P\Big(\Delta_n(\al)\backslash(\Delta^{(1)}_n(\al,\delta)
\cap \Delta_n^{(2)}(\al,\delta))\Big) \build{\la}_{\delta\to 0}^{} 0.
\end{equation}

Assuming the claim, it is easy to complete the proof of Proposition \ref{maximaldiscrete}.
Let $\eta>0$. Using Lemma \ref{maxitech} and (\ref{maxitech2}),
we can choose first $\alpha>0$, and then $\delta>0$ in such a way that
$$ \liminf_{n\to\infty} P\Big(\Delta_n(\al)
\cap (\Delta^{(1)}_n(\al,\delta)
\cap \Delta_n^{(2)}(\al,\delta))\Big)>1-\eta.$$
However, on the event $\Delta_n(\al))
\cap (\Delta^{(1)}_n(\al,\delta)
\cap \Delta_n^{(2)}(\al,\delta))$ there is at least one
point $a\in\mm_n$ such that ${\rm Mult}_{n,\delta}(\partial_n,a)\geq 3$.
Indeed we pick a triplet $(j_1,j_2,j_3)\in{\cal H}_n(\al)$ and set $a=v^n_{j_1}=v^n_{j_2}=v^n_{j_3}$.
We consider the ``simple geodesics'' $\ga_1,\ga_2,\ga_3$ defined by
\begin{eqnarray*}
&&\ga_i(\ell)=v^n_{m_i(\ell)}\quad\hbox{where }m_i(\ell)
=\inf\{j\geq j_i:\Lambda^n_j=\Lambda_{j_1}-\ell\}\hbox{ for } 0\leq \ell\leq \Lambda^n_{j_1}-1,\\
&&\ga_i(\Lambda^n_{j_1})=\partial_n
\end{eqnarray*}
for $i=1,2,3$.
It is then easy to verify that $\ga_1,\ga_2,\ga_3$ are three geodesics from 
$a$ to $\partial_n$ in the map $M_n$, and the property (\ref{maxitech1})
guarantees that $d(\ga_i,\ga_{i'})\geq \delta n^{1/4}$
for every $i,i'\in\{1,2,3\}$ with $i\not =i'$. 

We still have to prove our claim (\ref{maxitech2}). We verify that
\begin{equation}
\label{maxitech3}
\sup_{n\geq 1} P\Big(\Delta_n(\al)\backslash\Delta^{(1)}_n(\al,\delta)
\Big) \build{\la}_{\delta\to 0}^{} 0.
\end{equation}
and the same argument applies if $\Delta^{(1)}_n(\al,\delta)$ is replaced
by $\Delta^{(2)}_n(\al,\delta)$. We argue by contradiction and assume that
(\ref{maxitech3}) does not hold. Then we can find $\ve>0$ such that, for every
$k\geq 1$,
$$\sup_{n\geq 1} P\Big(\Delta_n(\al)\backslash\Delta^{(1)}_n(\al,\frac{1}{k})\Big)>\ve.$$
Hence, for every $k\geq 1$, we can find a positive integer $\ov n_k$
such that
$$P\Big(\Delta_{\ov n_k}(\al)\backslash\Delta^{(1)}_{\ov n_k}(\al,\frac{1}{k})\Big)>\ve.$$
The sequence $(\ov n_k)_{k\geq1}$ must converge to $+\infty$ (for any fixed $n$,
the condition (\ref{maxitech1}) holds automatically if $\delta$ is small enough).
So from this sequence, we can extract a monotone increasing subsequence $(n_k)_{k\geq 1}$
along which the almost sure convergence (\ref{basicconv}) holds for 
a suitable choice of the random maps $M_{n_k}$. We have then
$$P\Big(\limsup_{k\to\infty}\, (\Delta_{n_k}(\al)\backslash\Delta^{(1)}_{n_k}(\al,\frac{1}{k}))
\Big)\geq \ve.$$
Furthermore, on the event 
$$\limsup_{k\to\infty}\, (\Delta_{n_k}(\al)\backslash\Delta^{(1)}_{ n_k}(\al,\frac{1}{k}))$$
we can find a (random) sequence of values of $k$ along which the following holds:
There exists a triplet $(j^k_1,j^k_2,j^k_3)\in {\cal H}_{n_k}(\al)$ and two integers
$m_k\in[j^k_1,j^k_2]$ and $m'_k\in[j^k_2,pn_k]$ such that
$$\Lambda^{n_k}_{m_k}=\min_{j_1^k\leq \ell\leq m_k} \Lambda^{n_k}_\ell
\quad ,\quad \Lambda^{n_k}_{j_1^k}-\frac{2\al}{3}n_k^{1/4}\leq \Lambda^{n_k}_{m_k}
\leq \Lambda^{n_k}_{j_1^k}-\frac{\al}{3}n_k^{1/4}$$
and 
\begin{equation}
\label{maxitech4}
d_{n_k}(m_k,m'_k)\leq \frac{1}{k}\,n_k^{1/4}.
\end{equation}
By a compactness argument, we may assume that
$$\frac{j^k_1}{pn_k}\la s\ ,\quad \frac{j^k_2}{pn_k}\la t\ ,\quad \frac{m_k}{pn_k}\la r\,
\quad \frac{m'_k}{pn_k}\la r'$$
as $k\to\infty$. We have then $s\leq r\leq t\leq r'$, $s\sim t$ (and so in particular $\ov Z_s=\ov Z_t$)
and $D(r,r')=0$ from (\ref{basicconv}) and (\ref{maxitech4}). Furthermore, the
convergence (\ref{basicconv}) also implies that
$$\inf_{s\leq u\leq t} \ov Z_u\leq \ov Z_s -\al \kappa_p$$
and
$$\ov Z_r=\inf_{s\leq u\leq r} \ov Z_u\ ,\quad \ov Z_s-\frac{2\al}{3}\kappa_p
\leq \ov Z_r\leq  \ov Z_s-\frac{\al}{3}\kappa_p.$$
From these properties, it follows that neither of the conditions $r\sim r'$
or $r\approx r'$ can hold. This contradicts the equality $D(r,r')=0$, and this
contradiction completes the proof. \cq

\smallskip
\rem In the same way as Propositions \ref{uniquenessdiscrete} and \ref{multiplediscrete},
Proposition \ref{maximaldiscrete} still holds if the root vertex $\partial_n$
is replaced by a vertex chosen uniformly at random in $\mm_n$. A way to obtain this
result is to use
the variant of the Bouttier-Di Francesco-Guitter bijection that is presented
in Section 2.3 of \cite{MaMi}. We omit the details of the argument.

\smallskip
We conclude this section with a discrete version of Corollary \ref{coincid}.

\begin{proposition}
\label{coinciddiscret}
 Let $\chi >0$ and $\ve >0$, and let ${\rm Geo}_n^\chi$
 denote the set of all discrete geodesics starting from the root
 $\partial_n$ and with length greater than $\chi n^{1/4}$
 in the map $M_n$. Then we can choose a constant $\beta \in]0,\chi[$ 
so that, for every $\delta >0$,
 $$\liminf_{n\to\infty} 
 P\Big(\sup_{\ga,\ga'\in{\rm Geo}_n^\chi}\Big\{\sup_{0\leq i\leq 
 \lfloor\beta n^{1/4}\rfloor} d_n(\ga(i),\ga'(i))\Big\}
 \leq \delta n^{1/4}\Big) \geq 1-\ve.$$
 \end{proposition}
 
 Proposition \ref{coinciddiscret} is derived from Corollary \ref{coincid}
 and the Gromov-Hausdorff convergence in Theorem \ref{mainIM}
 in a way very similar to the above proof of Proposition
 \ref{uniquenessdiscrete}. To be specific, if $\alpha>0$ is the random
 variable appearing in the statement of Corollary \ref{coincid}
 when $\eta=\kappa_p\chi$, we choose the constant $\beta>0$
 so that $P(\kappa_p^{-1}\alpha >\beta)>1-\ve$. From the proof
 of Corollary \ref{coincid}, we see that the choice 
 of $\beta$ is determined from the value of $\chi$ and $\ve$ and
 the distribution of the process $\ov Z$, so that it does not depend on
 the particular subsequence used in Theorem \ref{mainIM}. The remaining part of the argument
 is straightforward.

\section*{Appendix}

In this Appendix, we prove Lemma \ref{excursion} and Lemma \ref{maxitech}.
As previously the integer $p\geq 2$ is fixed, and all constants that will
appear depend on $p$, even though this will not be mentioned explicitly.

\medskip

\noi{\bf Proof of Lemma \ref{excursion}:} The proof depends on an intermediate
estimate on the distribution of labels in a labelled $p$-tree. We say that a 
pair $(\tau,(\ell_v)_{v\in\tau^\circ})$ consisting of a $p$-tree $\tau$
and a collection of integer labels assigned to the white vertices of $\tau$
is a labelled $p$-tree if it satisfies property (b) in the definition of a $p$-mobile
(cf subsection 2.6). We say that it is a standard labelled $p$-tree if in
addition $\ell_\varnothing=0$. In a labelled $p$-tree, labels may take
negative values. 

Let $(\tau,(\ell_v)_{v\in\tau^\circ})$ be a standard labelled $p$-tree, and let
$v\in\tau^\bullet$. Let $v_{(0)},v_{(1)},\ldots,v_{(p-1)}$ be the neighboring vertices
of $v$ listed as in property (b) of the definition of a $p$-mobile. Simple combinatorics show that
given $\ell_{v_{(0)}}$ there are $2p-1\choose p-1$ possible choices for
$\ell_{v_{(1)}},\ldots,\ell_{v_{(p-1)}}$. As an immediate consequence,
given a $p$-tree $\tau$ with $n$ black vertices, there are ${2p-1\choose p-1}^n$
choices for the labels $\ell_v,v\in\tau^\circ$ that produce a standard labelled $p$-tree.

We denote by $\Pi^p_n$ the uniform distribution over the set
of all standard labelled $p$-trees with $n$ black vertices.

\medskip
\noi{\bf Lemma}
{\it For every integer $k\geq 1$, there exists a constant $C_k$ such that, for every $n\geq 1$
and every $x>0$,
\beq
\label{ap1}
\Pi^p_n\Big(\max_{v\in\tau^\circ} |\ell_v|\geq x\,n^{1/4}\Big)
\leq {C_k\over x^k}.
\eeq
Consequently, for every $\ve >0$, there exists a constant
$C_{(\ve)}$ such that, for every $n\geq 1$
and every $x\geq \ve$,
$$\Pi^p_n\Big(\min_{v\in\tau^\circ} \ell_v\leq -x\,n^{1/4}\Big)
\leq 1-\exp(-C_{(\ve)}\,x^{-4}).$$
}

\medskip
\noi{\bf Remark.} When $p=2$, the bound (\ref{ap1}) follows from the much stronger exponential
estimate in Proposition 4 of \cite{CS}. It is very plausible that such an exponential
estimate holds for general $p$, but we content ourselves with the weaker bound (\ref{ap1}).
This bound is sufficient to get the second assertion of the lemma, which is what we really
need.

\medskip
\noi{\bf Proof of the lemma:} We write 
$(\wt\tau_n,(\wt\ell^n_v)_{v\in\wt\tau_n^\circ})$
for a random labelled $p$-tree which is uniformly distributed
over the set
of all standard labelled $p$-trees with $n$ black vertices.
From the remarks preceding the lemma, it is clear that $\wt\tau_n$
is uniformly distributed over the set of all $p$-trees with $n$
black vertices. For our purposes, it will be convenient to view $\wt\tau_n$
as a conditioned multitype Galton-Watson tree. To this end, consider the random
$p$-tree $\tau_*$ whose distribution is informally specified as
follows:
\begin{description}
\item{(i)} for every $v\in\tau_*^\bullet$, $k_v(\tau_*)=p-1$;
\item{(ii)} for every $v\in\tau_*^\circ$, 
$k_v(\tau_*)$ is geometrically distributed with parameter $p^{-1}$, that is
$$P(k_v(\tau_*)=k\mid v\in\tau_*)=(1-p^{-1})p^{-k},\qquad k\geq 0.$$
\end{description}
Then $\wt\tau_n$ has the same distribution as $\tau_*$
conditioned on the event $\#\tau^\bullet_*=n$, or equivalently
on the event $\#\tau^\circ_*=(p-1)n+1$. 

Now note that $\tau_*^\circ$ can be viewed as a random
plane tree, by declaring that two vertices $v$
and $v'$ of $\tau_*^\circ$ are linked by an edge if
$v$ is the grandparent of $v'$ or conversely (and of course
keeping the root and the order structure). With this 
interpretation, $\tau^\circ_*$ is a Galton-Watson tree,
whose offspring distribution is the law of 
$p-1$ times a geometric random variable with parameter $p^{-1}$.

Recall that $\wt\tau^\circ_n$ has the same distribution
as $\tau^\circ_*$ conditioned on the event $\#\tau^\circ_*=(p-1)n+1$.
We denote by $\wt v^n_0,\ldots,\wt v^n_{2(p-1)n}$
the depth-first search sequence of $\wt\tau^\circ_n$
viewed as a plane tree, and we set
$${\cal C}^n_i=\frac{1}{2} |\wt v^n_i|$$
for every $i\in\{0,1,\ldots,2(p-1)n\}$. An estimate of Gittenberger \cite{Gi}
gives two positive constants $C_0$ and $C_1$, which do not depend
on $n$, such that, for every $i,j\in\{0,1,\ldots,2(p-1)n\}$
and every $x>0$,
\beq
\label{ap2}
P\Big({\cal C}^n_i+{\cal C}^n_j - 2\min_{i\wedge j\leq q \leq i\vee j}
{\cal C}^n_q\geq x\Big)
\leq \frac {C_0n}{|i-j|}\,\exp\Big(-\frac{C_1 x}{\sqrt{|i-j|}}\Big).
\eeq

Let $N>4$ be an integer. As an immediate consequence of the bound (\ref{ap2}), 
we get
\beq
\label{ap3}
E\Big[\Big( n^{-1/2}\Big({\cal C}^n_i+{\cal C}^n_j - 2\min_{i\wedge j\leq q \leq i\vee j}
{\cal C}^n_q\Big)\Big)^N\Big]
\leq C'_N \Big(\frac{|i-j|}{n}\Big)^{\frac{N}{2}-1},
\eeq
with a constant $C'_N$ that does not depend on $n$. 

Fix $i,j\in\{0,1,\ldots,2(p-1)n\}$ such that $i<j$, and let 
$m\in\{i,i+1,\ldots,j\}$ be such that $\wt v^n_m$ is the most recent
common ancestor to $\wt v^n_i$ and $\wt v^n_j$
in the tree $\wt\tau^\circ_n$. Then, conditionally 
on $\wt\tau_n$, $\wt\ell^n_{\wt v^n_i}-\wt\ell^n_{\wt v^n_m}$ is distributed as
the sum of $\frac{1}{2}(|\wt v^n_i|-|\wt v^n_m|)$ integer-valued
independent random variables,
which are centered and bounded by $p$ in absolute value. A similar property
holds for $\wt\ell^n_{\wt v^n_j}-\wt\ell^n_{\wt v^n_m}$. Noting that
$$\frac{1}{2}(|\wt v^n_i|+|\wt v^n_j|-2|\wt v^n_m|)
={\cal C}^n_i+{\cal C}^n_j - 2\min_{i\wedge j\leq q \leq i\vee j}
{\cal C}^n_q,$$
we can use standard arguments to derive the bound
\beq
\label{ap4}
E[|\wt\ell^n_{\wt v^n_i}-\wt\ell^n_{\wt v^n_j}|^{2N}\mid \wt\tau_n]
\leq C''_N \Big({\cal C}^n_i+{\cal C}^n_j - 2\min_{i\wedge j\leq q \leq i\vee j}
{\cal C}^n_q\Big)^N 
\eeq
where the constant $C''_N$ only depends on $N$.

Set ${\cal V}^n_i=\ell^n_{\wt v^n_i}$ for $i\in\{0,1,\ldots,2(p-1)n\}$ and extend
${\cal V}^n$ to $[0,2(p-1)n]$ by interpolating linearly between integers.
From (\ref{ap3}) and (\ref{ap4}), we have, for $s,t\in \{0,1,\ldots,2(p-1)n\}$,
\beq
\label{ap5}
E\Big[\Big|\frac{{\cal V}^n_s-{\cal V}^n_t}{n^{1/4}}\Big|^{2N}\Big]
\leq \ov C_N\,\Big(\frac{|s-t|}{n}\Big)^{\frac{N}{2}-1}
\eeq
where $\ov C_N=C'_NC''_N$. This bound remains valid for every reals $s,t\in[0,2(p-1)n]$,
with a possibly different constant $\ov C_N$ that still depends only on $N$. 

Set $\wh{\cal V}^n_s=n^{-1/4}{\cal V}^n_{ns}$ for every $s\in[0,2(p-1)]$. From
(\ref{ap5}), we have for every $s,t\in[0,2(p-1)]$,
$$E[|\wh{\cal V}^n_s-\wh{\cal V}^n_t|^{2N}]\leq \ov C_N\,|s-t|^{\frac{N}{2}-1}.$$
Since $\frac{N}{2}-1>1$, an application of the Kolmogorov criterion
(see e.g. Theorem I.2.1 in \cite{RY}) gives a
constant $\wt C_N$, which does not depend on $n$, such that
$$E\Big[\sup_{0\leq s\leq 2(p-1)} |\wh {\cal V}^n_s|^{2N}\Big] \leq \wt C_N.$$
The bound (\ref{ap1}) follows.  

To get the second assertion, first note that 
for every $\ve >0$ there exists a constant $C'_{(\ve)}>0$
such that, for every integer $n\geq 1$,
\beq
\label{bound5bis}
1-\Pi^p_n\Big(\min_{v\in\tau^\circ}\ell_v\leq -\ve n^{1/4}\Big)\geq C'_{(\ve)}.
\eeq
Indeed, it follows from Theorem 11 in \cite{MaMi} that
$$\limsup_{n\to\infty} \Pi^p_n\Big(\min_{v\in\tau^\circ}\ell_v\leq -\ve n^{1/4}\Big)
\leq P\Big(\min_{0\leq s\leq 1} Z^\eg_s\leq -\kappa_p\ve\Big),$$
and the limit is strictly less than $1$, as can be seen for instance
from the invariance of the CRT under uniform re-rooting.
Using (\ref{bound5bis}) and the first assertion of the lemma with $k=4$,
we get for every $x\geq \ve$,
$$\Pi^p_n\Big(\min_{v\in\tau^\circ}\ell_v\leq -x n^{1/4}\Big)\leq (1-C'_{(\ve)})
\wedge C_4x^{-4} \leq 1- \exp(-C_{(\ve)}x^{-4})$$
for some other constant $C_{(\ve)}$.
This completes the proof of the lemma.
\cq

\smallskip
We now turn to the proof of Lemma \ref{excursion}. We will in fact prove a dual version
of this lemma. Precisely, we will show that, for every choice of $U\in\,]0,1[$, $u_0\in\,]U,1[$,
$\delta\in\,]0,\frac{1-u_0}{2}[$ and $\eta,\eta'>0$, there exists a constant
$K=K(U,u_0,\eta,\eta',\delta)$ such that the following holds. Whenever
$u\in\,]U,u_0[$ and $\alpha>0 $, we have for all $n$ sufficiently large
\beq
\label{ap6}
P\Big(\inf_{\lfloor pnu\rfloor\leq i\leq \lfloor pn(u+\delta)\rfloor} C^n_i
\geq C^n_{\lfloor pnu \rfloor} -\al n^{1/2} \,\Big|\,\g^n(0,\lfloor pnu \rfloor)\Big)
\leq K\,\al
\eeq
on the event 
\beq
\label{ap7}
\Big\{C^n_{\lfloor pnu \rfloor}\geq \eta n^{1/2},
\inf_{\eta n^{1/2}/2\leq j \leq C^n_{\lfloor pnu \rfloor}} W^n_{\lfloor pnu \rfloor}(j)
\geq \eta' n^{1/4}\Big\}.
\eeq
Strictly speaking, this dual version is not equivalent to Lemma \ref{excursion},
because of the lack of symmetry in property (b) of the definition of a $p$-mobile.
However, the reader will easily check that the same arguments we will use
to prove the dual version also apply to the statement of Lemma \ref{excursion},
at the cost of a somewhat heavier notation.

In proving the bound (\ref{ap6}), we may and will assume that $\al\leq \eta/(4p)$.
Let us fix $u\in\,]U,u_0[$, and to simplify notation set $N=C^n_{\lfloor pnu\rfloor}$.
Recall our notation $v^n_0,v^n_1,\ldots,v^n_{pn}$ for the contour sequence of $\tau_n^\circ$, and
let
$$w^n_0=\varnothing,w^n_1,w^n_2,\ldots,w^n_N=v^n_{\lfloor pnu \rfloor}$$
be the white vertices of $\tau_n$ that are ancestors of $v^n_{\lfloor pnu \rfloor}$,
listed in such a way that $|w^n_i|=2i$ for $0\leq i\leq N$. Also let
$\wt w^n_1,\wt w^n_2,\ldots,\wt w^n_N$ be the black vertices that 
are ancestors of $v^n_{\lfloor pnu \rfloor}$, 
listed in such a way that $|\wt w^n_i|=2i-1$ for $1\leq i\leq N$.

Consider now all vertices $x\in\tau^\circ_n$
that are of one of the following two types:
\begin{description}
\item{$\bullet$} either $x\in\{w^n_0,w^n_1,\ldots,w^n_N\}$;
\item{$\bullet$} or there exists $i\in\{1,\ldots,N\}$ such that $x$ is a child of $\wt w^n_i$ that has not yet appeared in the 
sequence $v^n_0,v^n_1,\ldots,v^n_{\lfloor pnu \rfloor}$.
\end{description}
Let $x^n_0,x^n_1,\ldots,x^n_M$ be the sequence consisting of all such vertices $x$
listed in such a way that $|x^n_i|\leq |x^n_{i'}|$ if $i\leq i'$, and in lexicographical order
for vertices of the same generation. Notice that $N\leq M\leq (p-1)N$.

It is easily verified that the quantities $N,M$ and the random vertices
$x^n_0,x^n_1,\ldots,x^n_M$ are measurable with respect to the $\sigma$-field
$\g^n(0,\lfloor pnu \rfloor)$. This measurability property does not
hold for the labels $\ell^n_{x^n_i}$, $0\leq i\leq M$
(unless $p=2$). So we let $\wh\g^n(0,\lfloor pnu \rfloor)$ be the $\sigma$-field
generated by $\g^n(0,\lfloor pnu \rfloor)$ and the labels $\ell^n_{x^n_i}$, $0\leq i\leq M$,
and we will prove the stronger form of (\ref{ap6}) where the $\sigma$-field $\g^n(0,\lfloor pnu \rfloor)$
is replaced by $\wh\g^n(0,\lfloor pnu \rfloor)$. 

Let $k$ and $m$ be two positive integers, with $k\leq m\leq (p-1)k$, and let
$\beta_0,\beta_1,\ldots,\beta_{\lfloor pnu \rfloor}$, 
$\nu_0,\nu_1,\ldots,\nu_{\lfloor pnu \rfloor}$, $l_0,l_1,\ldots,l_m$
be nonnegative integers such that
$$P\Big(N=k,M=m;C^n_j=\beta_j,\Lambda^n_j=\nu_j,\forall j\in\{0,\ldots,{\lfloor pnu \rfloor}\};
\ell^n_{x^n_i}=l_i,\forall i\in\{0,\ldots,m\}\Big)>0.$$
Recalling the event considered in (\ref{ap7}), we also assume that
\begin{eqnarray*}
&& k\geq \eta\,n^{1/2},\\
&& l_i\geq \frac{1}{2}\eta' n^{1/4}\hbox{ \ if \ }i\in\{m-\lfloor\frac{1}{2}\eta n^{1/2}\rfloor,
\ldots,m\}.
\end{eqnarray*}
(To justify this last assumption, note that the label of a child of $\wt w^n_i$
can differ by at most $p$ from the label of $\wt w^n_i$, for any $i\in\{1,\ldots,N\}$.)
Consider then the conditional probability
$$P^{*,n}=P\Big(\cdot\,\Big|\, N=k,M=m;C^n_j=\beta_j,\Lambda^n_j=\nu_j,\forall j\in\{0,\ldots,{\lfloor pnu \rfloor}\};
\ell^n_{x^n_i}=l_i,\forall i\in\{0,\ldots,m\}\Big).$$
In order to get (\ref{ap6}), we need to verify the bound
\beq
\label{ap8}
P^{*,n}(B_n(u,\al,\delta))\leq K\,\al\,,
\eeq
where
$$B_n(u,\al,\delta)=\Big\{\inf_{\lfloor pnu\rfloor\leq i\leq \lfloor pn(u+\delta)\rfloor} C^n_i
\geq C^n_{\lfloor pnu \rfloor} -\al n^{1/2}\Big\},$$
and the constant $K$ depends only on $U,u_0,\eta,\eta'$ and $\delta$.

From now on, we argue under the conditional probability $P^{*,n}$. For every
$i\in\{0,1,\ldots,m\}$, we define a labelled $p$-tree 
$(\tau_{n,(i)},(\ell^{n,(i)}_v)_{v\in\tau_{n,(i)}^\circ})$ by the following:
\begin{description}
\item{$\bullet$} if $x^n_i\notin\{w^n_0,w^n_1,\ldots,w^n_k\}$, $\tau_{n,(i)}$
is just the subtree of descendants of $x^n_i$, meaning that
$$\tau_{n,(i)}=\{v\in\u:x^n_iv\in\tau_n\}$$
where $x^n_iv$ stands for the concatenation of the words $x^n_i$ and $v$;
\item{$\bullet$} if $x^n_i=w^n _j$ for some $j\in\{0,1,\ldots,k\}$, $\tau_{n,(i)}$
consists of $\varnothing$ and all $v\in\u$ such that $x^n_iv$ is a vertex of $\tau_n$
which has not yet appeared in $\{v^n_0,v^n_1,\ldots,v^n_{\lfloor pnu\rfloor}\}$;
\item{$\bullet$} $\ell^{n,(i)}_v=\ell^n_{x^n_iv}$ for every $i\in\{0,1,\ldots,m\}$
and $v\in\tau_{n,(i)}^\circ$.
\end{description}

We need to describe the distribution of the collection
$(\tau_{n,(i)},(\ell^{n,(i)}_v)_{v\in\tau_{n,(i)}^\circ})$, $0\leq i\leq m$
under $P^{*,n}$. To this end, introduce another collection
$$(T_i,(L^i_v)_{v\in T_i^\circ}),\quad 0\leq i\leq m$$
of independent random labelled $p$-trees defined under a
probability measure $\PP$ and such that, for $i\in\{0,1,\ldots,m\}$,
\begin{description}
\item{$\bullet$} $T_i$ is distributed as the multitype Galton-Watson tree $\tau_*$
considered in the proof of the lemma above;
\item{$\bullet$} given $T_i$, the labels $L^i_v,v\in T_i^\circ$, are distributed uniformly at random among 
admissible choices making $(T_i,(L^i_v)_{v\in T_i^\circ})$ a $p$-labelled tree with
$L^i_\varnothing =l_i$.
\end{description}
Consider also the two events
\begin{eqnarray*}
&&A_1=\{L^i_v>0\,,\hbox{ for every }i\in\{0,1,\ldots,m\}\hbox{ and }v\in T_i^\circ\},\\
&&A_2=\{p\sum_{i=0}^m \# T_i^\bullet + m =pn-\lfloor pnu \rfloor\}.
\end{eqnarray*}
Then the distribution of $(\tau_{n,(i)},(\ell^{n,(i)}_v)_{v\in\tau_{n,(i)}^\circ})_{0\leq i\leq m}$
under $P^{*,n}$ coincides with the distribution 
of $(T_i,(L^i_v)_{v\in T_i^\circ})_{0\leq i\leq m}$ under $\PP(\cdot\mid A_1\cap A_2)$.
This follows from the fact that, given the information provided by the 
conditioning event in $P^{*,n}$, the distribution of 
$(\tau_n,(\ell^n_v)_{v\in\tau_n})$ is uniform over the set of all
$p$-mobiles that are compatible with this information. 

Let us set $m'=m-\lfloor\frac{1}{2}\eta n^{1/2}\rfloor$ to simplify notation. Consider the
event
$$A'_1=\{ L^i_v>0\,,\hbox{ for every }i\in\{0,1,\ldots,m'\} \hbox{ and } v\in T^\circ_i\}.$$

\noindent{\bf Claim.} {\it There exists a constant $\ov K$, which depends only
on $\eta'$, such that for all $n$ sufficiently
large, $\PP(A'_1\cap A_2)\leq \ov K \,\PP(A_1\cap A_2)$.}

\smallskip
Let us verify the claim. To simplify notation, set $L(T_i)=\min\{L^i_v:v\in T_i\}$,
for every $i\in\{0,1,\ldots,m\}$. Then it suffices to prove that, for every
choice of the integers $\sigma_0,\sigma_1,\ldots,\sigma_m\geq 0$ such that
$p\sum_{i=0}^m \sigma_i+m\leq pn$,
$$\PP\Big(L(T_i)>0,\,\forall i\in\{m'+1,\ldots,m\}\,\Big|\, \# T_i^\bullet =\sigma_i,\,\forall i\in\{0,1,\ldots,m\}\Big)
\geq (\ov K)^{-1}.$$
with a constant $\ov K$ depending only on $\eta'$. By independence,
this conditional probability is equal to
$$\prod_{i=m'+1}^m \PP(L(T_i)>0\mid \# T_i^\bullet =\sigma_i)
= \prod_{i=m'+1}^m \Big(1-\Pi^p_{\sigma_i}\Big(\min_{v\in\tau^\circ} \ell_v\leq -l_i\Big)\Big),$$
with the notation of the lemma above. Now recall that
$l_i\geq \frac{1}{2}\eta' n^{1/4}$ for $i\in\{m'+1,\ldots,m\}$. 
Using the second assertion of the lemma, we get for $i\in\{m'+1,\ldots,m\}$,
$$\Pi^p_{\sigma_i}\Big(\min_{v\in\tau^\circ} \ell_v\leq -l_i\Big)
\leq 1-\exp(-\wt K \,\frac{\sigma_i}{l_i^4}),$$
where $\wt K$ only depends on $\eta'$. Hence,
\begin{eqnarray*}
\prod_{i=m'+1}^m \Big(1- \Pi^p_{\sigma_i}\Big(\min_{v\in\tau^\circ} \ell_v\leq -l_i\Big)\Big)
&\geq&\prod_{i=m'+1}^m \exp(-\wt K \,\frac{\sigma_i}{l_i^4})\\
&\geq&\exp\Big(-\frac{16\wt K}{\eta'^4n}\,\sum_{i=m'+1}^m \sigma_i\Big)\\
&\geq&\exp\Big(-\frac{16\wt K}{\eta'^4}\Big)
\end{eqnarray*}
since $\sum_{i=0}^m \sigma_i\leq n$. The claim now follows.

From the claim and the obvious property $A_1\subset A'_1$, we get
\beq
\label{ap9}
\PP(\cdot\mid A_1\cap A_2)\leq \ov K\,\PP(\cdot\mid A'_1\cap A_2).
\eeq

Under the probability measure $P^{*,n}$, if the event $B_n(u,\al,\delta)$ holds, then
we must have
$$\sum_{i=m-p(\lfloor \al n^{1/2}\rfloor+1)+1}^m (1+p\#\tau_{n,(i)}^\bullet )
\geq \lfloor pn(u+\delta)\rfloor -\lfloor pnu \rfloor$$
because, from the properties of the contour sequence and
our definitions, the ancestor of $v^n_{\lfloor pnu\rfloor}$
at generation $2(C^n_{\lfloor pnu\rfloor}-\lfloor \al n^{1/2}\rfloor-1)$ is visited by the
contour sequence after time $\lfloor pnu\rfloor$ and before time
$$\lfloor pnu\rfloor+\sum_{i=m-p(\lfloor \al n^{1/2}\rfloor+1)+1}^m 
(1+p\#\tau_{n,(i)}^\bullet ).$$
From this observation and (\ref{ap9}), we get the bound
\beq
\label{ap9bis}
P^{*,n}(B_n(u,\al,\delta))\leq \ov K\, 
\PP\Big(\sum_{i=m-p(\lfloor \al n^{1/2}\rfloor+1)+1}^m (1+p\# T_i^\bullet )
\geq \lfloor pn(u+\delta)\rfloor -\lfloor pnu \rfloor\,\Big|\, A'_1\cap A_2\Big).
\eeq
Notice that the event $A'_1$ does not involve the labelled trees
$(T_i,(L^i_v)_{v\in T_i^\circ})$ for $m'<i\leq m$. Hence the distribution
of the collection $(\# T_i^\bullet,m'<i\leq m)$ under $\PP(\cdot\mid A'_1\cap A_2)$
is exchangeable. 

On the other hand, put $q=\lfloor \frac{\eta}{2p\al}\rfloor -1\geq 1$ (recall that we assumed
$\al\leq \eta/(4p)$). Then $qp(\lfloor \al n^{1/2}\rfloor +1)<m-m'$ when $n$ is large. On the event
$A_2$, we have the trivial estimate
\begin{eqnarray*}
&&\#\Big\{j\in\{1,\ldots,q\}:\sum_{m-jp(\lfloor \al n^{1/2}\rfloor+1)<i\leq m
-(j-1)p(\lfloor \al n^{1/2}\rfloor+1)}
(1+p\# T_i^\bullet)
\geq \lfloor pn(u+\delta)\rfloor -\lfloor pnu\rfloor\Big\}\\
&&\qquad\leq \frac{pn-\lfloor pnu \rfloor}{ \lfloor pn(u+\delta)\rfloor -\lfloor pnu\rfloor}
\leq \frac{2(1-u)}{\delta}
\end{eqnarray*}
when $n$ is large.
By exchangeability, the expected value of the left-hand side under 
$\PP(\cdot\mid A'_1\cap A_2)$ equals
$$q\,\PP\Big(\sum_{i=m-p(\lfloor \al n^{1/2}\rfloor+1)+1}^m (1+p\# T_i^\bullet )
\geq \lfloor pn(u+\delta)\rfloor -\lfloor pnu \rfloor\,\Big|\, A'_1\cap A_2\Big)$$
and, using (\ref{ap9bis}), we finally obtain that 
$$P^{*,n}(B_n(u,\al,\delta))\leq
\ov K\,\frac{2(1-u)}{\delta}\,\Big(\lfloor \frac{\eta}{2p\al}\rfloor -1\Big)^{-1}.$$
The estimate (\ref{ap8}) follows, and this completes the proof. \cq

\medskip
\noi{\bf Proof of Lemma \ref{maxitech}:} As we already observed after the 
statement of Theorem \ref{mainIM}, we may choose
the uniformly distributed $p$-mobiles $(\tau_n,(\ell^n_v)_{v\in\tau^\circ_n})$
in such a way  that the convergence
\beq
\label{ap10}
\left(\lambda_p\,n^{-1/2}\,C^n_{\lfloor pnt\rfloor},
\kappa_p\,n^{-1/4} \Lambda^n_{\lfloor pnt\rfloor}\right)_{0\leq t\leq 1}\build{\la}_{}^{\rm a.s.} \left(\ov\eg_t,\ov Z_t
\right)_{0\leq t\leq 1}
\eeq
holds as $n\to\infty$ and not only along a subsequence $(n_k)$ as in (\ref{basicconv}).

Let $\beta>0$ and let $F_n(\beta)$ denote the event on which there 
exist  three integers $k^n_1\leq k^n_2\leq k^n_3$
in $[0,pn]$ such that $v^n_{k^n_1}=v^n_{k^n_2}=v^n_{k^n_3}$
and $k^n_2-k^n_1\geq \beta n$, $k^n_3-k^n_2\geq \beta n$. Then,
it is easy to see that
\beq
\label{ap11}
\lim_{\al \to 0} \Big(\limsup_{n\to\infty} P\left(F_n(\beta)\backslash \Delta_n(\al)\right)\Big)=0.
\eeq
Indeed, on the event
$$\limsup_{n\to\infty}\left(F_n(\beta)\backslash \Delta_n(\al)\right)$$
we can use (\ref{ap10}) and a compactness argument to find
three times $s_1\leq s_2\leq s_3$ in $[0,1]$ such that $s_1\sim s_2\sim s_3$, 
$s_2-s_1\geq \beta/p$, $s_3-s_2\geq \beta/p$ and furthermore
$$\min_{r\in[s_1,s_2]}\ov Z_r\geq \ov Z_{s_1}-\kappa_p\al
\hbox{ \ or \ }
\min_{r\in[s_2,s_3]} \ov Z_r\geq \ov Z_{s_1}-\kappa_p\al.$$
Lemma \ref{equirel} now shows that the probability of the existence of such
a triplet $(s_1,s_2,s_3)$ tends to $0$ as $\al\to 0$, which implies (\ref{ap11}).

Thanks to (\ref{ap11}), Lemma \ref{maxitech} will follow if we can verify that
\beq
\label{ap12}
\lim_{\beta\to 0} \Big(\liminf_{n\to\infty} P(F_n(\beta))\Big)=1\,.
\eeq
In the case $p=2$, the proof of (\ref{ap12}) is easy from (\ref{ap10}).
Indeed (\ref{ap10}) implies the existence of vertices of $\tau_n$
whose removal produces three components containing each
a number of vertices of order $n$. If $p=2$, these vertices must
be white vertices, but in the general case, they could possibly be black vertices,
and so we need a different argument.

A proof of (\ref{ap12}) can be given along the lines of the proof 
of Lemma \ref{excursion}. The arguments are relatively straightforward
but somewhat tedious, and for this reason we will leave certain details
to the reader.
We use the notation introduced 
in the proof of Lemma \ref{excursion}, with the particular choice $u=1/2$. In particular, we
let $N=C^n_{\lfloor pn/2\rfloor}$, $M\in\{N,\ldots,(p-1)N\}$, the random vertices
$x^n_0,\ldots,x^n_M$
and the labelled trees $(\tau_{n,(i)},(\ell^{n,(i)}_v)_{v\in \tau^\circ_{n,(i)}})$
for $0\leq i\leq M$ be defined as in the previous proof,
with $u=1/2$. For $\eta,\eta'>0$,
we also set 
$$M'=(M-2\lfloor \eta n^{1/2}\rfloor)^+\ ,\ M''=(M-\lfloor \eta n^{1/2}\rfloor)^+.$$

Let $\ve >0$. Thanks to the convergence (\ref{ap10}), for any fixed choice of $\eta$ and $\eta'$,
we can choose a constant $\beta_1=\beta_1(\eta,\eta')>0$ small enough so that
the probability that  the contour
sequence of $\tau^\circ_n$ 
visits one of the vertices $x^n_{M'},x^n_{M'+1},\ldots, x^n_{M''}$ between times
$\lfloor pn/2\rfloor$ and $\lfloor pn/2\rfloor + \beta_1n$ is bounded above by $\ve$ when 
$n$ is large.

We then denote by $A^n_{\eta,\eta'}$ the event
\begin{eqnarray*}A^n_{\eta,\eta'}&=&\Big\{N>4\eta n^{1/2};\;
W^n_{\lfloor pn/2\rfloor}(j) > \eta'n^{1/4}+p,\forall j\in \{N-2\lfloor \eta n^{1/2}\rfloor,\ldots, N\};\\
&&\qquad\qquad
\sup_{v\in\tau^\circ_{n,(i)}} |\ell^{n,(i)}_v-\ell^{n,(i)}_\varnothing|\leq \eta' n^{1/4},
\forall i\in \{M',\ldots,M''\}\Big\}.
\end{eqnarray*}
We fix $\eta$ and $\eta'$ so that
$P(A^n_{\eta,\eta'})>1-\ve$ for every $n$ large enough. 

By
making the distribution of the collection $(\tau_{n,(i)},(\ell^{n,(i)}_v)_{v\in \tau^\circ_{n,(i)}})_{
0\leq i\leq M}$ explicit as in the preceding proof, one verifies that the distribution
of the vector
$$(\# \tau^\circ_{n,(M'+i)})_{0\leq i\leq \lfloor \eta n^{1/2}\rfloor}$$
under the conditioned measure $P(\cdot\!\mid \! A^n_{\eta,\eta'})$
is exchangeable. Notice that the proportion of indices 
$i\in\{M',\ldots, M''\}$ such that $x^n_i$ is an ancestor of 
$v^n_{\lfloor pn/2\rfloor}$ is at least $1/p$ when $n$ is large.

Furthermore, from the convergence (\ref{ap10}) and properties
of the limiting random function $\ov\eg$ (recall from subsection 2.4 that 
$\ov\eg$ is defined by a simple transformation 
of the Brownian excursion $\eg$), one gets that, for every integer $K\geq 1$,
we can choose $\beta>0$ sufficiently small so that
$$\liminf_{n\to \infty} P\Big(
A^n_{\eta,\eta'}\cap\Big\{\#\{i\in\{0,1,\ldots,\lfloor \eta n^{1/2}\rfloor\}:
\# \tau^\circ_{n,(M'+i)}> \beta n\} \geq K\Big\}\Big)\geq 1-2\ve.$$
Using an exchangeability argument, we can then find a constant $\beta_0\in]0,\beta_1(\eta,\eta')]$ such that the
following holds:  With probability
greater than $1-3\ve$ when $n$ is large, the event $A^n_{\eta,\eta'}$ holds and there exists 
$i_0\in\{M',\ldots,M''\}$ such that $x^n_{i_0}$ is an ancestor of 
$v^n_{\lfloor pn/2\rfloor}$, and $\# \tau^\circ_{n,(i_0)}> \beta_0 n$. 

On the event where $i_0$ is defined, let $k^n_1$ stand for the time of the first occurence of $x^n_{i_0}$
in the contour sequence of $\tau^\circ_n$, let $k^n_2$ be the time of the first occurence 
of $x^n_{i_0}$ after time
$\lfloor pn/2\rfloor$ and let $k^n_3$ be the time of the last occurence
of $x^n_{i_0}$. Then the triplet $(k^n_1,k^n_2,k^n_3)$
will be defined and will satisfy the property of the definition of $F_n(\beta_0)$,
with a probability greater than $1-4\ve$ when $n$ is large. Our claim 
(\ref{ap12}) now follows. \cq

\medskip
D\'epartement de math\'ematiques, Universit\'e
Paris-Sud, 91405 ORSAY Cedex, FRANCE

 e-mail:
jean-francois.legall@math.u-psud.fr

\end{document}